\documentclass{article}
\usepackage{arxiv}
\usepackage[document]{ragged2e}
\usepackage[numbers]{natbib}
\usepackage{amsmath}
\usepackage{amsfonts}
\usepackage{amssymb}
\usepackage{bm}
\usepackage{enumitem}
\usepackage{booktabs}
\usepackage{siunitx}
\usepackage{graphicx}
\usepackage{hyperref}
\usepackage{fancyhdr}
\usepackage{afterpage}
\usepackage{comment}
\usepackage{multirow}
\usepackage{url}
\usepackage[aboveskip=2pt, labelfont=bf,labelsep=period,justification=raggedright,singlelinecheck=off, font={normalsize,it}]{caption}

\newcommand{\R}{\mathbb{R}}
\DeclareSIUnit\dyne{dyn}
\newtheorem{remark}{Remark}[section]
\graphicspath{ {img/} }

\newcommand\mytilde[1]{\stackrel{\sim}{\smash{#1}\rule{0pt}{1.05ex}}}

\newcommand{\widebar}[1]{\mkern 1.5mu\overline{\mkern-1.5mu#1\mkern-1.5mu}\mkern 1.5mu}
\renewcommand{\vec}[1]{\underline{#1}}

\title{Model order reduction of h{\ae}modynamics by space--time reduced basis and reduced fluid--structure interaction}

\author{ \href{https://orcid.org/0000-0002-2784-9261}{\includegraphics[scale=0.06]{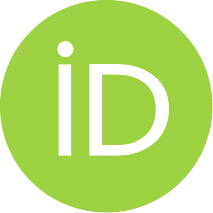} \hspace{1mm} Riccardo Tenderini}\\
	Institute of Mathematics\\
	Ècole Polytechnique Fédérale de Lausanne (EPFL)\\
	CH--1015 Lausanne, Switzerland \\
	\And
	\href{https://orcid.org/0000-0002-2832-6630}{\includegraphics[scale=0.06]{img/orcid.pdf} \hspace{1mm} Simone Deparis} \\
	Institute of Mathematics\\
	Ècole Polytechnique Fédérale de Lausanne (EPFL)\\
	CH--1015 Lausanne, Switzerland \\
}

\date{}

\renewcommand{\headeright}{Model order reduction of h{\ae}modynamics by ST--RB and reduced FSI}
\renewcommand{\undertitle}{}

\afterpage{\cfoot{\thepage}}

\numberwithin{equation}{section}

\begin{document}
	
\maketitle

\vspace{-0.5cm}

\begin{abstract}
	\justifying
	In this work, we apply the space--time Galerkin reduced basis (ST--GRB) method to a reduced fluid--structure interaction model, for the numerical simulation of h{\ae}modynamics in arteries. In essence, ST--GRB extends the classical reduced basis (RB) method, exploiting a data--driven low--dimensional linear encoding of the temporal dynamics to further cut the computational costs. The current investigation brings forth two key enhancements, compared to previous works on the topic. On the one side, we model blood flow through the Navier--Stokes equations, hence accounting for convection. In this regard, we implement a hyper--reduction scheme, based on approximate space--time reduced affine decompositions, to deal with nonlinearities effectively. On the other side, we move beyond the constraint of modelling blood vessels as rigid structures, acknowledging the importance of elasticity for the accurate simulation of complex blood flow patterns. To limit computational complexity, we adopt the Coupled Momentum model, incorporating the effect of wall compliance in the fluid's equations through a generalized Robin boundary condition. In particular, we propose an efficient strategy for handling the spatio--temporal projection of the structural displacement, which ultimately configures as a by--product. The performances of ST--GRB are assessed in three different numerical experiments. The results confirm that the proposed approach can outperform the classical RB method, yielding precise approximations of high--fidelity solutions at more convenient costs. However, the computational gains of ST--GRB vanish if the number of retained temporal modes is too large, which occurs either when complex dynamics arise or if very precise solutions are sought.   
\end{abstract}  

\keywords{Reduced order modeling, Space--time reduced basis, Navier--Stokes equations, H{\ae}modynamics,  Fluid--structure interaction}

\justifying

\section{Introduction}
\label{sec:introduction}
The use of physics--based computational models for the cardiovascular system has significantly grown in recent years, enabling the accurate simulation of h{\ae}modynamics, tissue mechanics, and of the cardiovascular system in general. As a matter of fact, these predictive models, originally designed for surgical planning, now serve purposes well beyond their initial application, thanks to their ability to combine clinical data with physical knowledge of the processes underlying the cardiac function \cite{schwarz2023beyond, fumagalli2024role}. High--fidelity personalized numerical simulations of h{\ae}modynamics can yield exceptionally accurate descriptions of the blood flow conditions by leveraging state--of--the--art computational fluid dynamics (CFD) techniques. Nevertheless, they require important computational resources and entail large computing times, which hinders their effective use in the daily clinical practice. As a way to mitigate this issue, Reduced Order Models (ROMs) allow to conduct simulations on standard consumer devices with computational times of the order of minutes or even seconds, hence representing a cost--effective alternative to the classical CFD approaches. Such substantial performance improvements come at the cost of accuracy losses, that remain nonetheless within the range of acceptability for many applications. Therefore, ROMs can be particularly useful in many--query scenarios, namely in cases where multiple simulations have to be carried out, considering different values of some parameters, related to the system properties. Uncertainty Quantification (UQ) is perhaps the most relevant example in this sense~\cite{fleeter2020multilevel}. \\

In h{\ae}modynamics, ROMs can be roughly subdivided into three categories: physics--based --- such as OD, 1D models \cite{shi2011review} --- projection--based --- such as the Reduced Basis (RB) method \cite{quarteroni2015reduced} --- and data--driven \cite{arzani2021data}. Hybrid approaches are also possible; for instance, several methods combining data availability with knowledge of the underlying physical laws --- such as Physics--Informed Neural Networks \cite{raissi2019physics,arzani2021data,moser2023modeling} --- have been introduced over the last decade. We refer to \cite{schwarz2023beyond, lombardi2023reduced} for comprehensive literature reviews of the topic. 
In this work, we focus on projection--based ROMs, and more precisely on the RB method. The latter stems from the orthogonal projection (in space) of the model equations onto linear subspaces, generated either with greedy algorithms~\cite{binev2011convergence, hesthaven2014efficient} or through Proper Orthogonal Decomposition (POD) from a database of high--fidelity solutions~\cite{kunisch2002galerkin}. Despite having been originally conceived for steady problems, the RB method can be seamlessly extended to time--dependent ones. This leads to the time integration of low--dimensional systems of ODEs, which is usually performed with the same scheme and timestep size used for the full--order model~(FOM)~\cite{hesthaven2022reduced}. In many cases, the dimensionality reduction in space introduced by the RB method is enough to drastically cut the computational costs. Indeed, the FOM spatial complexity typically exceeds the temporal one by a few orders of magnitude. However, in problems where either the simulation interval should be very large or the timestep size should be very small in order to capture some relevant physical behaviours, temporal complexity does represent a bottleneck and valuable computational gains may be difficult to realize with classical approaches. \\

To overcome the temporal complexity bottleneck, in this work we consider Space--Time Reduced Basis (ST--RB) methods \cite{urban2014improved, yano2014space, choi2019space, shimizu2021windowed, choi2021space, kim2021efficient, tenderini2022pde, tenderini2024space, mueller2024model}. In~\cite{tenderini2024space}, we already investigated the use of ST--RB methods based on spatio--temporal (Petrov--)Galerkin projections for h{\ae}modynamics. In particular, we proposed two stabilization procedures, that guarantee the well--posedness of the resulting low--dimensional problem. Both approaches attained remarkable computational gains compared to the baseline RB method and accuracies matching the theoretical expectations. However, these results were obtained in a simplified setting, where we neglected the presence of the nonlinear convective term in the momentum conservation equation (Stokes equations) and modelled arterial walls as rigid structures, via the imposition of no--slip boundary conditions. Solely focusing on the case of a Galerkin projection in space--time, the current investigation brings forth two key enhancements.
On the one hand, we model blood flow by the unsteady incompressible Navier--Stokes equations, recognizing the primary importance of convection for the precise simulation of h{\ae}modynamics. To effectively handle nonlinearities, we extend to space--time the hyper--reduction technique introduced in \cite{pegolotti2021model}, based on approximate affine decompositions which directly stem from the weak form of the convective term. Albeit being fairly standard, this choice is of crucial importance in our setting, since it allows to achieve non--intrusivity in the online phase, which is completely independent from the full--order code base. 
On the other hand, we move beyond the no--slip constraint, taking wall compliance into account for the accurate approximation of complex blood flow patterns \cite{nobile2001numerical, moireau2012external, pons2020fluid}. In principle, this would imply solving a fluid--structure interaction (FSI) problem, where the Navier--Stokes equations in a moving domain are coupled with elasto--mechanical models that describe structural displacements, by means of an Arbitrary Lagrangian--Eulerian (ALE) approach. Even if the development of efficient parallel algorithms \cite{crosetto2011parallel, nobile2013time} and of \emph{ad hoc} preconditioning strategies \cite{badia2009robin, deparis2016facsi} led to significant efficiency improvements, the computational cost of FSI simulations remains nonetheless high. In the h{\ae}modynamics context, this poses limitations to their application in the standard clinical practice and in particular to their use for real--time diagnoses. To mitigate this issue, we consider a reduced FSI model, called Coupled Momentum, originally proposed in~\cite{figueroa2006coupled} and successively further investigated e.g. in \cite{colciago2014comparisons, colciago2018reduced}. Under the hypothesis of small deformations --- reasonable in the majority of large blood vessels --- the coupled momentum model assumes that the fluid domain is fixed and incorporates the effect of wall compliance on the blood flow within the Navier--Stokes equations, by means of a generalized Robin boundary condition.
Furthermore, compared to \cite{tenderini2024space}, we impose absorbing resistance boundary conditions at the outlets, enforcing a linear relationship between pressure and volumetric flow rate. An appropriate choice of the outlet boundary conditions, able to adequately model the effect of the downstream vasculature, is important for performing accurate personalized simulations, and it becomes fundamental when accounting for compliant vessel walls. In particular, the choice of suitable and well--calibrated outlet boundary conditions prevents spurious reflecting pressure waves, thus both leading to more physiological results and enhancing the computational gains of model order reduction in space--time. 
Ultimately, to the best of our knowledge, this represents the first work in which projection--based model order reduction in space--time is combined with reduced FSI modelling for h{\ae}modynamics, enabling efficient, non--intrusive, and physics--based evaluation of physiological blood flow in compliant vessels under parametric variability. In particular, the linear compression of temporal dynamics onto suitably constructed subspaces is expected to yield superior runtime performance compared to sequential--in--time ROMs such as the canonical RB method, at the cost of marginal and \emph{a priori} quantifiable accuracy losses. As a result, this trade--off makes our approach particularly appealing for computationally demanding multi--query applications in cardiovascular modelling, characterized by non--negligible complexities in time. \\

The manuscript is structured as follows. In Section~\ref{sec:RFSI}, we introduce the unsteady incompressible Navier--Stokes equations and the coupled momentum model, and we discuss their weak formulation and high--fidelity discretization by the Finite Element (FE) method. In Section~\ref{sec:ST-RB for RFSI}, we describe the application of the ST--RB method based on a Galerkin projection to the problem at hand, detailing the assembling of the reduced convective term, the treatment of the generalized Robin boundary condition and the handling of non--zero initial conditions. In Section~\ref{sec: numerical results}, we present three different numerical tests (two idealized and one physiological) to empirically assess the method performances. Finally, in Section~\ref{sec:conclusions}, we summarize the main findings of the work and envision potential further developments.

\section{Model order reduction of FSI: the Coupled Momentum model} \label{sec:RFSI}

\subsection{Physical problem} 
\label{subsec: physical problem}
In sufficiently large vessels, blood can be approximated as a Newtonian fluid and its flow can be modelled by the incompressible Navier--Stokes equations \cite{formaggia2010cardiovascular}. Let $\Omega \subset \R^3$ be an open, simply connected and bounded domain and let $\partial\Omega$ be its boundary. Let $[0,T]$, with $T>0$, be the time interval of interest. Then, the unsteady incompressible Navier--Stokes equations read as follows:
\begin{equation} 
\label{eq: strong_form_NS}
\begin{cases}
\rho_f\vec{u}_t + \rho_f(\vec{u} \cdot \nabla)\vec{u} - \nabla\cdot (2\mu_f\nabla^s\vec{u}) + \nabla p = \vec{0}  & \text{in} \ \Omega \times (0,T]~, \\
\nabla \cdot \vec{u} = 0  & \text{in} \ \Omega \times (0,T]~.
\end{cases} 
\end{equation}
Here $\vec{u}: \Omega \times (0,T] \to \R^3$ and $p: \Omega \times (0,T] \to \R$ denote velocity and pressure, respectively; $\vec{u}_t$ is the partial derivative of the velocity with respect to time; $\rho_f, \mu_f \in \R^+$ are the blood density and viscosity, respectively; $\nabla^s\vec{u} = (\nabla\vec{u}+\nabla^T\vec{u}) \ / \ 2$ is the strain rate tensor. We remark that the first equation in Eq.\eqref{eq: strong_form_NS}, that enforces momentum conservation, is nonlinear, due to presence of the convective term $\vec{c}(\vec{u}) := \rho_f (\vec{u} \cdot \nabla)\vec{u}$. Eq.\eqref{eq: strong_form_NS} must be endowed with suitable initial and boundary conditions to guarantee well--posedness. Since the velocity partial derivative in time appears, the initial condition is
\begin{equation*}
\label{eq: initial conditions}
\vec{u} = \vec{u}_0 \quad \text{in} \quad \Omega \times \{0\}~,
\end{equation*}
where $\vec{u}_0: \Omega \to \R^3$ is a prescribed velocity field at $t=0$. For the purpose of defining the boundary conditions, let $\left\{\Gamma, \Gamma_{\operatorname{IN}}, \Gamma_{\operatorname{OUT}}\right\}$ be a partition of $\partial\Omega$, where $\Gamma$ denotes the vessel wall and $\Gamma_{\operatorname{IN}}$, $\Gamma_{\operatorname{OUT}}$ encode the ensemble of the inlet and outlet boundaries, respectively. In this work, we always prescribe velocity profiles at the inlets, whereas the nature of the outlet boundary conditions depends on the test case at hand. For this reason, we introduce $\Gamma_{\operatorname{D}} = \bigcup_{k=1}^{N_D} \Gamma_{\operatorname{D}}^k $ as the union of the $N_D$ inlet/outlet boundaries where velocity profiles are enforced (Dirichlet boundaries); analogously, we define $\Gamma_{\operatorname{N}} := \bigcup_{k=1}^{N_N} \Gamma_{\operatorname{N}}^k$ (Neumann boundaries) and $\Gamma_{\operatorname{R}} := \bigcup_{k=1}^{N_R} \Gamma_{\operatorname{R}}^k$ (resistance boundaries). Clearly, it holds that $\Gamma_{\operatorname{IN}} \cup \Gamma_{\operatorname{OUT}} = \Gamma_{\operatorname{D}} \cup \Gamma_{\operatorname{N}} \cup \Gamma_{\operatorname{R}}$.
Then, the inlet/outlet boundary conditions read as follows:
\begin{equation*}
\label{eq: IO boundary conditions}
\begin{cases}
\vec{u} = \vec{g}^k & \text{on} \ \Gamma_{\operatorname{D}}^k \times (0,T]~, \ \ k=1, \cdots, N_D~; \\
\sigma(\vec{u}, p) \ \vec{n} = \vec{0} & \text{on} \ \Gamma_{\operatorname{N}} \times (0,T]~; \\
\sigma(\vec{u}, p) \ \vec{n} + (R^{k'} Q^{k'}) \ \vec{n} = \vec{0}  & \text{on} \ \Gamma_{\operatorname{R}}^{k'} \times (0,T]~, \ \ k'=1, \cdots, N_R~; 
\end{cases}
\end{equation*}
where $\vec{g}^k: \Gamma_{\operatorname{D}}^k \times (0,T] \to \R^3$ is the Dirichlet boundary datum on $\Gamma_{\operatorname{D}}^k$, $\sigma(\vec{u}, p) = 2\mu_f\nabla^s\vec{u} - pI$ denotes the Cauchy stress tensor, $\vec{n}=\vec{n}(\vec{x})$ is the outward unit normal vector to $\partial\Omega$ at $\vec{x}$, and $Q^{k'} = \int_{\Gamma_{\operatorname{R}}^{k'}} \vec{u} \cdot \vec{n} \ $ is the volumetric flow rate at $\Gamma_{\operatorname{R}}^{k'}$. \\

To limit computational complexity, whilst taking into account the effect of the vessel wall compliance on blood flow, we employ the \emph{Coupled Momentum} model \cite{figueroa2006coupled}. To this aim, we consider a pseudo--compressible, linear elastic, isotropic, homogeneous material for the vessel structure and we suppose to work under a small deformations hypothesis. In this setting, by approximating the three--dimensional vessel structure as a thin elastic membrane, the coupled momentum model translates the wall mechanics into an \emph{ad hoc} generalized Robin boundary condition, enforced on $\Gamma$. Let $\rho_s \in \R^+$ be the structure density and let $\vec{d}: \Gamma \times (0,T] \to \R^3$ denote the membrane displacement. Following \cite{figueroa2009effect}, we consider a membrane model in the plane stress configuration with an augmented shear stress in the transverse directions. This leads to the following linear stress--strain relation:
\begin{equation}
\label{eq: stress-strain membrane}
\Pi_{\Gamma}(\vec{d}) = h_s \Big(\lambda_1 (\nabla_{\Gamma} \cdot \vec{d}) \ (I - \vec{n} \otimes \vec{n}) + 2\lambda_2 \nabla_{\Gamma}^s(\vec{d}) + 2\lambda_2(\gamma_{sh}-1)  \nabla_{\Gamma}^s(\vec{d}) \ \vec{n} \otimes \vec{n} \Big)~,
\end{equation}
where $\nabla_{\Gamma}: \R^3 \to \R^{3 \times 3}$ is the surface gradient operator on $\Gamma$, i.e. $\nabla_{\Gamma}(\vec{d}) := \nabla \vec{d} \ (I - \vec{n} \otimes \vec{n})$, and $\nabla_{\Gamma}^s(\vec{d})$ denotes its symmetric part; the symbol $\otimes$ denotes the tensor product and $I \in \R^{3 \times 3}$ is the identity matrix. The quantity $h_s \in \R^+$ represents the membrane thickness, $\gamma_{sh}=5/6$ is the transverse shear coefficient, while $\lambda_1, \lambda_2 \in \R$ are parameters whose expressions depend on the Young modulus $E$ and on the Poisson ratio $\nu$ of the material:
\begin{equation*}
\label{eq: lame parameters}
\lambda_1 = \frac{E \nu}{(1+\nu)(1-\nu)}~, \qquad \lambda_2 = \frac{E}{2(1+\nu)}~.
\end{equation*}
Furthermore, as done in \cite{malossi2012two}, we optionally take into account the elastic response of the surrounding tissues, to avoid non--physiological rigid vessel motions~\cite{liu2007surrounding}. We refer the reader to \cite{colciago2014comparisons} for further details on the coupled momentum model and on its relevance in the context of h{\ae}modynamics.

\begin{remark}
For simplicity, in the following we restrict to homogeneous initial conditions for both velocity and displacement. The full--order and reduced--order problem formulations with non--zero initial conditions are detailed in Appendix~\ref{app: incorporating inomogeneous initial conditions}.
\end{remark}

\subsection{Full--order algebraic formulation} \label{subsec: algebraic formulation}
We are interested in deriving a high--fidelity space--time algebraic formulation of the problem at hand. Concerning the spatial discretization, we employ the FE method and consider a stable couple $(\mathcal{V}_h, \mathcal{Q}_h)$ of FE spaces for velocity and pressure, built on a tetrahedralization of $\Omega$. We remark that $(\mathcal{V}_h, \mathcal{Q}_h)$ must satisfy the \emph{inf--sup} condition, in order to guarantee the well--posedness of the discrete problem~\cite{brezzi1974existence, boffi2013mixed}. For the temporal discretization, we partition the time interval $[0,T]$, by defining discrete time instants $\{t_n\}_{n=0}^{N^t}$ such that $t_n = n \Delta t$; $\Delta t \in \R^+$ is called the timestep size. To integrate the multidimensional ODE stemming from the application of the FE method, we rely on implicit BDF schemes. So, given a generic vector quantity $\phi=\phi(\cdot, t)$, we approximate its first--order partial derivative in time as
\begin{equation*}
\label{eq: BDF scheme}
\left.\frac{\partial\phi}{\partial t}\right\vert_{t=t_{n+1}} = \frac{1}{\beta \Delta t} \Bigg(\phi_{n+1} - \sum_{s=1}^S \alpha_s \phi_{n+1-s}\Bigg) \qquad \text{with} \quad \phi_m = \phi(\cdot, t_m)~,
\end{equation*} 
for suitably defined coefficients $\alpha_1, \dots, \alpha_S, \beta \in \R$. Additionally, we impose the inhomogeneous Dirichlet boundary conditions on $\Gamma_{\operatorname{D}}$ with Lagrange multipliers, to facilitate the extension of the proposed approach for coupling several subdomains \cite{pegolotti2021model}. To this aim, we introduce the (discrete) space of Lagrange multipliers $\mathcal{L}_h = \prod_{k=1}^{N_D} \mathcal{L}_h^k$. We refer to~\cite{pegolotti2021model} for further details on the choice of the space of Lagrange multipliers and on its discretization; we refer to~\cite{tenderini2024space} for a discussion on the well--posedness of the resulting problem, that exhibits a twofold saddle--point structure. Ultimately, the following finite--dimensional approximations of velocity, pressure, and Lagrange multipliers are considered:
\begin{equation*}
\label{eq: finite element approximation}
\vec{u}_h(\vec{x},t) = \sum_{i=1}^{N_u^s} \bm{u}_i(t) \ \vec{\varphi}_i^u(\vec{x}) \in \mathcal{V}_h~;
\qquad
p_h(\vec{x},t) = \sum_{i=1}^{N_p^s} \bm{p}_i(t) \ \varphi_i^p(\vec{x}) \in \mathcal{Q}_h~;
\qquad
\vec{\lambda}_h^k(\vec{x},t) = \sum_{i=1}^{N_\lambda^k} \bm{\lambda^k}_i(t) \ \vec{\eta}^k_i(\vec{x}) \in \mathcal{L}_h^k~;
\end{equation*}
where $N_u^s := \operatorname{dim}(\mathcal{V}_h)$, $N_p^s := \operatorname{dim}(\mathcal{Q}_h)$, $N_\lambda^k := \operatorname{dim}(\mathcal{L}_h^k)$, for $k = 1, \cdots, N_D$. Also, we define $N_\lambda := \prod_{k=1}^{N_d} N_\lambda^k$. The time--dependent vectors $\bm{u} \in \R^{N_u^s}$, $\bm{p} \in \R^{N_p^s}$, $\bm{\lambda^k} \in \R^{N_\lambda^k}$ collect the velocity, pressure, and Lagrange multipliers Degrees of Freedom (DOFs). \\

In view of the use of model order reduction techniques, we introduce a parametrization in the problem. Firstly, we parametrize the inflow/outflow rates at the Dirichlet boundaries. For every $k \in \{1, \cdots N_D\}$, we assume that the Dirichlet datum $\vec{g}^k$ on $\Gamma_{\operatorname{D}}^k$ can be factorized as
\begin{equation}
\label{eq: dirichlet datum}
\vec{g}^k(\vec{x}, t; \bm{\mu_f}) := \vec{g}_k^s(\vec{x}) \ g_k^t(t; \bm{\mu_f}) \qquad (\vec{x}, t) \in \Gamma_{\operatorname{D}}^k \times (0,T]~,
\end{equation}
where $\bm{\mu_f} \in \mathcal{P}_f \subset \R^{N_\mu^f}$ is a vector collecting parameters associated to the prescribed volumetric flow rates. Secondly, we introduce a source of parametrization that affects the vessel structure mechanics. To this aim, we define the parameter vector $\bm{\mu_m} := [h_s, \rho_s, E, \nu] \in \mathcal{P}_m \subset \R^{N_\mu^m}$, where $N_\mu^m = 4$; the entries of $\bm{\mu_m}$ have been defined in Section~\ref{subsec: physical problem}. We remark that, since the stress--strain relation reported in Eq.\eqref{eq: stress-strain membrane} is linear with respect to $\lambda_1(E,\nu)$, $\lambda_2(E, \nu)$, the considered parametrization is affine. The coefficients that multiply the affine components are stored in the vector $\bm{\tilde{\mu}_m} := [h_s \rho_s, h_s \lambda_1(E, \nu), 2h_s\lambda_2(E, \nu)] \in \R^{N_\mu^m-1}$. \\

We can now introduce the following FOM matrices and vectors:
\begin{equation*}
\label{eq: FOM matrices}
\begin{alignedat}{5}
&\bm{A} \in \R^{N_u^s \times N_u^s} :\quad &&\bm{A}_{ij} = 2\mu_f \int_{\Omega}\nabla^s(\vec{\varphi}_j^u):\nabla^s(\vec{\varphi}_i^u)~; 
\quad  
&&\bm{M} \in \R^{N_u^s \times N_u^s} :\quad  &&\bm{M}_{ij} = \rho_f \int_{\Omega}\vec{\varphi}_j^u \cdot \vec{\varphi}_i^u~; 
\\
&\bm{B} \in \R^{N_p^s \times N_u^s} :\quad &&\bm{B}_{ij} = -\int_{\Omega}\varphi_i^p \ \nabla\cdot\vec{\varphi}_j^u~;
\quad 
&&\bm{L}^k \in \R^{N_\lambda^k \times N_u^s} :\quad &&\bm{L}^k_{ij} = \int_{\Gamma_{\operatorname{D}}} \vec{\eta}_i^k \cdot \vec{\varphi}_j^u~;  
\\
&\bm{M_s} \in \R^{N_u^s \times N_u^s} :\quad &&(\bm{M_s})_{ij} = \int_{\Gamma} \vec{\varphi}_j^u \cdot \vec{\varphi}_i^u~; 
\quad  
&&\bm{A_s}^q \in \R^{N_u^s \times N_u^s} :\quad &&(\bm{A_s}^q)_{ij} = a_{mem}^q(\vec{\varphi}_j^u; \vec{\varphi}_i^u)~;
\\
&\bm{c}(\vec{u}_h) \in \R^{N_u^s} :\quad &&\bm{c}(\vec{u}_h)_i = \int_\Omega \ \rho_f \left( \left(\vec{u}_h \cdot \nabla \right) \vec{\varphi}_j^u \right) \cdot \vec{\varphi}_i^u~; 
\quad
&&\bm{g}_k^s \in \R^{N^k_\lambda} :\quad &&(\bm{g}_k^s)_i = \int_{\Gamma_{\operatorname{D}}^k} \vec{g}_k^s \cdot \vec{\eta}_i^k~;  
\\
&\bm{q}^{k'} \in \R^{N_u^s} :\quad &&(\bm{q}^{k'})_i =  \int_{\Gamma_{\operatorname{R}}^{k'}} \ \vec{\varphi}_i^u \cdot \vec{n}~;
\quad
&&\bm{R}^{k'} \in \R^{N_u^s \times N_u^s} :\quad &&(\bm{R}^{k'})_{ij} = R^{k'} (\bm{q}^{k'})_i \  (\bm{q}^{k'})_j~;
\end{alignedat}
\end{equation*}
where $k = 1, \cdots, N_D$, $k' = 1, \cdots, N_R$, and $q = 1,2$. The function $\vec{g}_h^k \in \mathcal{V}_h$ is the finite element approximation of $\vec{g}_k^s \in L^2(\Gamma_{\operatorname{D}}^k)$ introduced in Eq.\eqref{eq: dirichlet datum}, while the following bilinear forms are considered for the assembling of the matrices $\{\bm{A_s}^q\}_{q=1}^2$:
\begin{equation*}
\label{eq: affine components membrane}
\begin{alignedat}{2}
a_{mem}^1\left(\vec{\varphi}_j^u; \vec{\varphi}_i^u \right) &= \beta \Delta t \int_{\Gamma} \left(\nabla_{\Gamma} \cdot \vec{\varphi}_j^u\right) \ (I - \vec{n} \otimes \vec{n}) : \nabla_{\Gamma} \vec{\varphi}_i^u~; \\
a_{mem}^2\left(\vec{\varphi}_j^u; \vec{\varphi}_i^u  \right) &= \beta \Delta t \int_{\Gamma} \nabla_{\Gamma}^s \vec{\varphi}_j^u : \nabla_{\Gamma} \vec{\varphi}_i^u
\ + \ (\gamma_{sh} - 1) \left(~\nabla_{\Gamma}^s \vec{\varphi}_j^u \vec{n} \otimes \vec{n}\right) : \nabla_{\Gamma} \vec{\varphi}_i^u~.
\end{alignedat}
\end{equation*}
To ease the notation, we aggregate the information relative to the Lagrange multipliers by introducing 
\begin{equation*}
\label{eq: FOM matrices Lagrange multipliers}
\bm{L} = \left[ \left(\bm{L}^1\right)^T \lvert \ \cdots \ \rvert \left(\bm{L}^{N_D}\right)^T \right]^T \ \in \R^{N_\lambda \times N_u^s}~; \qquad
\bm{g}(t; \bm{\mu_f}) = \left[ \big(\bm{g}^s_1\big)^T g^t_1(t; \bm{\mu_f}) \ \lvert \ \cdots \ \rvert \left(\bm{g}^s_{N_D}\right)^T g^t_{N_D}(t; \bm{\mu_f}) \right]^T \ \in \R^{N_\lambda}.
\end{equation*} 
and we define the following cumulative resistance matrix:
\begin{equation*}
\label{eq: FOM matrices resistance BCs}
\bm{R} = \sum_{k'=1}^{N_R} \bm{R}^{k'} \ \in \R^{N_u^s \times N_u^s}~.
\end{equation*} 

Let $\bm{w}_n := \left[\bm{u}_n, \bm{p}_n, \bm{\lambda}_n\right]^T \in \R^{N^s}$ be the discrete solution at time $t_n$, with $N^s := N_u^s + N_p^s + N_\lambda$; let $\bm{d}_n \in \R^{N_u^s}$ be the discrete representation of the displacement at time $t_n$, expressed with respect to trace of the velocity basis functions on $\Gamma$. Furthermore, let us define $\bm{c}(\bm{u}_n) := \bm{c}\big(\sum_k \bm{u}_{n,k} \vec{\varphi}_k^u\big) \in \R^{N_u^s}$. Then, given $\bm{w}_{n+1-s}$ for $s=1,\cdots,S$, the solution $\bm{w}_{n+1}$ at time $t_{n+1}$ is such that:
\begin{equation}
\label{eq: residual BDF}
\bm{r}(\bm{w}_{n+1}) := \bm{H} \bm{w}_{n+1} - \sum_{s=1}^S \alpha_s \ \bm{H} \bm{w}_{n+1-s} - \beta \Delta t \ \bm{F}(t_{n+1}, \bm{w}_{n+1}, \bm{d}_{n+1}) = \bm{0}~,
\end{equation}
where $\bm{H} \bm{w}_n = \left[ \big(\left(\bm{M} + (\bm{\tilde{\mu}_m})_1 \ \bm{M_s}\right)\bm{u}_n\big)^T, \ {\bm{0}_{N_p^s}}^T, \ {\bm{0}_{N_\lambda}}^T \right]^T$ and
\begin{equation*}
\label{eq: FOM right-hand side}
\bm{F}(t_n, \bm{w}_n, \bm{d}_n) = 
\begin{bmatrix}
\bm{0} \\ \bm{0} \\ \bm{g}(t_n; \bm{\mu_f})
\end{bmatrix}
- 
\begin{bmatrix}
\bm{A} + \bm{R} & \bm{B}^T & \bm{L}^T \\
\bm{B} & \bm{0}   & \bm{0} \\
\bm{L} & \bm{0}   & \bm{0}
\end{bmatrix}
\begin{bmatrix}
\bm{u}_n \\ \bm{p}_n \\ \bm{\lambda}_n
\end{bmatrix} 
-
\begin{bmatrix}
\bm{c}\left(\bm{u}_n\right) \\ \bm{0} \\ \bm{0} \\
\end{bmatrix}
-
\begin{bmatrix}
\left(\bm{A_s}(\bm{\tilde{\mu}_m}) + c_s \ \bm{M_s} \right) \ \bm{d}_n \\ \bm{0} \\ \bm{0} \\
\end{bmatrix}~,
\end{equation*}
with $\bm{A_s}(\bm{\tilde{\mu}_m}) := \sum_{q=1}^2 (\bm{\tilde{\mu}_m})_{q+1} \ \bm{A_s}^q$. Here, $c_s \in \R^+$ is a parameter that encodes the elastic response of the surrounding tissues. \\

Eq.\eqref{eq: residual BDF} can be conveniently rewritten as a monolithic $N^{st}$--dimensional nonlinear system --- with $N^{st} := N^sN^t$ --- as follows:
\begin{equation}
\label{eq: monolitic_FOM_system}
\begin{bmatrix} 
\bm{A}^{st}_1 + (\bm{\tilde{\mu}_m})_1 \bm{M_s}^{st} & \bm{A}^{st}_2 & \bm{A}^{st}_3 \\ 
\bm{A}^{st}_4 & \bm{0} & \bm{0}  \\ 
\bm{A}^{st}_7 & \bm{0} & \bm{0}  
\end{bmatrix}
\begin{bmatrix}
\bm{u}^{st} \\ \bm{p}^{st} \\ \bm{\lambda}^{st}
\end{bmatrix}
\ + \ 
\begin{bmatrix}
\bm{c}^{st}(\bm{u}^{st}) \\ \bm{0} \\ \bm{0} \\
\end{bmatrix}
\ + \
\begin{bmatrix} 
\bm{A_s}^{st}(\bm{\tilde{\mu}_m}) + \bm{E_s}^{st}  & \bm{0} & \bm{0}  \\ 
\bm{0} & \bm{0} & \bm{0} \\ 
\bm{0} & \bm{0} & \bm{0}
\end{bmatrix}
\begin{bmatrix}
\bm{d}^{st} \\ \bm{0} \\ \bm{0} \\
\end{bmatrix}
\ = \ 
\begin{bmatrix}  
\bm{0} \\ \bm{0} \\ \bm{F}^{st}_3
\end{bmatrix}~,
\end{equation} 
where the vectors $\bm{u}^{st} \in \R^{N_u^{st}}$, $\bm{p}^{st} \in \R^{N_p^{st}}$, $\bm{\lambda}^{st} \in \R^{N_\lambda^{st}}, \bm{d}^{st} \in \R^{N_u^{st}}$ collect, respectively, the full--order velocity, pressure, Lagrange multipliers and displacement discrete solutions, at the time instants $\{t_n\}_{n=1}^{N^t}$. Here, $N_u^{st}, N_p^{st}, N_\lambda^{st} \in \mathbb{N}$ denote the total number of spatio--temporal FOM DOFs for velocity, pressure, and Lagrange multipliers, respectively, so that $N^{st} = N_u^{st} + N_p^{st} + N_\lambda^{st}$. The matrices $\bm{A}_i^{st}$ ($i \in \{2, 3, 4, 7\}$) and the vector $\bm{F}_3^{st}$ are identical to the ones defined in \cite{tenderini2024space} for the Stokes problem (see Eqs.~(2.10),~(2.11)). The matrix $\bm{A}_1^{st} \in \R^{N_u^{st} \times N_u^{st}}$ differs from the one reported in~\cite{tenderini2024space} due to the resistance boundary conditions. Indeed, we have
\begin{equation*}
\bm{A}_1^{st} = \
\operatorname{diag}\biggl(\underbrace{\bm{M}, \cdots, \bm{M}}_{N^t}\biggr) + 
\beta \Delta t \ \operatorname{diag}\biggl(\underbrace{\bm{A+R}, \cdots, \bm{A+R}}_{N^t}\biggr) - 
\sum_{s=1}^S \alpha_s \  \operatorname{subdiag}^{(s)}\biggl(\underbrace{\bm{M}, \cdots, \bm{M} }_{N^t-s}\biggr)~, 
\end{equation*}
where the operator $\operatorname{diag}: \R^{r_1 \times c_1} \times \cdots \times \R^{r_K \times c_K} \to \R^{(r_1+\cdots+r_K)\times(c_1+\cdots+c_K)}$ builds a block diagonal matrix from a set of $K$ input matrices; $\operatorname{subdiag}^{(n)}$ ($n \in \mathbb{N}$) is equivalent to $\operatorname{diag}$, but with respect to the $n$--th subdiagonal. 
The term $\bm{c}^{st}(\bm{u}^{st}) \in \R^{N_u^{st}}$ represents the full--order nonlinear convective term evaluated at all the discrete time instants, multiplied by the positive constant $\beta \Delta t$.
Finally, the matrices $\bm{A_s}^{st}(\bm{\tilde{\mu}_m})$, $\bm{M_s}^{st}$, and $\bm{E_s}^{st}$, all of size $N_u^{st} \times N_u^{st}$ and stemming from the coupled momentum model, are defined as follows:
\begin{equation}
\label{eq: space-time FOM membrane matrices}
\begin{alignedat}{2}
\bm{A_s}^{st} (\bm{\tilde{\mu}_m}) &= \sum_{q=1}^2 (\bm{\tilde{\mu}_m})_{q+1} \ \bm{A_s}^{st, q} \quad \text{with} \quad
\bm{A_s}^{st, q} = \beta \Delta t \ \operatorname{diag}\biggl(\underbrace{\bm{A_s}^q, \cdots, \bm{A_s}^q}_{N^t}\biggr)~; \\
\bm{M_s}^{st} &=  \operatorname{diag}\biggl(\underbrace{\bm{M_s}, \cdots, \bm{M_s}}_{N^t}\biggr) - \sum_{s=1}^S \alpha_s \ \operatorname{subdiag}^{(s)}\biggl(\underbrace{\bm{M_s}, \cdots, \bm{M_s}}_{N^t-s}\biggr)~; \\
\bm{E_s}^{st} &= \beta \Delta t \ c_s \ \operatorname{diag}\biggl(\underbrace{\bm{M_s}, \cdots, \bm{M_s}}_{N^t}\biggr)~.
\end{alignedat}
\end{equation}
Adopting a more concise notation, we can rewrite Eq.\eqref{eq: monolitic_FOM_system} as follows:
\begin{equation}
\label{eq: monolithic FOM system concise}
\left(\bm{A}^{st} + \bm{A}^{\Gamma, st}(\bm{\tilde{\mu}_m})\right) \bm{w}^{st} + \mathcal{E}_u\left(\bm{c}^{st}(\bm{u}^{st})\right) + \mathcal{E}_u\left(\left(\bm{A_s}^{st} (\bm{\tilde{\mu}_m}) + \bm{E_s}^{st}\right) \ \bm{d}^{st}\right) = \bm{F}^{st}(\bm{\mu_f}) + \bm{F}^{\Gamma, st}(\bm{\tilde{\mu}_m})~,
\end{equation}
where $\bm{w}^{st} := \left[\bm{u}^{st}, \bm{p}^{st}, \bm{\lambda}^{st} \right] \in \R^{N^{st}}$ and $\bm{A}^{\Gamma, st}(\bm{\tilde{\mu}_m}) := \mathcal{E}_{u,u}(\bm{M_s}^{st}(\bm{\tilde{\mu}_m})) \in \R^{N^{st} \times N^{st}}$. The quantities $\bm{A}^{st} \in \R^{N^{st} \times N^{st}}$ and $\bm{F}^{st} \in \R^{N^{st}}$ are the left-hand side matrix and right-hand side vector components appearing in Eq.\eqref{eq: monolitic_FOM_system} that do not stem from the coupled momentum model. Here, the operators $\mathcal{E}_u: \R^{N_u^{st}} \to \R^{N^{st}}$ and $\mathcal{E}_{u,u}: \R^{N_u^{st} \times N_u^{st}} \to \R^{N^{st} \times N^{st}}$ denote the velocity zero--padding in one dimension and in two dimensions, respectively.

\begin{remark}
Eq.\eqref{eq: monolitic_FOM_system} is generic in what concerns the treatment of the membrane displacement. Such a formulation is convenient at this stage; in Section~\ref{sec:ST-RB for RFSI}, we will propose different strategies to express $\bm{d}^{st}$ in terms of the velocity field $\bm{u}^{st}$, depending on the chosen model order reduction approach.
\end{remark}

\begin{remark}
The definitions of the matrices $\bm{E_s}^{st}$ and $\bm{A_s}^{st} (\bm{\tilde{\mu}_m})$ in Eq.\eqref{eq: space-time FOM membrane matrices} entail that the membrane contributions are treated implicitly. In fact, in Eq.\eqref{eq: monolitic_FOM_system} we are not exploiting the kinematic coupling condition that links the displacement with the fluid velocity on $\Gamma$:  
\begin{equation}
\label{eq: velocity displacement relation}
\vec{d}_h^{n+1} = \beta \Delta t \  \vec{u}_h^{n+1} + \sum_{s=1}^S \alpha_s \ \vec{d}_h^{n+1-s} \qquad \forall n \in \{0, \cdots, N^t-1\}~.
\end{equation}
If the problem is solved iteratively over time, Eq.\eqref{eq: velocity displacement relation}  allows to obtain, at every time instant, the displacement as a model by--product, combining the current velocity value and the displacement values available at $S$ previous timesteps. We refer the reader to~\cite{figueroa2006coupled, figueroa2009effect, colciago2014comparisons} for further details on this aspect.
\end{remark}

\section{ST--RB methods for the Coupled Momentum model}
\label{sec:ST-RB for RFSI}
The traditional RB method operates only on the spatial dimensionality of the problem. Indeed, the projection is performed onto spatial subspaces, where the solution is sought at every time instant. Therefore, the temporal dimensionality of the problem does not change upon the reduction step, and the ODE stemming from the semi--discrete ROM is solved through iterative numerical integration. While in many applications this strategy is enough to guarantee a drastic downscaling on the computational costs, the efficacy of RB may not be entirely satisfactory in problems characterized by a large temporal complexity. For example, this is the case if either the simulation interval should be very large or the timestep size should be very small, to properly capture some relevant physical behaviours. In order to address such a temporal--complexity bottleneck, in this work we consider ST--RB methods~\cite{choi2019space, choi2021space}.

\subsection{Space--time model order reduction} \label{subsec: space-time model order reduction}
ST--RB methods are based on the assumption that good approximations of high--fidelity solutions over the entire spatio--temporal domain can be found within a linear subspace, spanned by a small number of spatio--temporal basis functions. To reduce the computational costs associated to the projection operations and the overall memory requirements, these methods express the space--time basis functions as multiplications of spatial and temporal modes. Let 
\begin{equation*}
\label{eq: reduced bases space}
\bm{\Phi}^u := \left[ \bm{\varphi}^u_1 \big\vert \cdots \big\vert \bm{\varphi}^u_{n_u^s} \right] \ \ \in \R^{N_u^s \times n_u^s}~; \qquad
\bm{\Phi}^p := \left[ \bm{\varphi}^p_1 \big\vert \cdots \big\vert \bm{\varphi}^p_{n_p^s} \right] \ \ \in \R^{N_p^s \times n_p^s}~; \qquad
\bm{\Phi}^\lambda := \left[ \bm{\varphi}^\lambda_1 \big\vert \cdots \big\vert \bm{\varphi}^\lambda_{n_\lambda^s} \right] \ \ \in \R^{N_\lambda \times n_\lambda^s}~;
\end{equation*}
encode the discrete reduced basis functions in space and
\begin{equation*}
\label{eq: reduced bases time}
\bm{\Psi}^u := \left[ \bm{\psi}^u_1 \big\vert \cdots \big\vert \bm{\psi}^u_{n_u^t} \right] \ \ \in \R^{N^t \times n_u^t}~; \qquad
\bm{\Psi}^p := \left[ \bm{\psi}^p_1 \big\vert \cdots \big\vert \bm{\psi}^p_{n_p^t} \right] \ \ \in \R^{N^t \times n_p^t}~; \qquad
\bm{\Psi}^\lambda := \left[ \bm{\psi}^\lambda_1 \big\vert \cdots \big\vert \bm{\psi}^\lambda_{n_\lambda^t} \right] \ \ \in \R^{N^t \times n_\lambda^t}~;
\end{equation*}
encode the discrete reduced basis functions in time, for velocity, pressure, and Lagrange multipliers, respectively. All the reduced bases are computed following the ST--HOSVD approach proposed in~\cite{choi2019space}. In particular, we apply a randomized truncated POD algorithm~\cite{halko2011finding} to the matrices storing high--fidelity solutions, corresponding to the parameter values $\{\bm{\mu}_k\}_{k=1}^{M}$. Incidentally, we remark that the ST--GRB framework does not require the offline full--order snapshots to be generated with a monolithic space--time model; standard sequential--in--time approaches can be employed for this purpose.
The generic parameter vector $\bm{\mu}_k$ is such that $\bm{\mu}_k := [\bm{\mu_f}^k, \bm{\mu_m}^k] \in \mathcal{P} \subset \R^{N_\mu^f + N_\mu^m}$, with $\bm{\mu_f}^k \in \mathcal{P}_f \subset \R^{N_\mu^f}$ and $\bm{\mu_m}^k \in \mathcal{P}_m \subset \R^{N_\mu^m}$. We highlight that all the temporal bases are orthonormal with respect to the Euclidean norm in $\R^{N^t}$, whereas the spatial bases for velocity and pressure are orthonormal with respect to the $[H^1(\Omega)]^3$--norm and to the $L^2(\Omega)$--norm, respectively. We refer the reader to~\cite{tenderini2024space} for a detailed explanation of the reduced bases generation process.\\

We can now define the spatio--temporal discrete reduced basis functions as follows:
\begin{equation*}
\label{eq: basis functions}
\begin{alignedat}{7}
&\bm{\pi}_\ell^u \in \R^{N_u^{st}} \quad &&\text{such that} \quad \bm{\pi}_\ell^u &&= \bm{\varphi}^u_{\ell_s} &&\otimes \bm{\psi}^u_{\ell_t} \quad
&&\text{with} \ \ \ell &&= (\ell_s-1)n_u^t + \ell_t  \ &&=: \ \mathcal{F}_u(\ell_s, \ell_t)~;\\
&\bm{\pi}_k^p \in \R^{N_p^{st}} \quad &&\text{such that} \quad \bm{\pi}_k^p &&= \bm{\varphi}^p_{k_s} &&\otimes \bm{\psi}^p_{k_t} \quad 
&&\text{with} \ \ k &&= (k_s-1)n_p^t + k_t \ &&=: \ \mathcal{F}_p(k_s, k_t)~;\\
&\bm{\pi}_i^\lambda \in \R^{N_\lambda^{st}} \quad &&\text{such that} \quad \bm{\pi}_i^\lambda &&= \bm{\varphi}^\lambda_{i_s} &&\otimes \bm{\psi}^\lambda_{i_t} \quad
&&\text{with} \ \ i &&= (i_s-1)n_\lambda^t + i_t \ &&=: \ \mathcal{F}_\lambda(i_s, i_t)~.
\end{alignedat}
\end{equation*}
Furthermore, we define the space--time--reduced dimensions $n_u^{st} := n_u^s n_u^t$, $n_p^{st} := n_p^s n_p^t$, $n_\lambda^{st} := n_\lambda^s n_\lambda^t$ and $n^{st} := n_u^{st} + n_p^{st} + n_\lambda^{st}$. We remark that the orthonormality of the reduced bases in space and in time induces the orthonormality of the spatio--temporal reduced bases, with respect to the natural norms in space--time. An approximate spatio--temporal basis for FOM solutions manifold can be then encoded in the matrix:
\begin{equation*}
\label{eq: reduced basis matrix}
\bm{\Pi} := \operatorname{diag} \left(\left[\bm{\pi}^u_1\lvert...\rvert\bm{\pi}^u_{n^{st}_u}\right], \left[\bm{\pi}^p_1\lvert...\rvert\bm{\pi}^p_{n^{st}_p}\right], \left[\bm{\pi}^\lambda_1\lvert...\rvert\bm{\pi}^\lambda_{n^{st}_\lambda}\right] \right) = \operatorname{diag} \left( \bm{\Pi}^u, \bm{\Pi}^p, \bm{\Pi}^\lambda \right) \quad \in \R^{N^{st} \times n^{st}}~.
\end{equation*}
The orthonormality of the spatio--temporal reduced bases qualifies $\bm{\Pi}^T$ as representative of an orthogonal projection operator from the FOM space of dimension $N^{st}$ to the ST--ROM space of dimension $n^{st} \ll N^{st}$. Ultimately, we can define the space--time--reduced velocity $\widehat{\bm{u}} \in \R^{n_u^{st}}$, pressure $\widehat{\bm{p}} \in \R^{n_p^{st}} $ and Lagrange multipliers $\widehat{\bm{\lambda}} \in \R^{n_\lambda^{st}}$ such that:
\begin{equation*}
\label{eq: space-time reduced solutions}
\bm{u}^{st} \approx \sum_{\ell=1}^{n_u^{st}} \widehat{\bm{u}}_\ell \bm{\pi}_\ell^u~; \qquad
\bm{p}^{st} \approx \sum_{k=1}^{n_p^{st}} \widehat{\bm{p}}_k \bm{\pi}_k^p~; \qquad 
\bm{\lambda}^{st} \approx \sum_{i=1}^{n_\lambda^{st}} \widehat{\bm{\lambda}}_i \bm{\pi}_i^\lambda~.
\end{equation*}
For the sake of simplicity, from now on we will adopt the following notation for the indices of space--time--reduced quantities:
\begin{equation*}
\label{eq: space-time indices}
\begin{alignedat}{6}
&\ell = \mathcal{F}_u(\ell_s, \ell_t), \ \ &&m = \mathcal{F}_u(m_s, m_t) \quad 
&&\text{with} \ \ \ell_s,m_s &&\in \{1,\cdots, n_u^s\}, \ \ &&\ell_t,m_t &&\in \{1,\cdots, n_u^t\}~; \\
&k = \mathcal{F}_p(k_s, k_t), \ \ &&r = \mathcal{F}_p(r_s, r_t) \quad 
&&\text{with} \ \  k_s, r_s &&\in \{1, \cdots, n_p^s\}, \ \ &&k_t, r_t &&\in \{1, \cdots, n_p^t\}~; \\
&i = \mathcal{F}_\lambda(i_s, i_t), \ \ &&j = \mathcal{F}_\lambda(j_s, j_t) \quad 
&&\text{with} \ \ i_s, j_s &&\in \{1, \cdots, N_\lambda\}, \ \ &&i_t, j_t &&\in \{1, \cdots, n_\lambda^t\}~.
\end{alignedat}
\end{equation*}

Compared to the Navier--Stokes equations, the coupled momentum model features the presence of an additional unknown: the vessel wall displacement $\vec{d}$. However, the latter can be embodied in the Navier--Stokes equations, exploiting its kinematic coupling with the fluid velocity. As discussed in Section~\ref{subsec: algebraic formulation}, if the problem is solved iteratively over time, this can be achieved exploiting Eq.\eqref{eq: velocity displacement relation}. However, such a strategy does not produce any computational advantage if the temporal dynamics are projected onto a low--dimensional subspace. Indeed, in this setting the reduced solution encodes the displacement values at all the discrete time instants and the displacement does not be seamlessly configure as a model by--product. However, to strongly enforce the kinematic coupling condition under space--time model order reduction, we can approximate the spatio--temporal displacement field $\bm{d}^{st} \in \R^{N_u^{st}}$ as
\begin{equation}
\label{eq: space-time displacement}
\bm{d}^{st} \approx \sum_{\ell=1}^{n_u^{st}} \widehat{\bm{u}}_\ell \bm{\pi}_\ell^d \ = \ 
\sum_{\ell=1}^{n_u^{st}} \widehat{\bm{u}}_\ell \left(\bm{\varphi}_{\ell_s}^u \otimes \mathcal{P}_0\left( \mathcal{E}^t_S\left(\bm{\psi}_{\ell_t}^u\right)\right) \right)~.
\end{equation}
In Eq.\eqref{eq: space-time displacement}, $\mathcal{E}^t_S: \R^{N^t} \to \R^{N^t+S}$ is a zero--padding operator, adding $S$ zero entries at the beginning of its input. The operator $\mathcal{P}_0 : \R^{N^t+S} \to \R^{N^t}$ denotes instead the numerical approximation of the primitive through the chosen multistep method, under the assumption that its input is equal to zero at the first $S$ timesteps. In particular, the following recursive relation holds:
\begin{equation*}
\label{eq: discrete primitive}
\big(\mathcal{P}_0(\bm{v})\big)_{n+1} = \beta \Delta t \ \bm{v}_{n+1} + \sum_{s=1}^S \alpha_s \left(\mathcal{P}_0(\bm{v})\right)_{n+1-s} \qquad \forall n \in \{0, \cdots, N^t-1\}~.
\end{equation*}
Eq.\eqref{eq: space-time displacement} allows to express the displacement as a function of the reduced velocity coefficients, while strongly enforcing the kinematic coupling condition. Therefore, the displacement no longer configures as a problem unknown and can be retrieved as a by--product during the post--processing phase. In fact, this is achieved through a projection onto a convenient spatio--temporal subspace. The latter is defined as the product space, generated by the velocity reduced subspace in space and by the subspace spanned by the primitives of the velocity modes in time. The vectors $\{\bm{\pi}_\ell^d\}_{\ell=1}^{n_u^{st}}$ are stored as the columns of the matrix 
\begin{equation}
\label{eq: space-time displacement 2}
\bm{\Pi}^d := \bm{\Phi}^u \otimes \mathcal{P}_0(\mathcal{E}_S^t(\bm{\Psi}^u)) \ \ \in \R^{N_u^{st} \times n_u^{st}}~,
\end{equation}
that encodes an approximate low--dimensional basis for the FOM displacement solutions manifold. In Eq.\eqref{eq: space-time displacement 2}, the operators $\mathcal{P}_0$ and $\mathcal{E}_S^t$ are defined as in Eq.\eqref{eq: space-time displacement}, but they are applied columnwise. Thanks to Eq.\eqref{eq: space-time displacement 2}, we can write $\bm{d}^{st} \approx \bm{\Pi}^d \widehat{\bm{u}}$. \\

\subsection{ST--RB problem definition} \label{subsec: ST-RB problem definition}
From an algebraic standpoint, the application of the ST--RB method to Eq.\eqref{eq: monolithic FOM system concise} requires to find $\widehat{\bm{w}} = [\widehat{\bm{u}}, \widehat{\bm{p}}, \widehat{\bm{\lambda}}]\in \R^{n^{st}}$ such that
\begin{equation}
\label{eq: ST-RB definition}
\widehat{\bm{r}}(\widehat{\bm{w}}) := \ \mytilde{\bm{\Pi}}^T \Bigg( \bm{F}^{st}(\bm{\mu_f}) - 
\left( \bm{A}^{st} + \bm{A}^{\Gamma, st}(\bm{\tilde{\mu}_m})\right) \bm{\Pi}\widehat{\bm{w}} - \mathcal{E}_u \Big(\bm{c}^{st}\left(\bm{\Pi}^u\widehat{\bm{u}}\right) + \left(\bm{A_s}^{st} (\bm{\tilde{\mu}_m}) + \bm{E_s}^{st}\right) \bm{\Pi}^d \widehat{\bm{u}} \Big) \Bigg) = \bm{0}~.
\end{equation}
Here $\mytilde{\bm{\Pi}} \ \in \R^{N^{st}\times n^{st}}$ is a possibly parameter--dependent matrix, while $\mathcal{E}_u$ is the zero--padding operator used in Eq.\eqref{eq: monolithic FOM system concise}. 
Eq.\eqref{eq: ST-RB definition} prescribes the orthogonality between the approximate FOM residual --- evaluated at the high--fidelity reconstruction of the space--time--reduced solution $\widehat{\bm{w}}^{st} := \bm{\Pi}\widehat{\bm{w}} \in \R^{N^{st}}$ --- and a low--dimensional subspace, spanned by the columns of $\mytilde{\bm{\Pi}}$. Different ST--RB methods stem from different choices of $\mytilde{\bm{\Pi}}$. In \cite{tenderini2024space}, two special cases have been investigated: the first stems from the trivial choice $\mytilde{\bm{\Pi}} \ = \bm{\Pi}$ (Galerkin projection, ST--GRB method), whereas the second arises from the least--squares minimization of the FOM residual in a suitable norm (least--squares Petrov--Galerkin projection, ST--PGRB method). In this work, we focus only on the ST--GRB method. Nonetheless, we highlight that a similar analysis can be carried out whenever the matrix $\mytilde{\bm{\Pi}}$ is parameter--independent. It is worth remarking that, in the case of Galerkin projections, suitable augmentations of the spatial and temporal velocity reduced basis are required to guarantee the \emph{inf--sup} stability, and hence the well--posedness, of the problem. We refer the reader to \cite{rozza2005optimization} for the supremizers enrichment procedure in space and to \cite{tenderini2024space} for the stabilizers enrichment procedure in time. \\

Eq.\eqref{eq: ST-RB definition} describes a $n^{st}$--dimensional nonlinear system, whose nonlinearity is due to the presence of the convective term $\bm{c}(\bm{\Pi}^u \widehat{\bm{u}})$. If the Navier--Stokes equations are solved iteratively over time, it is a popular strategy to linearize the convective term, by extrapolation in time of the velocity solutions available at previous time instants. However, as already discussed concerning the treatment of the coupled momentum model contributions, the use of a monolithic space--time formulation to enable the compression of the temporal dynamics nullifies the computational gains of this approach. Therefore, the nonlinear nature of the problem cannot be easily circumvented and Eq.\eqref{eq: ST-RB definition} should be solved adopting \emph{ad hoc} numerical techniques. In this work, we rely on the Newton--Raphson method. Let $\widehat{\bm{w}}^{(0)} \in \R^{n^{st}}$ be a given initial guess; the superscript denotes the iteration index. Then, the solution to Eq.\eqref{eq: ST-RB definition} is iteratively found as follows:
\begin{equation}
\label{eq: Newton iterations}
\widehat{\bm{w}}^{(l+1)} = \widehat{\bm{w}}^{(l)} - \left(\widehat{\bm{J}}_{\widehat{\bm{r}}}\big(\widehat{\bm{w}}^{(l)}\big)\right)^{-1}\widehat{\bm{r}}\big(\widehat{\bm{w}}^{(l)}\big) \qquad \text{with} \quad 
\left(\widehat{\bm{J}}_{\widehat{\bm{r}}}\big(\widehat{\bm{w}}^{(l)}\big)\right)_{ij} := \left.\frac{\partial \widehat{\bm{r}}_i}{\partial \widehat{\bm{w}}_j}\right\vert_{\widehat{\bm{w}} = \widehat{\bm{w}}^{(l)}}~,
\end{equation}
where the matrix $\widehat{\bm{J}}_{\widehat{\bm{r}}} \in \R^{n^{st} \times n^{st}}$ denotes the Jacobian of the residual $\widehat{\bm{r}} \in \R^{n^{st}}$ with respect to the space--time--reduced solution $\widehat{\bm{w}}$. From Eq.\eqref{eq: ST-RB definition}, setting $\mytilde{\bm{\Pi}} = \bm{\Pi}$ (ST--GRB method), we have that
\begin{equation}
\label{eq: reduced Jacobian matrix}
\widehat{\bm{J}}_{\widehat{\bm{r}}}\big(\widehat{\bm{w}}\big) = - \bm{\Pi}^T \left(\bm{A}^{st} + \bm{A}^{\Gamma, st}(\bm{\tilde{\mu}_m})\right) \bm{\Pi} - \bm{\Pi}^T \mathcal{E}_{u,u}\Big( \bm{J_c}^{st}(\bm{\Pi}^u \widehat{\bm{u}}) \Big) \bm{\Pi} - \bm{\Pi}^T \mathcal{E}_{u,u} \Big( \bm{A_s}^{st} (\bm{\tilde{\mu}_m}) + \bm{E_s}^{st}\Big) \bm{\Pi}^d~,
\end{equation} 
where $\bm{J_c}^{st}(\bm{\Pi}^u \widehat{\bm{u}}) \in \R^{N_u^{st} \times N_u^{st}}$ is a block--diagonal matrix containing the full--order convective Jacobian matrices, evaluated at the reconstructed reduced velocity solution for all the discrete time instants and multiplied by $\beta \Delta t$. Also, $\mathcal{E}_{u,u}$ is the zero--padding operator used in Eq.\eqref{eq: monolithic FOM system concise}. The iterations in Eq.\eqref{eq: Newton iterations} are executed until the following condition is satisfied:
\begin{equation}
\label{eq: Newton stopping criterion}
\left. \left\Vert \widehat{\bm{r}}\big(\widehat{\bm{w}}^{(l)}\big) \right\Vert_2 \ \middle/ \ \left\Vert \widehat{\bm{r}}\big(\widehat{\bm{w}}^{(0)}\big) \right\Vert_2 \right. \leq \tau_{NR}~,
\end{equation}
where $\tau_{NR} \in \R^+$ is a user--provided tolerance. If the stopping criterion in Eq.\eqref{eq: Newton stopping criterion} is not reached after $K_{NR} \in \mathbb{N}$ iterations (at most), then we say that the Newton--Raphson method does not converge.

\subsection{Assembling the reduced system} \label{subsec: assembling the reduced system}
In this section, we discuss how the space--time--reduced residual $\widehat{\bm{r}}(\widehat{\bm{w}})$ (see Eq.\eqref{eq: ST-RB definition}) and its Jacobian $\widehat{\bm{J}}_{\widehat{\bm{r}}}\big(\widehat{\bm{w}}\big) $ (see Eq.\eqref{eq: reduced Jacobian matrix}) can be efficiently assembled. It is worth remarking that the affine parametrization of the problem enables the offline computation of several space--time--reduced quantities. This allows to drastically lower the online computational costs and to maximize speedups at no loss in accuracy. In the general case, approximate affine decompositions can be retrieved via the Discrete Empirical Interpolation Method (DEIM)~\cite{chaturantabut2010nonlinear} or its matricial version MDEIM~\cite{negri2015efficient}. \\

To simplify the notation, we rewrite the space--time--reduced residual and its Jacobian as follows:
\begin{subequations}
\label{eq: space-time-reduced residual and Jacobian}
\begin{align}
\label{eq: space-time-reduced residual}
\widehat{\bm{r}}(\widehat{\bm{w}}) &= 
\widehat{\bm{F}}(\bm{\mu_f})
- \left(\widehat{\bm{A}} + \widehat{\bm{A}}^{\Gamma}(\bm{\tilde{\mu}_m})\right) \widehat{\bm{w}} 
- \widehat{\mathcal{E}}_u \big(\widehat{\bm{c}}(\widehat{\bm{u}})\big) - \widehat{\mathcal{E}}_u \big( (\widehat{\bm{A}}_{\bm{s}} (\bm{\tilde{\mu}_m}) + \widehat{\bm{E}}_{\bm{s}}) \ \widehat{\bm{u}} \big)~; \\
\label{eq: space-time-reduced Jacobian}
\widehat{\bm{J}}_{\widehat{\bm{r}}}\big(\widehat{\bm{w}}\big) &=
- \left(\widehat{\bm{A}} + \widehat{\bm{A}}^{\Gamma}(\bm{\tilde{\mu}_m})\right) - \widehat{\mathcal{E}}_{u,u}\big(\widehat{\bm{J}}_{\bm{c}}(\widehat{\bm{u}})\big) - \widehat{\mathcal{E}}_{u,u}\big(\widehat{\bm{A}}_{\bm{s}} (\bm{\tilde{\mu}_m}) + \widehat{\bm{E}}_{\bm{s}}\big)~;
\end{align}
\end{subequations} 
where $\mathcal{\widehat{E}}_u: \R^{n_u^{st}} \to \R^{n^{st}}$ and $\mathcal{\widehat{E}}_{u,u}: \R^{n_u^{st}\times n_u^{st}} \to \R^{n^{st}\times n^{st}}$ are the vectorial and the matricial space--time--reduced velocity zero--padding operators, respectively. The quantities appearing in Eq.\eqref{eq: space-time-reduced residual and Jacobian} will be suitably defined throughout this section. Before detailing how the reduced system is assembled, let us introduce the following space--reduced quantities:
\begin{itemize}[itemsep=1pt]
\item the matrices $\widebar{\bm{M}} \in \R^{n_u^s \times n_u^s}$, $\widebar{\bm{A}} \in \R^{n_u^s \times n_u^s}$, $\widebar{\bm{B}} \in \R^{n_u^s \times n_p^s}$, $\widebar{\bm{L}} \in \R^{n_u^s \times N_\lambda}$, $\widebar{\bm{B}}^T \in \R^{n_p^s \times n_u^s}$, $\widebar{\bm{L}}^T \in \R^{N_\lambda \times n_u^s}$, $\widebar{\bm{R}} \in \R^{n_u^s \times n_u^s}$, that are the space--reduced counterparts of the FOM matrices $\bm{M}$, $\bm{A}$, $\bm{B}$, $\bm{L}$, $\bm{B}^T$, $\bm{L}^T$, $\bm{R}$, respectively. For example, $\widebar{\bm{M}} := (\bm{\Phi}^u)^T \bm{M} \ \bm{\Phi}^u$.
\item the matrices $\widebar{\bm{M}}_{\bm{s}} \in \R^{n_u^s \times n_u^s}$ and $\{\widebar{\bm{A}}^q_{\bm{s}} \}_{q=1}^2 \in \R^{n_u^s \times n_u^s}$, that are the space--reduced counterparts of the FOM structure matrices $\bm{M_s}$ and $\{\bm{A_s}^q\}_{q=1}^2$, respectively. In particular, $\widebar{\bm{M}}_{\bm{s}} := (\bm{\Phi}^u)^T \bm{M_s} \bm{\Phi}^u$;  $\widebar{\bm{A}}^q_{\bm{s}} := (\bm{\Phi}^u)^T \bm{A}^q_{\bm{s}} \bm{\Phi}^u$ for $q=1,2$.
\end{itemize}

\subsubsection{Space--time--reduced linear terms} \label{subsubsec: ST-reduced matrices} 
The assembling of the parameter--independent space--time--reduced left--hand side matrix $\widehat{\bm{A}} := \bm{\Pi}^T \bm{A} \ \bm{\Pi} \in \R^{n^{st} \times n^{st}}$ is analogous to the Stokes case~\cite{tenderini2024space}, besides including the resistance boundary conditions contributions. Hence, $\widehat{\bm{A}}$ has the following block structure
\begin{equation*}
\label{eq: stokes matrix ST-RB}
\widehat{\bm{A}} 
=
\begin{bmatrix}
\widehat{\bm{A}}_1 & \widehat{\bm{A}}_2
& \widehat{\bm{A}}_3 \\
\widehat{\bm{A}}_4 & \bm{0} & \bm{0} \\
\widehat{\bm{A}}_7 & \bm{0} & \bm{0}
\end{bmatrix}~,
\end{equation*}
and the entries of the matrices in the non--zero blocks are
\begin{equation*}
\label{eq: stokes ROM spacetime blocks}
\begin{alignedat}{2}
\left(\widehat{\bm{A}}_1\right)_{\ell m} &= \ \left(\widebar{\bm{M}} + \beta \Delta t \ (\widebar{\bm{A}} + \widebar{\bm{R}})\right)_{\ell_s m_s}\delta_{\ell_t,m_t} \ 
- \ \sum_{s=1}^S \alpha_s \ \widebar{\bm{M}}_{\ell_s m_s} \left((\bm{\psi}^u_{\ell_t})_{s+1:}^T(\bm{\psi}^u_{m_t})_{:-s}\right) , \\
\left(\widehat{\bm{A}}_2\right)_{\ell r} &= \ \beta \Delta t \ \widebar{\bm{B}}^T_{\ell_s r_s} \left( (\bm{\psi}^u_{\ell_t})^T (\bm{\psi}^p_{r_t}) \right) , 
\qquad \qquad  
\left(\widehat{\bm{A}}_4\right)_{k m} = \ \widebar{\bm{B}}_{k_s m_s}  \left( (\bm{\psi}^p_{k_t})^T (\bm{\psi}^u_{m_t}) \right) , \\ 
\left(\widehat{\bm{A}}_3\right)_{\ell j} &=  \ \beta \Delta t \ \widebar{\bm{L}}^T_{\ell_s j_s} \left( (\bm{\psi}^u_{\ell_t})^T (\bm{\psi}^\lambda_{j_t}) \right) ,
\qquad \qquad  \
\left(\widehat{\bm{A}}_7\right)_{i m} = \ \widebar{\bm{L}}_{i_s m_s} \left( (\bm{\psi}^\lambda_{i_t})^T (\bm{\psi}^u_{m_t}) \right)~. \\
\end{alignedat}
\end{equation*}
Here the notations $\bm{v}_{i:}$, $\bm{v}_{:-i}$ denote the sub--vector of a given vector $\bm{v}$ containing all the entries from the $i$--th to the last one and from the first one to the $i$--th from last, respectively. From Eq.\eqref{eq: monolitic_FOM_system}, we notice that the coupled momentum model adds a parametrized boundary mass term in the velocity--velocity block. Therefore, we define the matrix $\widehat{\bm{A}}^{\Gamma}(\bm{\tilde{\mu}_m}) := (\bm{\tilde{\mu}_m})_1 \ \widehat{\mathcal{E}}_{u,u}( \widehat{\bm{M}}_{\bm{s}})$; exploiting the tridiagonal block structure of $\bm{M}_{\bm{s}}^{st}$ (see Eq.\eqref{eq: space-time FOM membrane matrices}), we have 
\begin{equation*}
\label{eq: space--time reduced boundary mass}
\left(\widehat{\bm{M}}_{\bm{s}}\right)_{\ell m} = (\widebar{\bm{M}}_{\bm{s}})_{\ell_s m_s}\delta_{\ell_t,m_t}
- \sum_{s=1}^S \alpha_s (\widebar{\bm{M}}_{\bm{s}})_{\ell_s m_s}\left((\bm{\psi}^u_{\ell_t})_{s+1:}^T(\bm{\psi}^u_{m_t})_{:-s}\right)~.
\end{equation*}

The coupled momentum model introduces a boundary stiffness matrix featuring a single non--zero block in the velocity--velocity compartment. Since the space--time--reduced basis matrix $\bm{\Pi}$ is block--diagonal, the same block structure is maintained upon the projection. Leveraging the affine parametrization and the block--diagonal structure of $\bm{A}_{\bm{s}}^{st}$ (see Eq.\eqref{eq: space-time FOM membrane matrices}), and exploiting the space--time approximation of the displacement field reported in Eq.\eqref{eq: space-time displacement}, we have that
\begin{equation*}
\label{eq: space--time reduced boundary stiffness}
\widehat{\bm{A}}_{\bm{s}} (\bm{\tilde{\mu}_m}) = \beta \Delta t \ \sum_{q=1}^2 (\bm{\tilde{\mu}_m})_{q+1} \ \widehat{\bm{A}}_{\bm{s}}^{q} \quad \text{with} \ \ \left(\widehat{\bm{A}}_{\bm{s}}^{q}\right)_{\ell m} = \Big(\widebar{\bm{A}}_{\bm{s}}^q\Big)_{\ell_s m_s} \left( \left(\bm{\psi}^u_{\ell_t}\right)^T \left(\mathcal{P}_0\left(\mathcal{E}_S^t\left(\bm{\psi}^u_{m_t}\right)\right)\right) \right).
\end{equation*}
Similarly, the surrounding tissues contribution is encoded in the matrix $\widehat{\bm{E}}_{\bm{s}}$, whose entries are
\begin{equation*}
\label{eq: space--time reduced surrounding tisues}
\left(\widehat{\bm{E}}_{\bm{s}}\right)_{\ell m} = \beta \Delta t \ c_s \ \Big(\widebar{\bm{M}}_{\bm{s}}\Big)_{\ell_s m_s} \left( \left(\bm{\psi}^u_{\ell_t}\right)^T \left(\mathcal{P}_0\left(\mathcal{E}_S^t\left(\bm{\psi}^u_{m_t}\right)\right)\right) \right).
\end{equation*}

 
The space--time--reduced residual $\widehat{\bm{r}}(\widehat{\bm{w}})$ features a single non--zero block $\widehat{\bm{F}}_3(\bm{\mu_f}) := \bm{\Pi}_\lambda^T \bm{F}_3^{st}(\bm{\mu_f})$ in right--hand side vector. 
The non--affine parameter--dependent nature of $\bm{F}_3^{st}$ implies that $\widehat{\bm{F}}_3$ cannot be pre--assembled during the offline phase. However, we can take advantage of the space--time factorization of the Dirichlet datum in Eq.\eqref{eq: dirichlet datum} and, in particular, of the fact that the parametric dependency characterizes only the (discretized) flow rates $\bm{g}_k^t(\bm{\mu_f}) \in \R^{N^t}$, with $k=1, \cdots, N_D$. Indeed, the $k$--th block of the space--time reduced right--hand side $\widehat{\bm{F}}_{3,k}(\bm{\mu_f})$ is such that
\begin{equation*}
\label{eq: reduced RHS local} 
\widehat{\bm{F}}_{3,k}(\bm{\mu_f}) \in \R^{n^{st}_{\lambda_k}} \quad \text{with} \ \ \left(\bm{\widehat{F}}_{3,k}\right)_{i^k} = \beta \Delta t \ \left(\left( \varphi_{i_s^k}^{\lambda_k} \right)^T \bm{g}_k^s \right) \  \left(\left(\bm{\psi}_{i^k_t}^{\lambda_k}\right)^T \bm{g}_k^t(\bm{\mu_f})\right)~,
\end{equation*}
where $i^k = \mathcal{F}_{\lambda_k}(i^k_s, i^k_t) := (i^k_s-1)n_{\lambda_k}^t + i^k_t$. Ultimately, leveraging the product structure of the Lagrange multipliers subspace, we have
\begin{equation*}
\label{eq: reduced RHS global}
\widehat{\bm{F}}_{3}(\bm{\mu_f}) = \left[\left(\widehat{\bm{F}}_{3,1}(\bm{\mu_f})\right)^T, \cdots, \left(\widehat{\bm{F}}_{3,N_D}(\bm{\mu_f})\right)^T\right]^T \in \R^{n^{st}_{\lambda}}~.
\end{equation*}
We remark that this assembling step is extremely cheap, since it involves $N_D$ matrix--vector products of dimension $N^t$. 

\subsubsection{Space--time--reduced nonlinear convective term} \label{subsubsec: ST-reduced convective term} 
The presence of the nonlinear convective term is what differentiates the Navier--Stokes equations from the Stokes equations. As discussed in Section~\ref{subsec: ST-RB problem definition}, the linearization of the convective term by means of temporal extrapolation does not load to significant computational gains within the ST--RB framework; so, a fully--implicit approach is compulsory. The most straightforward strategy is to assemble the full--order convective term and project it onto the velocity reduced subspace. However, this inherently comes with several drawbacks. First of all, the associated computational costs are high, since they inherently depend on the high--fidelity dimensions in space--time. Indeed, at every Newton iteration, one should expand the space--time--reduced solution to the FOM space, compute the FOM convective term at all the discrete time instants and project it back again to the space--time--reduced subspace. Secondly, this procedure is intrusive at both the offline and the online stage, where by ``intrusivity'' we refer to the use of the FOM software during the ROM pipeline. While offline intrusivity characterizes the RB method, online intrusivity can represent a big limitation. Indeed, efficient communication routines between the FOM and the ROM softwares should be implemented, which not be trivial if closed--source programs are employed for the high--fidelity simulations. \\

To overcome the pitfalls of the baseline assembling approach, in this work we proceed as in ~\cite{pegolotti2021model}, directly exploiting the bilinearity of the convective term. Incidentally, we refer the reader to~\cite{pegolotti2021model} for a detailed comparison of the proposed approach with the non--hyper--reduced model, within a canonical space--only RB setting. For $\vec{u}_h, \vec{v}_h, \vec{w}_h \in \mathcal{V}_h$, the following affine decomposition holds:
\begin{equation}
\label{eq: convective term affine decomposition FOM}
\rho_f \int_\Omega \left( \left(\vec{u}_h \cdot \nabla \right) \vec{w}_h \right) \cdot \vec{v}_h =: a_{con}\left(\vec{u}_h, \vec{w}_h; \vec{v}_h\right) = \rho_f \sum_{i=1}^{N_u^s} \sum_{j=1}^{N_u^s} \bm{u}_i \ \bm{w}_j \ a_{con}\left(\vec{\varphi}^u_i, \vec{\varphi}^u_j; \vec{v}_h\right)~,
\end{equation}
where $\bm{u}, \bm{w} \in \R^{N_u^s}$ are the FOM discrete representations of $\vec{u}_h$, $\vec{w}_h$.
Even though Eq.\eqref{eq: convective term affine decomposition FOM} is written considering high--fidelity discretizations of the velocity, its extension to the space--reduced framework is trivial. Indeed, adopting an algebraic perspective, we can write the space--reduced convective term $\bm{\widebar{c}} \in \R^{n_u^s}$ as follows:
\begin{equation}
\label{eq: convective term affine decomposition ROM}
\bm{\widebar{c}}(\bm{\widebar{u}}) = \beta \Delta t \ \sum_{\ell_s'=1}^{n_u^s} \sum_{\ell_s''=1}^{n_u^s} \bm{\widebar{u}}_{\ell_s'} \bm{\widebar{u}}_{\ell_s''} \bm{\widebar{k}}_{\ell_s' \ell_s''} \quad \text{with} \ \  \left(\bm{\widebar{k}}_{\ell_s' \ell_s''}\right)_{m_s} =  a_{con}\left(\vec{\phi}^u_{\ell_s'}, \vec{\phi}^u_{\ell_s''}; \vec{\phi}^u_{m_s}\right)~,
\end{equation}
where $\vec{\phi}^u_{\ell_s} \in \mathcal{V}_h$ is such that $\vec{\phi}^u_{\ell_s}(\vec{x}) := \sum_{i=1}^{N_u^s} (\bm{\Phi}^u)_{i \ell_s} \ \vec{\varphi}^u_i(\vec{x})$ is the $\ell_s$--th velocity reduced basis function in space and $\bm{\widebar{u}} \in \R^{N^u_s}$ is the discrete representation of the space--reduced velocity field at an arbitrary time instant. We remark that the vectors $\{\bm{\widebar{k}}_\ell\}_{\ell=1}^{(n_u^s)^2}$ are independent from the solution and thus can be pre--assembled during the offline phase. Hence, Eq.\eqref{eq: convective term affine decomposition ROM} shows that the space--reduced convective term exhibits an exact affine decomposition, whose weights stem from the outer product of the space--reduced solution $\bm{\widebar{u}}$ with itself. A key implication of Eq.\eqref{eq: convective term affine decomposition ROM} is that the number of affine components depends quadratically on the number of velocity reduced bases, possibly limiting the computational advantages of the method. In this regard, to prevent efficiency losses, we can leverage the fact that the spatial velocity bases are sorted in decreasing order of ``importance'', since they have been computed through POD. Therefore, we can approximate the reduced nonlinear convective term with the following truncated expansion:
\begin{equation}
\label{eq: convective term approximate affine decomposition ROM}
\bm{\widebar{c}}(\bm{\bar{u}}) \ \approx \ \beta \Delta t \ \sum_{\ell_s'=1}^{n_c} \sum_{\ell_s''=1}^{n_c} \bm{\widebar{u}}_{\ell_s'} \ \bm{\widebar{u}}_{\ell_s''} \ \bm{\widebar{k}}_{\ell_s' \ell_s''}~, \quad \text{where} \ \ 0 \leq n_c \leq n_u^s~.
\end{equation}
The effect of the choice of $n_c$ on the efficiency and accuracy of the method will be empirically investigated in Section~\ref{subsec: navier-stokes test}. \\

To embed the nonlinear convective term in the ST--RB problem formulation of Eq.\eqref{eq: ST-RB definition}, we need to project Eq.\eqref{eq: convective term approximate affine decomposition ROM} onto the space--time--reduced velocity subspace. To this aim, we observe that the space--reduced solution $\bm{\widebar{u}}_n$ at time $t_n$ can be written as a function of the space--time--reduced solution $\widehat{\bm{u}}$ as follows:
\begin{equation}
\label{eq: space space-time velocity relation}
\bm{\widebar{u}}_n \in \R^{n_u^s} : \ \ (\bm{\widebar{u}}_n)_{\ell_s} = \sum_{\ell_t=1}^{n_u^t} \widehat{\bm{u}}_\ell \ (\bm{\psi}^u_{\ell_t})_n \quad \text{with} \ \ \ell = \mathcal{F}_u(\ell_s, \ell_t)~.
\end{equation}
Therefore, the space--time--reduced convective term $\widehat{\bm{c}}(\widehat{\bm{u}}) \in \R^{n_u^{st}}$ can be approximated as follows:
\begin{equation}
\label{eq: convective term affine decomposition ST-ROM}
\begin{split}
\left(\widehat{\bm{c}}(\widehat{\bm{u}})\right)_m &\approx 
\beta \Delta t \ \sum_{n=1}^{N^t} \left(\sum_{\ell_s'=1}^{n_c} \sum_{\ell_s''=1}^{n_c}  \left( \sum_{\ell_t'=1}^{n_u^t} \widehat{\bm{u}}_{\ell'} (\bm{\psi}^u_{\ell_t'})_n \right) \left( \sum_{\ell_t''=1}^{n_u^t} \widehat{\bm{u}}_{\ell''} (\bm{\psi}^u_{\ell_t''})_n \right) \bm{\widebar{k}}_{\ell_s'\ell_s''} \right)_{m_s} (\bm{\psi}^u_{m_t})_n \\
&=
\beta \Delta t \ \sum_{\ell_s'=1}^{n_c} \sum_{\ell_s''=1}^{n_c} \left(\bm{\widebar{k}}_{\ell_s'\ell_s''}\right)_{m_s} \sum_{\ell_t'=1}^{n_u^t} \sum_{\ell_t''=1}^{n_u^t} \widehat{\bm{u}}_{\ell'} \widehat{\bm{u}}_{\ell''} \ \psi^{u,3}_{\ell_t',\ell_t'',m_t}~,
\end{split}
\end{equation}
where $\psi^{u,3}_{\ell_t',\ell_t'',m_t} := \sum_n (\bm{\psi}^u_{\ell_t'})_n (\bm{\psi}^u_{\ell_t''})_n (\bm{\psi}^u_{m_t})_n$. 
Since the quantities $\bm{\widebar{k}}_{\ell_s' \ell_s''}$ and $\psi^{u,3}_{\ell_t',\ell_t'',m_t}$ can be pre--assembled offline, exploiting tensorial calculus the computation of $ \ \widehat{\bm{c}}(\widehat{\bm{u}})$ reduces to a linear combination and becomes relatively cheap. \\

Solving the ST--RB problem with the Newton--Raphson method requires to assemble the space--time--reduced convective Jacobian matrix $\widehat{\bm{J}}_{\bm{c}} \in \R^{n_u^{st} \times n_u^{st}}$. To this aim, we exploit that the high--fidelity convective Jacobian operator $J_{con}: \mathcal{V}_h \to \R$ is linear; indeed for $\vec{u}_h, \ \vec{v}_h', \ \vec{v}_h'' \ \in \mathcal{V}_h$, we have
\begin{equation*}
\label{eq: convective Jacobian operator}
J_{con}(\vec{u}_h; \vec{v}_h', \vec{v}_h'') = \int_\Omega \rho_f \left( \left(\vec{v}_h'' \cdot \nabla \right) \vec{u}_h \right) \cdot \vec{v}_h' + \int_\Omega \rho_f \left( \left(\vec{u}_h \cdot \nabla \right) \vec{v}_h'' \right) \cdot \vec{v}_h'~.
\end{equation*}
Hence, the following affine decomposition trivially holds for the space--reduced convective Jacobian $\bm{\widebar{J}_c}(\bm{\bar{u}}) \in \R^{n_u^s \times n_u^s}$:
\begin{equation*}
\label{eq: convective Jacobian affine decomposition ROM}
\bm{\widebar{J}_c}(\bm{\bar{u}}) = \sum_{\ell_s=1}^{n_u^s} \bm{\widebar{u}}_{\ell_s} \bm{\widebar{K}}_{\ell_s} \approx \sum_{\ell_s=1}^{n_{c,J}} \bm{\widebar{u}}_{\ell_s} \bm{\widebar{K}}_{\ell_s} \quad \text{with} \ \ \left(\bm{\widebar{K}}_{\ell_s}\right)_{m_s'm_s''} = J_{con}\left(\vec{\phi}_{\ell_s}^u; \vec{\phi}^u_{m_s'}, \vec{\phi}^u_{m_s''}\right)~,
\end{equation*}
where $0\leq n_{c,J} \leq n_u^s$. Exploiting Eq.\eqref{eq: space space-time velocity relation}, we can then approximate the space--time--reduced convective Jacobian matrix $\widehat{\bm{J}}_{\bm{c}}(\widehat{\bm{u}}) \in \R^{n_u^{st} \times n_u^{st}}$ as follows:
\begin{equation}
\label{eq: convective Jacobian affine decomposition ST-ROM}
\begin{split}
\left(\widehat{\bm{J}}_{\bm{c}}(\widehat{\bm{u}})\right)_{m'm''} &\approx 
\beta \Delta t \ \sum_{n=1}^{N^t} \left(\sum_{\ell_s=1}^{n_{c,J}} \left(\sum_{\ell_t=1}^{n_u^t} \widehat{\bm{u}}_\ell (\bm{\psi}^u_{\ell_t})_n \right) \bm{\widebar{K}}_{\ell_s}  \right)_{m_s'm_s''} (\bm{\psi}^u_{m_t'})_n (\bm{\psi}^u_{m_t''})_n \\   
&=
\beta \Delta t \ \sum_{\ell_s=1}^{n_{c,J}} \left(\bm{\widebar{K}}_{\ell_s}\right)_{m_s' m_s''} \sum_{\ell_t=1}^{n_u^t} \widehat{\bm{u}}_\ell \ \psi^{u,3}_{\ell_t, m_t', m_t''}~,
\end{split}
\end{equation}
where the coefficients $\psi^{u,3}_{\ell_t, m_t', m_t''}$ are defined as in Eq.\eqref{eq: convective term affine decomposition ST-ROM}. Once again, the possibility to pre--compute the quantities $\bm{\bar{K}}_{\ell_s}$ and $\psi^{u,3}_{\ell_t, m_t', m_t''}$ makes the assembling of $\widehat{\bm{J}}_{\bm{c}}(\widehat{\bm{u}})$ relatively inexpensive.

\begin{remark}
The assembly of the space--reduced convective term (and of its Jacobian) could also be performed using alternative hyper--reduction techniques. In this context, a widely adopted method in model order reduction is the Discrete Empirical Interpolation Method (DEIM)~\cite{barrault2004empirical, chaturantabut2010nonlinear}, whose extension to matrices (MDEIM) has also been proposed~\cite{negri2015efficient}. Notably, the application of (M)DEIM in a space--time--reduced setting has already been thoroughly investigated in \cite{mueller2024model}, considering linear problems with non--affine parametrizations. More recently, the S--OPT method was introduced in \cite{lauzon2024sopt}, and was shown to outperform DEIM in terms of accuracy for representative nonlinear model problems.
While both DEIM and S--OPT offer computational advantages, they require the evaluation of high--fidelity quantities at a selected set of mesh nodes. This introduces intrusivity at the online stage, necessitating communication between the FOM and ROM softwares. Although non--intrusive pipelines may be recovered via interpolation techniques --- as demonstrated, for instance, in \cite{chasapi2023localized} for the DEIM method --- in the present work we opted to rely exclusively on the approximate affine decompositions presented in Eqs.~\eqref{eq: convective term affine decomposition ST-ROM},~\eqref{eq: convective Jacobian affine decomposition ST-ROM}. Nevertheless, we acknowledge that a detailed comparison with alternative hyper--reduction approaches constitutes an important direction for future work. In all numerical experiments presented in Section~\ref{sec: numerical results}, we note that the increase in computational cost when retaining more temporal modes is largely independent of the chosen hyper--reduction strategy. Indeed, a code profiling analysis revealed that this phenomenon primarily stems from the resolution of dense linear systems during Newton iterations, which accounts for the majority of the computational budget when the temporal reduced basis size becomes large.
\end{remark}

\subsection{Online phase} \label{subsec: online phase}
The goal of the online phase is to efficiently compute the solution to Eq.\eqref{eq: ST-RB definition} for a parameter instance 
\begin{equation*}
\bm{\mu}^* = [\bm{\mu_f}^*, \bm{\mu_m}^*] \in \mathcal{P} \subset \R^{N_\mu^f+N_\mu^m}~, \quad \text{where} \ \ \bm{\mu_f}^* \in \mathcal{P}_f \subset \R^{N_\mu^f}, \ \bm{\mu_m}^* \in \mathcal{P}_m \subset \R^{N_\mu^m}~,
\end{equation*}
that was not considered offline. As detailed in \cite{tenderini2024space}, it comprises three steps: the assembling of the reduced system, the computation of the reduced solution and its full--order reconstruction. The first step has fundamentally already been detailed in Section~\ref{subsec: assembling the reduced system}. \\

Dealing with a nonlinear problem, the computation of the space--time--reduced solution involves the sequential resolution of some ``small'' dense linear systems, whose number depends on the convergence of the Newton iterations. This entails that another means of achieving efficiency is by reducing the number of iterations needed to satisfy the stopping criterion in Eq.\eqref{eq: Newton stopping criterion}. To this aim, the choice of the initial guess $\widehat{\bm{w}}^{(0)} \in \R^{n^{st}}$ is crucial. In fact, the convergence properties of the Newton--Raphson method indicate that the closer $\widehat{\bm{w}}^{(0)}$ is to the reduced solution $\widehat{\bm{w}}(\bm{\mu}^*)$, the more likely and the faster the iterations converge. If no dimensionality reduction in time takes place, the initial guess can be selected either as the solution at the previous timestep or by extrapolation \cite{carlberg2015decreasing}. However, these are not viable strategies with ST--RB methods, for which the choice of $\widehat{\bm{w}}^{(0)}$ becomes less straightforward. Following~\cite{choi2019space}, we decided to employ non--intrusive interpolation--based ROM approaches. So, let $\{\bm{\mu}_k\}_{k=1}^{M}$ be the offline parameter values and $\{\widehat{\bm{w}}_k\}_{k=1}^{M} \in \R^{n^{st}}$ be the space--time projections of the corresponding FOM solutions, which can be pre--computed. Then, for a given $\bm{\mu}^* \in \mathcal{P}$, the initial guess is calculated via interpolation in the parameter space. In this work, we compare two approaches:
\begin{itemize}[itemsep=1pt]
\item Nearest Neighbours Interpolation (NNI): 
\vspace{-0.2cm} 
\begin{equation*}
\label{eq: NNI formula}
\widehat{\bm{w}}^{(0)}(\bm{\mu}^*) = \sum_{k=1}^{K} \frac{\lvert\lvert\bm{\mu}^* - \bm{\mu}_{\mathcal{S}(k)}\rvert\rvert_2}{D} \widehat{\bm{w}}_{\mathcal{S}(k)} \quad \text{with} \ \ D := \sum_{k=1}^{K} \lvert\lvert\bm{\mu}^* - \bm{\mu}_{\mathcal{S}(k)}\rvert\rvert_2~,
\end{equation*}
where $\mathcal{S}(k)$ denotes the index of the $k$--th closest parameter to $\bm{\mu}^*$ within $\{\bm{\mu}_k\}_{k=1}^M$ in Euclidean norm. The quantity $k \in \mathbb{N}$ is the number of considered nearest neighbours in the parameters' space.
\item Proper Orthogonal Decomposition Interpolation (PODI) \cite{bui2004aerodynamic, girfoglio2022non}: 
\vspace{-0.2cm} 
\begin{equation}
\label{eq: PODI formula}
\begin{split}
\left(\widehat{\bm{w}}^{(0)}(\bm{\mu}^*)\right)_\ell = &\sum_{k=1}^{M} \bm{\vartheta}^{\ell}_k \zeta\left(\lvert\lvert\bm{\mu}^* - \bm{\mu}_k\rvert\rvert_2\right) \\
\qquad &\text{with} \ \ \bm{\vartheta}^{\ell} \in \R^{M} \ : \ \bm{Z}\bm{\vartheta}^\ell = \bm{a}^\ell, \quad 
\big(\bm{Z}\big)_{k_1k_2} = \zeta\left(\lvert\lvert\bm{\mu}_{k_1} - \bm{\mu}_{k_2}\rvert\rvert_2\right), \ \left(\bm{a}^\ell\right)_{k} = \left(\widehat{\bm{w}}_k\right)_\ell.
\end{split}
\end{equation}
Here $\zeta: \R^+ \to \R$ is a kernel function; in this work we relied on thin--plate spline interpolation, hence choosing $\zeta(x) = x^2 \ln(x)$. We underline that this approach entails solving, during the offline phase, $n^{st}$ independent linear systems of dimension $M$. Since the left--hand side matrix $\bm{Z}$ is the same for all the systems, convenient factorizations can be exploited to speedup the computations.
\end{itemize}
Section~\ref{subsec: navier-stokes test} features an empirical comparison of the performances of NNI and PODI; the \emph{naive} choices of selecting $\widehat{\bm{w}}^{(0)}$ either as the zero vector or as the average  of the offline reduced solutions are also considered.

\begin{remark}
	To further reduce the online computational costs, one can neglect the reduced Jacobian of the convective term $\widehat{\bm{J}}_{\bm{c}}(\widehat{\bm{w}})$ in Eq.\eqref{eq: reduced Jacobian matrix}, thus resorting to a quasi--Newton method. In the laminar flow regime, where viscous forces dominate over inertial effects, this should not significantly impact on the accuracy of the method, since the change in the reduced Jacobian $\widehat{\bm{J}}_{\widehat{\bm{r}}}$ is expected to be small. Nonetheless, this technique is particularly convenient because it keeps the Jacobian constant across the Newton iterations, hence enabling the use of LU factorization. The largest computational gains can be realized when the parametrization is confined to the right--hand side vector; an example is represented by test case reported in Section~\ref{subsec: navier-stokes test}. Indeed, in such cases the LU factorization of the reduced Jacobian can be pre--computed, and only triangular systems are left to be solved online.
\end{remark}

Finally, the FOM expansion $\widehat{\bm{w}}^{st}(\bm{\mu}^*)$ of the reduced solution $\widehat{\bm{w}}(\bm{\mu}^*)$ can be computed as
\begin{equation}
\label{eq: FOM expansion}
\widehat{\bm{w}}^{st}(\bm{\mu}^*) = \bm{\Phi} \ \widehat{\bm{w}}^M(\bm{\mu}^*) \ \bm{\Psi}^T \quad \in \R^{N^s \times N^t}~,
\end{equation}
where $\widehat{\bm{w}}^M(\bm{\mu}^*) \in \R^{n^s \times n^t}$ is a suitable two--dimensional reshaping of $\widehat{\bm{w}}(\bm{\mu}^*)$, while $\bm{\Phi} \in \R^{N^s \times n^s}$ and $\bm{\Psi} \in \R^{N^t \times n^t}$ are block--diagonal matrices, whose blocks contain the velocity, pressure, and Lagrange multipliers reduced bases in space and in time, respectively. Here we define $n^s := n_u^s + n_p^s + N_\lambda$ and $n^t := n_u^t + n_p^t + n_\lambda^t$. We remark that in Eq.\eqref{eq: FOM expansion} $\widehat{\bm{w}}^{st}(\bm{\mu}^*)$ is considered a two--dimensional matrix for convenience; however, it should be conveniently flattened into a $N^{st}$--dimensional array to be compatible with Eq.\eqref{eq: monolithic FOM system concise}.

\section{Numerical results} \label{sec: numerical results}

In this section, we analyse the performances of the ST--GRB method on three different test cases. The ``standard'' RB method (denoted as SRB--TFO) serves as a baseline. All simulations were performed on the \emph{Jed} cluster of the \emph{Scientific IT and Application Support (SCITAS)}\footnote{\url{https://www.epfl.ch/research/facilities/scitas/hardware/}} at EPFL.

\subsection{Test cases setup} \label{subsec: test cases setup}
\begin{figure}[t!]
\centering
\includegraphics[width=0.99\textwidth]{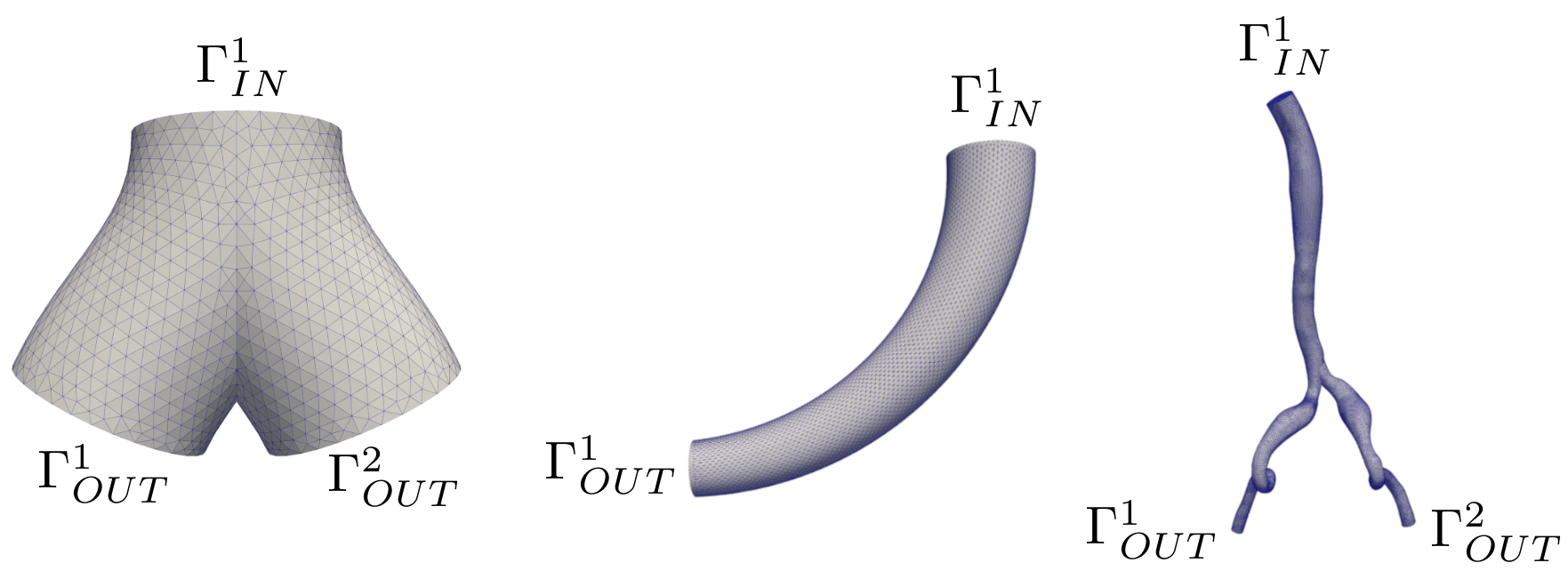}
\caption{Geometries and computational meshes used in the three numerical experiments.}
\label{fig: geometries}
\end{figure}

We solve the unsteady incompressible Navier--Stokes equations in three--dimensional vessel anatomies. In the first test case, we neglect the effect of wall compliance, enforcing no--slip boundary conditions on the lateral wall of an idealized symmetric bifurcation geometry. In the other two test cases, instead, we focus on the physics--based reduction of the FSI problem through the coupled momentum model, considering either an elementary cylindrical geometry or the anatomy of a descending aorta with external iliac arteries~\cite{wilson2013vascular}. The three geometries and the associated computational meshes are reported in Figure~\ref{fig: geometries}. In this work, we always adopt the \emph{cgs} (centimeter--gram--second) unit system. So, for instance, the velocities are expressed in $\si{\centi\metre~\per\second}$ and the pressures in $\si{\dyne~\per\centi\metre\squared}$, where $\si{\dyne} := \si{\gram\cdot\centi\metre~\per\square\second}$. \\

\begin{figure}[t!]
\centering
\includegraphics[width=0.475\textwidth]{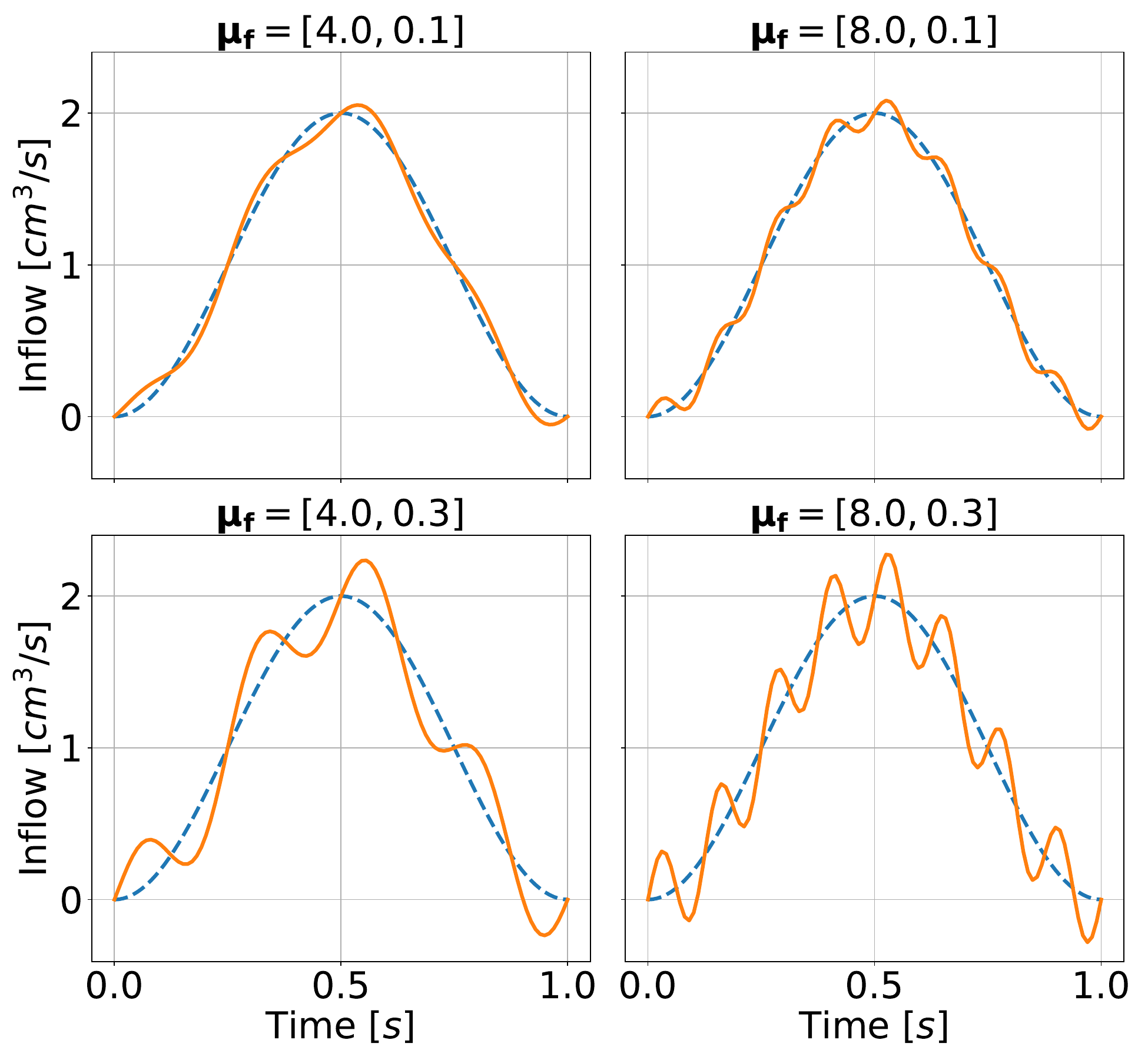}
\hspace{.03\textwidth}
\includegraphics[width=0.475\textwidth]{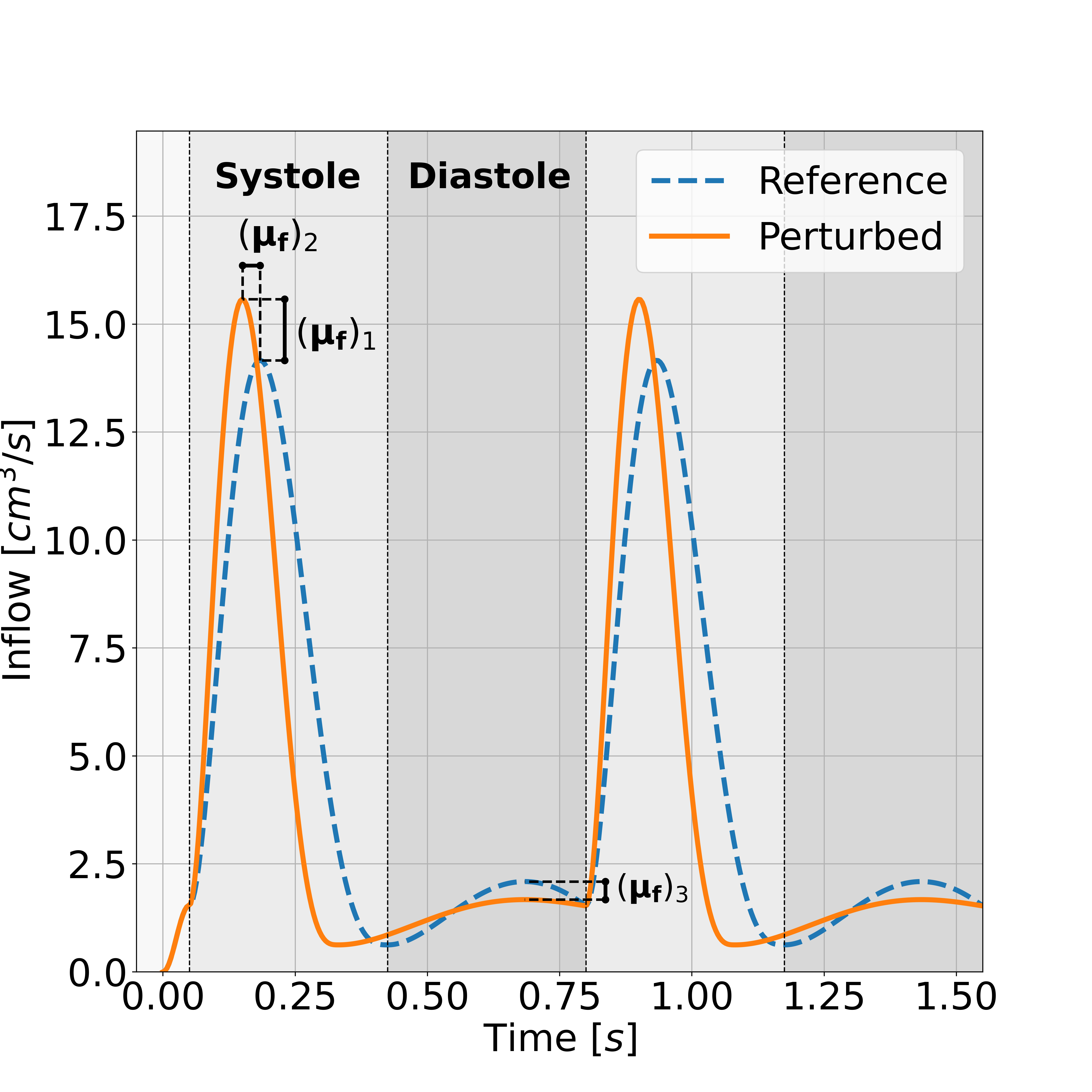}
\caption{Volumetric blood flow waveforms considered in the three numerical experiments. In particular: (left) perturbed sinusoidal waveform, used in test cases 1 and 2; (right) parametrized supraceliac aortic waveform, obtained from patient--specific PC--MRI measurements.}
\label{fig: flow rates}
\end{figure}

We set blood density to $\rho = 1.06 \ \si{\gram~\per\centi\metre\cubed}$ and blood viscosity to $\mu = 3.5 \cdot 10^{-3} \ \si{\gram~\per{(\centi\metre \cdot \second)}}$. At the inlet/outlet boundaries, depending on the test case, we either prescribe a time--varying parabolic velocity profile that matches a given volumetric flow rate or we impose homogeneous Neumann boundary conditions. As discussed in Section~\ref{subsec: algebraic formulation}, Dirichlet boundary conditions on $\Gamma_D$ are weakly imposed using Lagrange multipliers, whose space is discretized by orthonormal basis functions, built from Chebyshev polynomials \cite{pegolotti2021model}. We selected a maximal polynomial degree of $5$ for the inlets, to guarantee a good approximation of the velocity profiles, and of $0$ for the outlets, where we are just interested in enforcing the flow rate. More in detail, we work under the assumption that all the inlet/outlet Dirichlet boundaries $\{\Gamma_{\operatorname{D}}^k\}_{k=1}^{N_D}$ are planar circular surfaces. So, let $r_k \in \R^+$, $\bm{c_k} \in \R^3$, $\bm{n_k} \in \R^3$ be the radius, the center and the outward unit normal vector of $\Gamma_{\operatorname{D}}^k$, respectively. The parabolic velocity profile, that coincides with the space--dependent factor in Eq.\eqref{eq: dirichlet datum}, is then defined as follows:
\begin{equation}
\label{eq: parabolic profile}
\vec{g}^s_k(\vec{x}) := \pm \frac{2}{\pi r_k^2} \left( 1 - \frac{\left\Vert \vec{x} - \bm{c_k} \right\Vert_2^2}{r_k^2} \right) \bm{n_k} \qquad \text{with} \ \vec{x} \in \Gamma_{\operatorname{D}}^k~,
\end{equation}
where the leading sign is negative at the inlets and positive at the outlets. We remark that, since $\vec{g}^s_k$ integrates to $1$ on $\Gamma_{\operatorname{D}}^k$, the term $g_t^k(t; \bm{\mu_f})$ in Eq.\eqref{eq: dirichlet datum} coincides with the flow rate prescribed at $\Gamma_{\operatorname{D}}^k$. The expressions of $\{g_t^k(t; \bm{\mu})\}_{k=1}^{N_D}$ will be defined separately for each test case.\\

During the offline phase, we generate a training dataset of $M$ full--order solutions, that correspond to different parameter values $\{\bm{\mu}_i\}_{i=1}^{M}$, sampled uniformly at random from $\mathcal{P} := \mathcal{P}_f \times \mathcal{P}_m$. The definitions of $\mathcal{P}_f$ (flow rate parameters space) and $\mathcal{P}_m$ (membrane parameters space) will be detailed separately for the three test cases. The high--fidelity snapshots are computed with \emph{LifeV}, a C++ FE library with support for high--performance computing~\cite{bertagna2017lifev}, exploiting parallelization over $72$ cores. For what concerns the discretization in space, we employ either $\mathbb{P}_2-\mathbb{P}_1$ Taylor--Hood Lagrangian finite elements (test cases 1 and 2) or $\mathbb{P}_1^b-\mathbb{P}_1$ Crouzeix--Raviart elements (test case 3) for the discretization of velocity and pressure. Time integration is performed using the implicit BDF2 method, a second--order multistep scheme with coefficients $\beta=\frac23$, $\alpha_1=-\frac43$, $\alpha_2=\frac13$. The FOM sparse linear systems are solved using the preconditioned GMRES method~\cite{saad1986gmres}, considering the saddle point block preconditioner proposed in~\cite{pegolotti2021model}. Concerning the Newton iterations, we select a tolerance $\tau_{NR}=10^{-5}$ and a maximal number of iterations $K_{NR}=10$. \\

The performances of the reduced order models are assessed considering additional $M^*$ high--fidelity simulations, corresponding to parameter values $\{\bm{\mu}_i^*\}_{i=1}^{M^*}$ that have not been used for the reduced bases generation. The accuracy of the methods is quantified by computing the average relative errors on velocity, pressure, and displacement in the natural spatio--temporal norms. For $i = 1, \cdots, M^*$, let $\widehat{\bm{u}}_i^* \in \R^{n_u^{st}}, \ \widehat{\bm{p}}_i^* \in \R^{n_p^{st}}$ be the space--time--reduced velocity and pressure obtained with ST--GRB for the parameter value $\bm{\mu}_i^*$ and let $\bm{u}_i^* \in \R^{N_u^{st}}, \ \bm{p}_i^* \in \R^{N_p^{st}}, \bm{d}_i^* \in \R^{N_d^{st}}$ denote the corresponding space--time high--fidelity solution. Here $\bm{d}_i^*$ represents the full--order displacement field and we define $N_d^{st} := N_d^s N^t$, where $N_d^s < N_u^s$ is the number of FOM DOFs that lie on the lateral wall $\Gamma$. Also, let $\bm{X}_u^{st} \in \R^{N_u^{st} \times N_u^{st}}$, $\bm{X}_p^{st} \in \R^{N_p^{st} \times N_p^{st}}$, $\bm{X}_d^{st} \in \R^{N_d^{st} \times N_d^{st}}$ be the matrices that encode the natural spatio--temporal norms for velocity, pressure, and displacement, respectively. We refer to~\cite{tenderini2024space} for the precise definitions of $\bm{X}_u^{st}$ and $\bm{X}_p^{st}$. The matrix $\bm{X}_d^{st}$, instead, writes as follows:
\begin{equation*}
\label{eq: displacement norm matrix}
\bm{X}_d^{st} := \operatorname{diag} \Bigl( \underbrace{\bm{X}_d, \ \cdots, \bm{X}_d}_{N^t} \Bigr) \quad \text{with} \ \ \bm{X}_d \in \R^{N_d^s \times N_d^s} \ : \ \ \left(\bm{X}_d\right)_{ij} = \int_\Gamma \vec{\varphi}_i^u \cdot \vec{\varphi}_j^u~.
\end{equation*}
The average relative errors are then
\begin{equation}
\label{eq: error metrics}
E_u = \frac{1}{M^*}\sum_{i=1}^{M^*} \frac{\lVert \bm{\Pi}^u \widehat{\bm{u}}_i^* - \bm{u}_i^* \rVert_{\bm{X}_u^{st}}}{\lVert \bm{u}_i^* \rVert_{\bm{X}_u^{st}}}~; \quad 
E_p = \frac{1}{M^*}\sum_{i=1}^{M^*} \frac{\lVert \bm{\Pi}^p \widehat{\bm{p}}_i^* - \bm{p}_i^* \rVert_{\bm{X}_p^{st}}}{\lVert \bm{p}_i^* \rVert_{\bm{X}_p^{st}}}~; \quad 
E_d = \frac{1}{M^*}\sum_{i=1}^{M^*} \frac{\lVert \bm{\Pi}^d \widehat{\bm{u}}_i^* - \bm{d}_i^* \rVert_{\bm{X}_d^{st}}}{\lVert \bm{d}_i^* \rVert_{\bm{X}_d^{st}}}~. 
\end{equation}
To guarantee a fair comparison of the performances across different choices of the tolerances in the truncated POD algorithm, we express accuracy by means of the ratios $E_u / \varepsilon_u$, $E_p / \varepsilon_p$, $E_d / \varepsilon_u$, where $\varepsilon_u, \varepsilon_p \in \R^+$ are the user--defined POD tolerances for velocity and pressure (in both space and time), respectively. In the following, we will refer to such quantities as the \emph{normalized average relative errors}. Incidentally, we define $\varepsilon_\lambda^s, \ \varepsilon_\lambda^t \in \R^+$ as the POD tolerances for the Lagrange multipliers, respectively, in space and in time. To quantify the computational efficiency, instead, we consider two different indicators: the \emph{speedup} (SU), defined as the ratio between the average duration of FOM and ROM simulations, and the \emph{reduction factor} (RF), computed as the ratio between the number of FOM and ROM DOFs. \\

\begin{table}[!t]
\centering 
\caption{Average mesh size ($h_{\operatorname{avg}}$), timestep size, number of space and time FOM DOFs for velocity, pressure, and Lagrange multipliers, and average duration of a FOM simulation (on $72$ cores) in the three test cases.}
\label{tab: fem info}
\resizebox{.95\textwidth}{!}{
\begin{tabular}{cccccccccc}
	\toprule
	
	\multicolumn{3}{c}{\textbf{Test case}}
	&\multicolumn{4}{c}{\textbf{Space}}
	&\multicolumn{2}{c}{\textbf{Time}} & \\
	
	\cmidrule(lr){1-3} \cmidrule(lr){4-7} \cmidrule(lr){8-9} 
	
	\textbf{Test \#}
	&\textbf{Problem}
	&\textbf{Geometry}
	&\textbf{$\bm{h_{\operatorname{avg}}}$}
	&\textbf{$\bm{N_u^s}$}
	&\textbf{$\bm{N_p^s}$}
	&\textbf{$\bm{N_\lambda}$}
	&\textbf{$\bm{\Delta t}$}
	&\textbf{$\bm{N^t}$}
	&\textbf{Duration} \\
	
	\midrule
	
	\textbf{1} & \textbf{Navier--Stokes} & \textbf{Bifurcation} & 1.35 \si{\milli\meter}\  & \phantom{0}76'974  & \phantom{0}3'552 & 66 & 1.0 \si{\milli\second} & 1'000 & \phantom{0}1\si{\hour} 44\si{\min} \\
	
	\textbf{2} &\textbf{RFSI} & \textbf{Tube} & 2.18 \si{\milli\meter} & \phantom{0}86'802  & \phantom{0}4'350  & 63 & 1.0 \si{\milli\second} & 1'000 & \phantom{0}1\si{\hour} 47\si{\min} \\
	
	\textbf{3} &\textbf{RFSI} & \textbf{Aorta} & 2.21 \si{\milli\meter} & 394'527 & 24'687 & 63 & 1.0 \si{\milli\second} & 1'550 & 10\si{\hour} 24\si{\min} \\
	
	\bottomrule 
\end{tabular}
}
\end{table}

\begin{table}[t!]
\caption{Offline computational costs, expressed in terms of CPU--time (in \si{\second}), in the three test cases, for $\varepsilon_u=\varepsilon_p=\varepsilon_\lambda^t=10^{-3}$ and $\varepsilon_\lambda^s=10^{-5}$.}
\label{tab: offline phase info}
\centering 
\resizebox{.85\textwidth}{!}{
\begin{tabular}{cccccccccc}
	\toprule
	
	&&& \multicolumn{4}{c}{\textbf{Bases construction}} & & & \\
	
	\cmidrule(lr){4-7}  
	
	\multicolumn{1}{c}{\textbf{Test \#}}
	&\multicolumn{1}{c}{$\bm{M}$}
	&\multicolumn{1}{c}{\textbf{Method}}
	&\multicolumn{2}{c}{$\bm{u}$}
	&\multicolumn{2}{c}{$\bm{p}$}
	&\multicolumn{3}{c}{\textbf{Assembling}} \\
	
	\cmidrule(lr){4-5} \cmidrule(lr){6-7} \cmidrule(lr){8-10}
	
	&&& \textbf{Space} & \textbf{Time} &\textbf{Space} & \textbf{Time} & \textbf{\emph{IG}} & \textbf{\emph{Conv}} \textsuperscript{$(\bm{n_c})$} & \textbf{\emph{Other}} \\
	
	\midrule
	
	\multirow{2}{*}{\textbf{1}} & \multirow{2}{*}{$\bm{50}$} & \textbf{SRB--TFO} &
	\multirow{2}{*} {550 \si{\second}} & // & \multirow{2}{*} {20 \si{\second}} & // & // & \multirow{2}{*} {231 \si{\second} \textsuperscript{(39)}} & 2 \si{\second} \\ 
	
	& &\textbf{ST--GRB} &
	& 1 \si{\second} & & < 1 \si{\second} &  38 \si{\second} &  & 2 \si{\second} \\ 
	
	\midrule
	
	\multirow{2}{*}{\textbf{2}} & \multirow{2}{*}{$\bm{50}$} & \textbf{SRB--TFO} &
	\multirow{2}{*} {547 \si{\second}} & // &  \multirow{2}{*} {37 \si{\second}} &  // & // & \multirow{2}{*} {65 \si{\second}} \textsuperscript{(10)} & 3 \si{\second} \\
	
	& & \textbf{ST--GRB} &
	& 1 \si{\second} &  & < 1 \si{\second} & 32 \si{\second} &  & 3 \si{\second} \\  
	
	\midrule
	
	\multirow{2}{*}{\textbf{3}} & \multirow{2}{*}{$\bm{25}$} &  \textbf{SRB--TFO} &
	\multirow{2}{*} {497 \si{\second}} & // & \multirow{2}{*} {75 \si{\second}} &  // & // & \multirow{2}{*} {2'928 \si{\second} \textsuperscript{(93)}} & 9 \si{\second} \\
	
	& &  \textbf{ST--GRB} &
	& 6 \si{\second} &  & < 1  \si{\second} &   162 \si{\second} &  & 12 \si{\second} \\ 
	
	\bottomrule 
\end{tabular}
}
\end{table}

Table~\ref{tab: fem info} reports the average mesh size, the timestep size, the number of FOM DOFs in space and in time and the average duration of a single high--fidelity simulation (performed in parallel over $72$ cores) for the three test cases. By comparing test cases 1 and 2, that feature a similar number of DOFs, we remark that the use of the coupled momentum model --- denoted for convenience as RFSI --- does not entail important increases in the computational cost. \\

Table~\ref{tab: offline phase info} summarizes the offline costs of SRB--TFO and ST--GRB, using a POD tolerance of $\varepsilon=10^{-3}$ for all the fields, both in space and in time. Three remarks are worth to follow. 
Firstly, as also observed in~\cite{tenderini2024space}, all the wall--times reported in Table~\ref{tab: offline phase info} are a small fraction of a single FOM solve. This reveals that the high--fidelity data generation (or availability) is the true bottleneck of the offline phase.
Neglecting the construction of the FOM solutions cohort, the computation of the velocity reduced basis in space appears as the most expensive offline operation, with a duration of roughly $10 \ \si{\minute}$ for all the test cases. Nonetheless, we remark that assembling the affine components for the convective term and its jacobian (see Section~\ref{subsubsec: ST-reduced convective term}~) also entails an important cost. The associated wall--time is reported in the column \emph{Conv}, considering $n_c$ equal to the number of velocity basis functions prior to the supremizers' enrichment; such value is reported as a superscript. Notably, the affine components assembly takes approximately six times longer than the reduced bases generation for the aorta test case. However, it is worth mentioning that parallelization was not exploited at this stage, despite the task at hand being embarrassingly parallelizable.
Compared to SRB--TFO, the computational cost of ST--GRB is increased by the use of interpolation--based techniques to determine a good initial guess for the Newton iterations. Indeed, this step requires the offline projection of the spatio--temporal FOM snapshots onto the reduced subspace. This operation carries a non--negligible computational burden (see column \emph{IG}), that scales linearly with the number of snapshots $M$.
Lastly, we underline that the offline cost of model order reduction in space is largely dominant not only for the reduced bases generation, but also for the assembling of parameter--independent quantities (see Section~\ref{subsubsec: ST-reduced matrices}); the associated wall--time is reported in the column \emph{Other}. In fact, the offline assembling phase of ST--GRB is substantially as expensive as the one of SRB--TFO, for all the test cases. As a result, the calculation of the generalized coordinates of the FOM snapshots configures as the primary differentiating factor in efficiency between the offline phases of SRB--TFO and ST--GRB.

\subsection{Test case 1: rigid Navier--Stokes in a symmetric bifurcation} \label{subsec: navier-stokes test}
We investigate the performances of ST--GRB in solving the Navier--Stokes equations in an idealized symmetric bifurcation geometry (see Figure~\ref{fig: geometries} -- left). At the inlet $\Gamma_{\operatorname{IN}}^1$, we impose a parabolic velocity profile (see Eq.\eqref{eq: parabolic profile}) that matches the following parametrized volumetric flow rate (see Figure~\ref{fig: flow rates} -- left):
\begin{equation}
\label{eq: flow rate 1}
g_1^t(t; \bm{\mu_f}) = 1 - \cos \left(\frac{2 \pi t}{T}\right) + (\bm{\mu_f})_2 \sin \left(\frac{2 \pi (\bm{\mu_f})_1 t}{T}\right),
\end{equation} 
where $T = 1 \ \si{\second}$ is the final time of the simulation. To enhance variability within the solution manifold, we unbalance the flow in the two outlet branches by (weakly) enforcing the outflow rate $g_2^t(t; \bm{\mu_f}) := (\bm{\mu_f})_3 \ g_1^t(t; \bm{\mu_f})$ on $\Gamma_{\operatorname{OUT}}^1$. Hence, we have that $\bm{\mu_f} \in \R^3$; we define the space of flow rate parameters as $\mathcal{P}_f = [4.0,8.0] \times [0.1, 0.3] \times [0.2, 0.8]$. Incidentally, we note that the first two parameters mostly affect the temporal diversity of the solution manifold, while the last parameter is related to spatial diversity.
On $\Gamma_{\operatorname{OUT}}^2$ we impose homogeneous Neumann boundary conditions. We underline that in this test case we neglect wall compliance; so, $\mathcal{P}_m = \emptyset$ and no--slip boundary conditions are strongly enforced at the lateral wall $\Gamma$. Unless otherwise specified, in the numerical tests we consider $M=50$ snapshots, we choose the POD tolerances $\varepsilon_u=\varepsilon_p=\varepsilon_\lambda^t=10^{-3}$, $\varepsilon_\lambda^s = 10^{-5}$ and $n_{c,J}=0$ affine components for the Jacobian of the convective term, and we compute the initial guess $\widehat{\bm{w}}^{(0)}$ through PODI (see Eq.\eqref{eq: PODI formula}). \\

\begin{figure}[t!]
\centering
\includegraphics[width=0.495\textwidth]{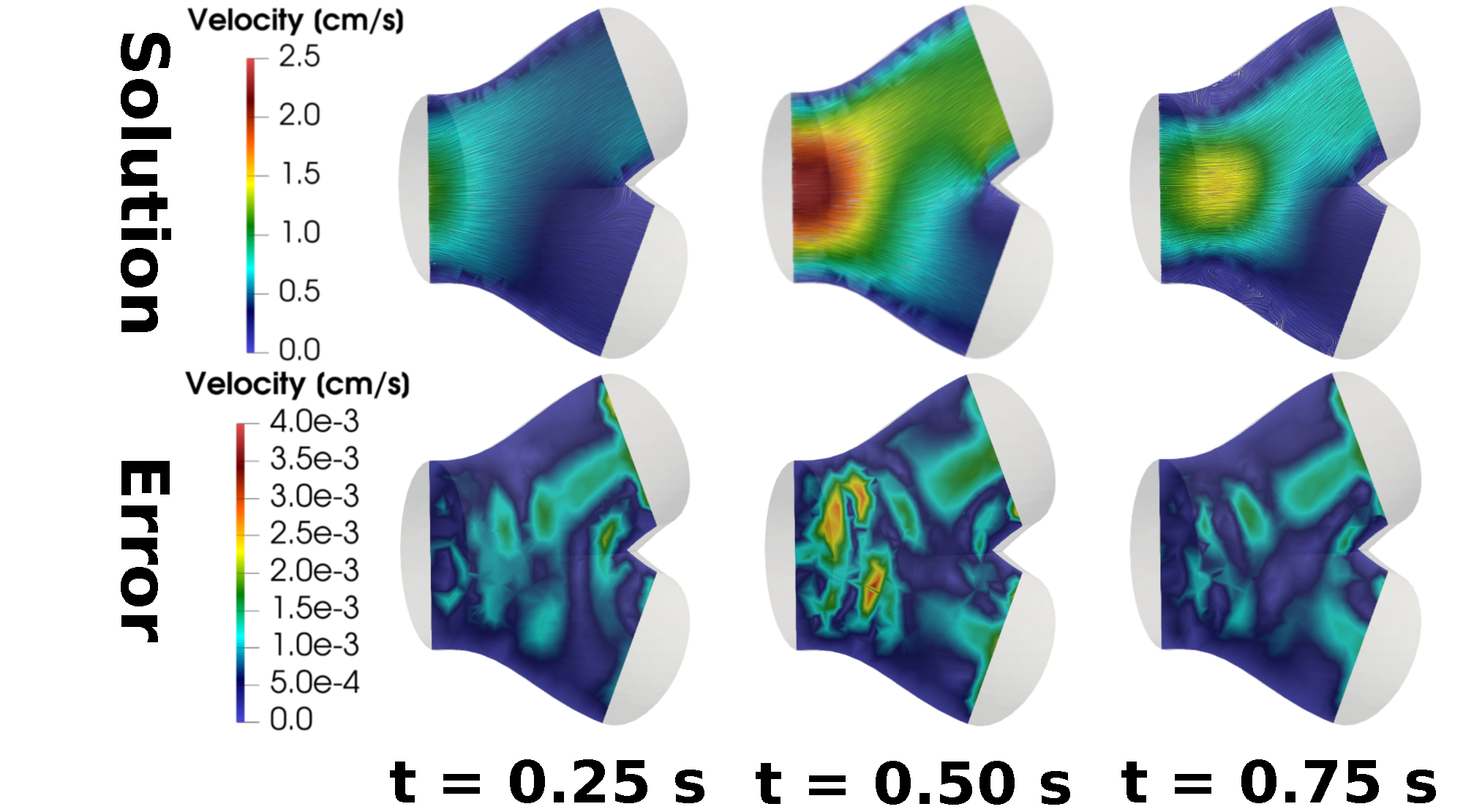}
\includegraphics[width=0.495\textwidth]{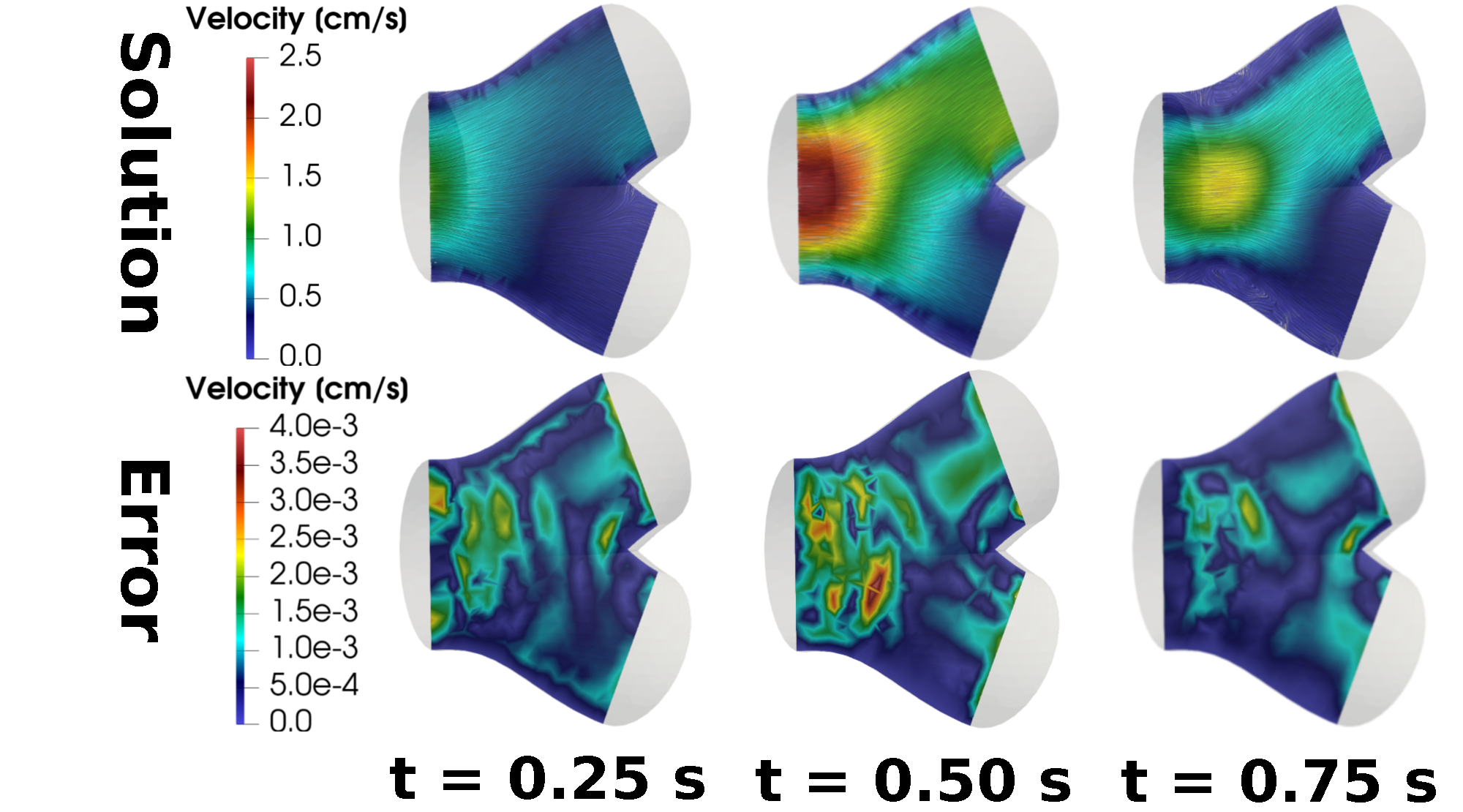}
\caption{Line integral convolution of the velocity field on the median slice (top) and corresponding pointwise absolute error (bottom) at three time instants, obtained with SRB-TFO (left) and ST--GRB right for $\bm{\mu_f}^* = [7.56, 0.14, 0.74]$.}
\label{fig: NS solution}
\end{figure}

Figure~\ref{fig: NS solution} shows the line integral convolution of the velocity field on the median slice of the domain obtained for $\bm{\mu_f}^* = [7.56, 0.14, 0.74]$ with SRB--TFO and ST--GRB, and the corresponding pointwise error with respect to the FOM solution. We underline that this parameter value was not considered to generate the reduced basis during the offline phase of the method. The velocity relative errors (see Eq.\eqref{eq: error metrics}) for this particular experiment are $0.28 \%$ for SRB--TFO and $0.61 \%$ for ST--GRB. From a qualitative standpoint, the two solutions are almost undistinguishable at the three displayed time instants. Notably, both SRB--TFO and ST--GRB are able to capture the flow recirculation phenomenon that occurs after the flow peak, during the deceleration phase. \\

\begin{figure}[t!]
\centering
\includegraphics[width=0.95\textwidth]{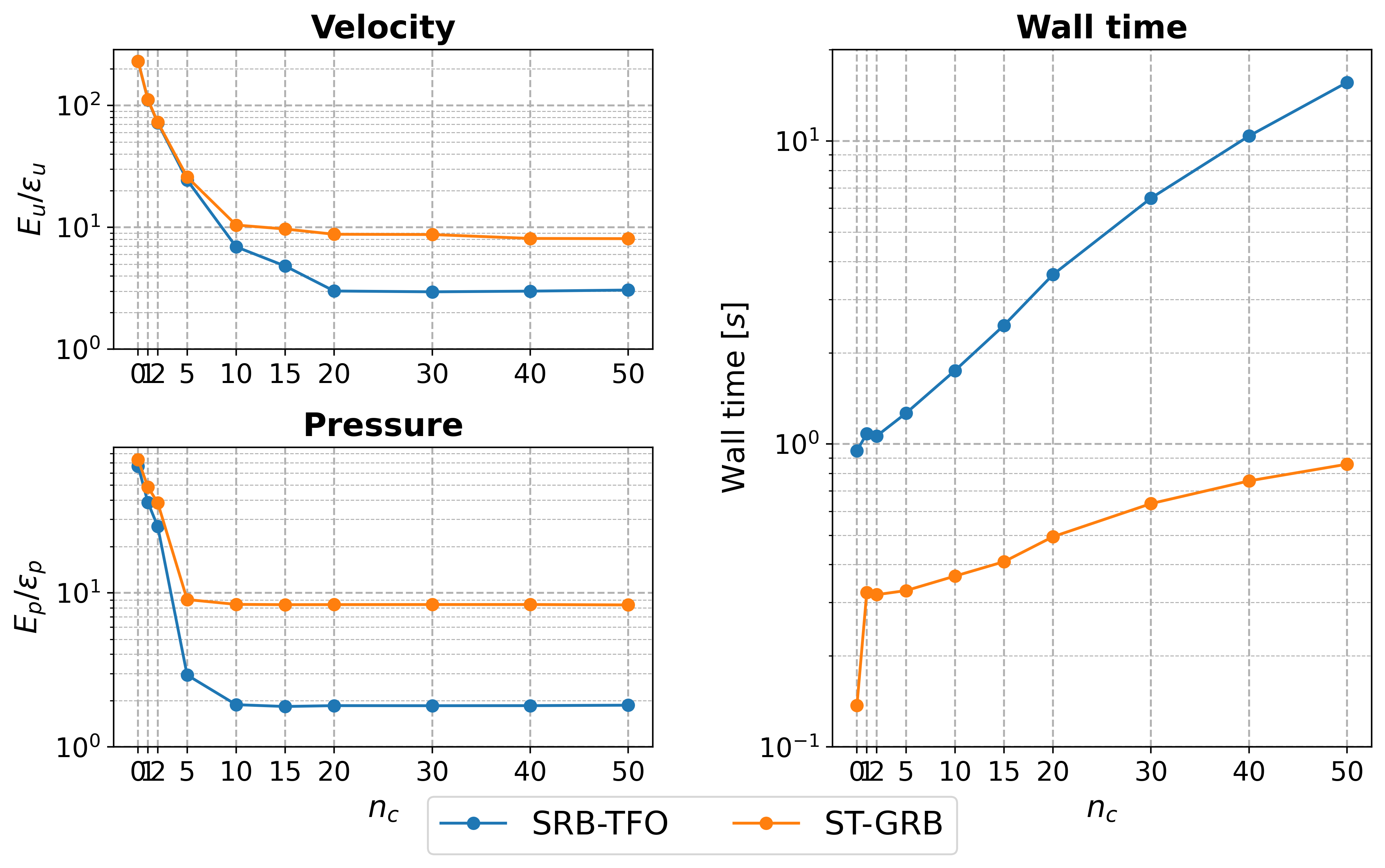}
\caption{Trend of velocity and pressure normalized average test relative errors (left) and speedups (right) with respect to $n_c$, the number of affine components for the nonlinear reduced convective term.}
\label{fig: GS errors NS}
\end{figure}

Figure~\ref{fig: GS errors NS} shows the trend of the normalized average relative test errors and of the computational wall times with respect to the number of affine components $n_c$, considered for the approximation of the convective term (see Eq.\eqref{eq: convective term approximate affine decomposition ROM}). Both SRB--TFO and ST--GRB exhibit a decreasing error trend for small values of $n_c$, while a plateau is reached once a sufficiently large number of convective modes is included.
On the one side, this behaviour tells that an imprecise approximation of the nonlinearity severely deteriorates the solution quality. However, on the other side, it also highlights that good quality solutions can be obtained considering a relatively small number of affine components. More in detail, we notice that SRB--TFO is always more accurate than ST--GRB, and that the plateau in the error is reached later with SRB--TFO and for the velocity field. These behaviours can be explained considering that ST--GRB introduces an additional reduction layer compared to SRB--TFO and that the definition of the convective term does not explicitly depend on pressure. 
Notably, neither SRB--TFO nor ST--GRB achieve superior accuracy levels when increasing $n_c$ from 50 to $n_c = n_u^s = 74$, at which point our approach exactly coincides with the non--hyper--reduced one. 
Concerning efficiency, we observe that the wall--time is strongly related to $n_c$ with SRB--TFO, while such dependency is weaker with ST--GRB. This result can be justified considering the impact of the convective term assembly on the overall online phase duration. Indeed, with ST--GRB, the most expensive online operation consists in solving dense linear systems of dimension $n^{st} \approx 10^3$, and assembling the convective term (once per Newton iteration) carries a smaller cost, at least for sufficiently small values of $n_c$. Conversely, with SRB--TFO, the nonlinear term assembly cost (once per Newton iteration at every timestep) eclipses the one of solving dense linear systems of dimension $n^s \approx 10^2$, which explains the observed markedly increasing wall--time trend. Ultimately, the value $n_c=20$ appears to provide the best trade--off between accuracy and efficiency for both methods. Therefore, we used it in all the numerical experiments reported for this test case. \\

As pointed out in Section~\ref{subsec: online phase}, the choice of the initial guess for the Newton iterations is not straightforward in the context of space--time model order reduction. For this reason, we investigate six different alternatives. The results of the numerical tests, conducted both neglecting and including the convective Jacobian modes, are summarized in Table~\ref{tab:bifurcation ST--RB IG}, and demonstrate that ST--GRB benefits from a smart choice of the initial guess.
Indeed, naively selecting $\widehat{\bm{w}}^{(0)} = \bm{0}$ (\emph{Zero}) yields rather imprecise approximations, particularly for the velocity field. Conversely, fast convergence and accurate solutions are attained when employing non--intrusive ROMs, such as \emph{K--NN} (with $K = 1, 3, 5$) and \emph{PODI}. In particular, for both $n_{c,J}=0$ and $n_{c,J}=20$, the best results are obtained with PODI, even though $3$--NN also showcases similar performances. Three remarks are worth following. Firstly, selecting $\widehat{\bm{w}}^{(0)} = \frac{1}{M} \ \sum_{i=1}^{M} \widehat{\bm{w}}(\bm{\mu}_i)$ (\emph{Average}) leads to convergent Newton iterations and low errors (at least for $n_{c,J}=20$), despite this choice is not tailored to the value of $\bm{\mu}^*$. Secondly, taking into account the convective Jacobian modes enhances accuracy at the cost of larger computational times, even if the average number of Newton iterations is lowered. In this regard, we underline that LU factorization of the left--hand side matrix was not exploited in these tests, to guarantee a fair comparison of the different alternatives. Finally, we note that both the number of Newton iterations and the online times--to--solution are similar for all the considered initial guesses, except for the case $\widehat{\bm{w}}^{(0)} = \bm{0}$ with $n_{c,J} = 0$. This result indicates that the online query of the non--intrusive ROMs entails negligible computational costs. \\

\begin{table}[t!]
\centering 
\caption{Normalized average relative test errors on the solution and on the initial guess, average number of Newton iterations and wall times. The results have been obtained considering the ST--GRB method with different choices of the initial guess in the Newton iterations and different numbers of affine components for the convective Jacobian. The quantities $E_u^{(0)}$, $E_p^{(0)}$ denote the average relative test error on velocity and pressure, considering the initial guess as reduced solution.}
\label{tab:bifurcation ST--RB IG}
\resizebox{.9\textwidth}{!}{
\begin{tabular}{ccccccccccccc} 
	\toprule
	& \multicolumn{2}{c}{\large\textbf{Zero}} & \multicolumn{2}{c}{\large\textbf{Average}} & \multicolumn{2}{c}{\large\textbf{1--NN}} & \multicolumn{2}{c}{\large\textbf{3--NN}} & \multicolumn{2}{c}{\large\textbf{5--NN}} & \multicolumn{2}{c}{\large\textbf{PODI}} \\
	
	\cmidrule(lr){2-3} \cmidrule(lr){4-5} \cmidrule(lr){6-7} \cmidrule(lr){8-9} \cmidrule(lr){10-11} \cmidrule(lr){12-13} 
	
	\normalsize$\bm{n_{c,J}}$ & \normalsize$\bm{0}$ & \normalsize$\bm{20}$ & \normalsize$\bm{0}$ & \normalsize$\bm{20}$ & \normalsize$\bm{0}$ & \normalsize$\bm{20}$ & \normalsize$\bm{0}$ & \normalsize$\bm{20}$ & \normalsize$\bm{0}$ & \normalsize$\bm{20}$ & \normalsize$\bm{0}$ & \normalsize$\bm{20}$ \\ 
	
	\midrule
	
	\normalsize$\bm{E_u / \varepsilon_u}$ & 40.6 & 23.2 & 9.82 & 8.23 & 9.15 & 8.22 & 8.81 & 8.21 & 8.99 & 8.21 & 8.78 & 8.20 \\
	\normalsize$\bm{E_p / \varepsilon_p}$ & 8.47 & 10.1 & 8.37 & 8.43 & 8.36 & 8.41 & 8.38 & 8.40 & 8.37 & 8.40 & 8.39 & 8.39 \\
	
	\midrule
	
	\normalsize$\bm{E_u^{(0)} / \varepsilon_u}$ & \multicolumn{2}{c}{1000} & \multicolumn{2}{c}{294} & \multicolumn{2}{c}{179} & \multicolumn{2}{c}{161} & \multicolumn{2}{c}{174} & \multicolumn{2}{c}{146} \\
	\normalsize$\bm{E_p^{(0)} / \varepsilon_p}$ & \multicolumn{2}{c}{1000} & \multicolumn{2}{c}{738} & \multicolumn{2}{c}{610} & \multicolumn{2}{c}{588} & \multicolumn{2}{c}{603} & \multicolumn{2}{c}{599} \\
	
	\midrule
	
	\normalsize\textbf{Iterations} & 7 & 2.2 & 3.8 & 2.1 & 3.6 & 2.1 & 3.6 & 2 & 3.6 & 2 & 3.5 & 2 \\
	\normalsize\textbf{Time [s]} & 0.85 & 1.13 & 0.48 & 1.25 & 0.47 & 1.27 & 0.49 & 1.14 & 0.53 & 0.94 & 0.50 & 1.03\\
	
	\bottomrule
\end{tabular}
}
\end{table}

\begin{table}[t!]
\centering 
\caption{Results obtained with SRB--TFO and ST--GRB on the Navier--Stokes symmetric bifurcation test case, for different POD tolerances $\varepsilon$. In particular: (left) size of the spatial and temporal reduced bases for velocity, pressure, and Lagrange multipliers; (center) reduction factor and average speedup; (right) normalized average relative test errors on velocity and pressure. \label{tab:bifurcation ST--RB pod tolerances}}
\resizebox{.95\textwidth}{!}{
\begin{tabular}{ccccccccc} 
	\toprule
	& &
	\multicolumn{3}{c}{\large\textbf{ROM size}}&
	\multicolumn{2}{c}{\large\textbf{Efficiency}}&
	\multicolumn{2}{c}{\large\textbf{Accuracy}}\\
	
	\cmidrule(lr){3-5} \cmidrule(lr){6-7} \cmidrule(lr){8-9}
	
	&
	\normalsize$\bm{\varepsilon}$ & 
	\normalsize$\bm{(n_u^s,n_u^t)}$ & 
	\normalsize$\bm{(n_p^s,n_p^t)}$ &
	\normalsize$\bm{(\{N_{\lambda_k}\},\{n_{\lambda_k}^t\})}$ &
	\normalsize\textbf{RF} & 
	\normalsize\textbf{SU} & 
	\normalsize$\bm{E_u / \varepsilon_u}$ & 
	\normalsize$\bm{E_p / \varepsilon_p}$ \\
	
	\midrule
	
	\multirow{5}{*}{\rotatebox[origin=c]{90}{\centering \textbf{\normalsize SRB--TFO}}} 
	& \normalsize\textbf{1e-2} & (32,$\ \cdot \ $) & (4,$\ \cdot \ $) & \{(13,$\ \cdot \ $), (2,$\ \cdot \ $)\} & 1.58e3 & 2.14e3 & 2.92 & 0.99 \\
	& \normalsize\textbf{5e-3} & (41,$\ \cdot \ $) & (4,$\ \cdot \ $) & \{(13,$\ \cdot \ $), (2,$\ \cdot \ $)\} & 1.30e3 & 2.05e3 & 5.87 & 1.17  \\
	& \normalsize\textbf{1e-3} & (74,$\ \cdot \ $)& (9,$\ \cdot \ $) & \{(23,$\ \cdot \ $), (3,$\ \cdot \ $)\} & 7.39e2 &  1.88e3 & 3.01 & 1.85 \\
	& \normalsize\textbf{5e-4} & (93, $\ \cdot \ $) & (12,$\ \cdot \ $) & \{(27,$\ \cdot \ $), (3,$\ \cdot \ $)\} & 5.97e2 & 1.54e3 & 2.87 & 1.54 \\
	& \normalsize\textbf{1e-4} & (159,$\ \cdot \ $) & (21,$\ \cdot \ $) & \{(41,$\ \cdot \ $), (3,$\ \cdot \ $)\} & 3.60e2 & 1.39e3 & 11.3 & 1.97 \\
	
	\midrule
	
	\multirow{5}{*}{\rotatebox[origin=c]{90}{\centering \textbf{\normalsize ST--GRB}}} 
	& \normalsize\textbf{1e-2} & (32,15) & (4,13) & \{(13,13), (2,12)\} & 1.11e5 & 4.50e4 & 9.62 & 6.99 \\
	& \normalsize\textbf{5e-3} & (41,16) & (4,13) & \{(15,14), (2,12)\} & 8.55e4 & 4.20e4 & 18.4 & 16.3 \\
	& \normalsize\textbf{1e-3} & (74,21) & (9,19) & \{(23,19), (3,15)\} & 3.65e4 & 1.62e4 & 8.78 & 8.39 \\
	& \normalsize\textbf{5e-4} & (93,23) & (12,21) & \{(27,22), (3,15)\} & 2.66e4 & 4.94e3 & 8.76 & 5.51 \\
	& \normalsize\textbf{1e-4} & (159,33) & (21,32) & \{(41,32), (3,21)\} & 1.10e4 & 8.96e2 & 32.4 & 13.7 \\
	
	\bottomrule 
\end{tabular}
}

\bigskip

\centering
\includegraphics[width=.95\textwidth]{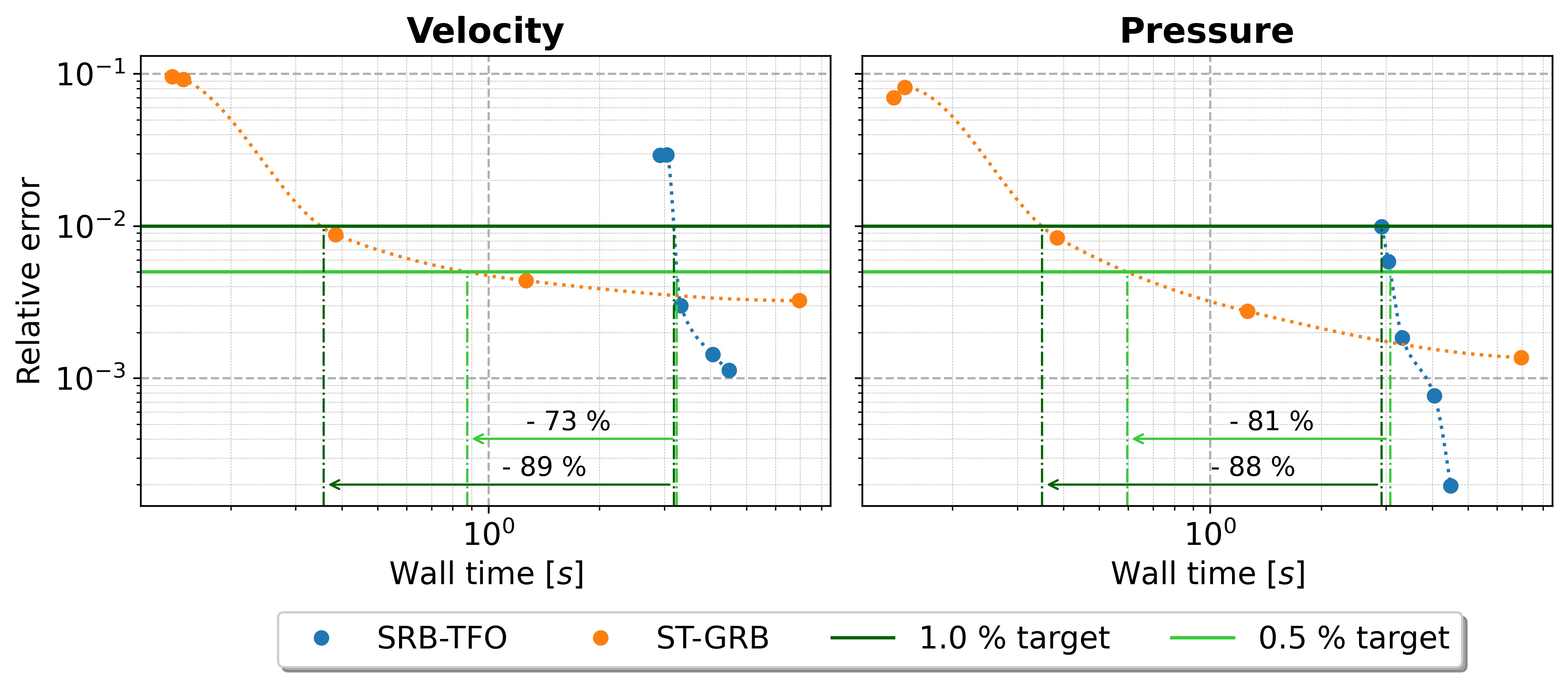}
\captionof{figure}{Average relative test error on velocity (left) and pressure (right) vs. wall time (in \si{\second}), for SRB--TFO and ST--GRB, considering the same POD toleraces reported in Table~\ref{tab:bifurcation ST--RB pod tolerances}. The dotted lines are obtained by piecewise cubic Hermite polynomial shape--preserving interpolation. To ease interpretability, the $1\%$ and $0.5\%$ accuracy targets are depicted for both fields, and the associated computational gains are indicated.}
\label{fig: performance plot NS}

\end{table}

In Table~\ref{tab:bifurcation ST--RB pod tolerances}, we report the reduced bases sizes and the performance metrics of SRB--TFO and ST--GRB, considering four different POD tolerances. The same data are also displayed in Figure~\ref{fig: performance plot NS}, where we show the average relative test errors of the two methods with respect to their computational times, and specifically investigate their behaviour to match the $1\%$ and $0.5\%$ accuracy targets.
From a general standpoint, model order reduction in time makes ST--GRB more efficient, but less accurate than SRB--TFO. For instance, in the case of $\varepsilon = 10^{-3}$, ST--GRB is roughly $10$ times faster than SRB--TFO (average simulation duration: $0.38 \ \si{\second}$ vs. $4.04 \ \si{\second}$), but its errors are approximately $3$ times larger for the velocity and $4$ for the pressure. Nonetheless, the computational advantages of ST--GRB vanish as smaller values of $\varepsilon$ are considered. In fact, for $\varepsilon=10^{-4}$, SRB--TFO is not just more accurate but also faster than ST--GRB (average simulation duration: $4.46 \ \si{\second}$ vs. $6.94 \ \si{\second}$). This behaviour is mainly related to the computational cost of solving linear systems. Indeed, solving a single $n^{st}$--dimensional (with $n^{st} = n^s n^t$) dense linear system is cheaper than sequentially solving $N^t$ distinct $n^s$--dimensional dense linear systems only if $n^t \ll N^t$. Since direct linear systems solvers have at most a cost of $\mathcal{O}(N^3)$, we can claim that ST--GRB is cheaper than SRB--TFO if $n^t \lesssim \sqrt[3]{N^t}$. For this particular test case, where $N^t=10^3$, computational gains are realized up to $n^t \approx 20$ temporal reduced bases, which corresponds to setting $\varepsilon=5\cdot 10^{-4}$. 
Figure~\ref{fig: performance plot NS} further allows interpreting the obtained results from an accuracy target perspective. Notably, we observe that ST--GRB outperforms the canonical RB method for coarse enough precision levels. Indeed, with reference to the velocity errors, it cuts the average computational time of SRB--TFO by $89\%$ and $73\%$ for the $1\%$ and $0.5\%$ accuracy targets, respectively.
Lastly, we highlight that the results reported in Table~\ref{tab:bifurcation ST--RB pod tolerances} and in Figure~\ref{fig: performance plot NS} have been obtained for $n_{c,J}=0$, i.e. neglecting the presence of the convective Jacobian in the Newton iterations. This implies that the left--hand side matrix of the linear systems to be solved is parameter--independent, since the flow rate parametrization solely affects the right--hand side vector. Therefore, a pre--computed LU factorization can be exploited to further boost the efficiency of both methods, at no loss in accuracy. Remarkably, the computational gains stemming from this strategy are more important, in relative terms, with ST--GRB than with SRB--TFO, because larger linear systems are solved under space--time model order reduction. For example, in the case of $\varepsilon=10^{-3}$, leveraging LU decomposition allows to lower the average simulation wall time of $7\%$ with SRB--TFO (from $3.31$ \si{\second} to $3.08$ \si{\second}) and of $66\%$ with ST--GRB (from $0.38$ \si{\second} to $0.13$ \si{\second}). \\

\begin{table}[t!]
\centering 
\caption{Results obtained with SRB--TFO and ST--GRB on the Navier--Stokes symmetric bifurcation test case for a different number of offline snapshots $M$, setting $\varepsilon=10^{-3}$. In particular: (left) size of the spatial and temporal reduced bases for velocity and pressure; (center) average velocity and pressure relative test errors and simulation wall times (WT) for $10$ parameters sampled within the parametric domain $\mathcal{P}_f$; (right) velocity relative test errors and average simulation wall times for $4$ parameters sampled outside the parametric domain $\mathcal{P}_f$, specifically $\bm{\bar{\mu}_1} = [\underline{2.0}, 0.2, 0.6]$, $\bm{\bar{\mu}_2} = [6.0, \underline{0.4}, 0.6]$, $\bm{\bar{\mu}_3} = [\underline{2.0}, \underline{0.4}, 0.5]$, $\bm{\bar{\mu}_4} = [6.0, 0.2, \underline{0.9}]$, where underlined values are those that lie outside of the training envelope. \label{tab:bifurcation ST--RB snapshots}}
\resizebox{.975\textwidth}{!}{
	\begin{tabular}{cccccccccccc} 
		\toprule
		& &
		\multicolumn{2}{c}{\large\textbf{ROM size}}&
		\multicolumn{3}{c}{\large\textbf{Interpolation}}&
		\multicolumn{5}{c}{\large\textbf{Extrapolation}}\\
		
		\cmidrule(lr){3-4} \cmidrule(lr){5-7} \cmidrule(lr){8-12}
		
		&
		\multirow{2}{*}{\normalsize$\bm{M}$} & 
		\multirow{2}{*}{\normalsize$\bm{(n_u^s,n_u^t)}$} & 
		\multirow{2}{*}{\normalsize$\bm{(n_p^s,n_p^t)}$} 
		&
		\multirow{2}{*}{\normalsize$\bm{E_u}$} & 
		\multirow{2}{*}{\normalsize$\bm{E_p}$} &
		\multirow{2}{*}{\normalsize\textbf{WT}}
		&
		\multicolumn{4}{c}{\normalsize$\bm{E_u}$} &
		\multirow{2}{*}{\normalsize\textbf{WT}} \\
		
		\cmidrule(lr){8-11}
		
		& & & & & & &
		\normalsize$\bm{\bar{\mu}_1}$ &
		\normalsize$\bm{\bar{\mu}_2}$ & 
		\normalsize$\bm{\bar{\mu}_3}$ &
		\normalsize$\bm{\bar{\mu}_4}$ & \\
		
		\midrule
		
		\multirow{4}{*}{\rotatebox[origin=c]{90}{\centering \textbf{\normalsize SRB--TFO}}}
		& \normalsize\textbf{5} & (59,$\ \cdot \ $) & (9,$\ \cdot \ $) & 2.39\% & 0.58\% & 3.54 \si{\second} & 2.36\% & 2.38\% & 2.44\% & 3.07\% & 3.10 \si{\second}\\ 
		& \normalsize\textbf{10} & (73,$\ \cdot \ $) & (9,$\ \cdot \ $) & 0.35\% & 0.18\% & 3.12 \si{\second} & 0.39\% & 0.32\% & 0.57\% & 0.83\% & 3.74 \si{\second} \\
		& \normalsize\textbf{20} & (73,$\ \cdot \ $) & (9,$\ \cdot \ $) & 0.31\% & 0.19\% & 3.21 \si{\second} & 0.34\% & 0.27\% & 0.51\% & 0.84\% & 3.42 \si{\second} \\
		& \normalsize\textbf{50} & (74,$\ \cdot \ $) & (9,$\ \cdot \ $) & 0.30\% & 0.19\% & 3.31 \si{\second} & 0.35\% & 0.27\% & 0.51\% & 0.92\% & 3.24 \si{\second} \\
		
		\midrule
		
		\multirow{4}{*}{\rotatebox[origin=c]{90}{\centering \textbf{\normalsize ST--GRB}}}
		& \normalsize\textbf{5} & (59,14) & (9,14) & 8.46\% & 37.0\% & 0.18 \si{\second} & 15.0\% & 8.83\% & 24.4\% & 5.29\% & 0.17 \si{\second} \\ 
		& \normalsize\textbf{10} & (73,19) & (9,17) & 3.51\% & 14.9\% & 0.22 \si{\second} & 15.7\% & 3.24\% & 29.5\% & 5.03\% & 0.37 \si{\second} \\
		& \normalsize\textbf{20} & (73,21) & (9,19) & 1.17\% & 1.68\% & 0.31 \si{\second} & 15.0\% & 1.35\% & 28.7\% & 1.86\% & 0.41 \si{\second} \\
		& \normalsize\textbf{50} & (74,21) & (9,19) & 0.84\% & 0.88\% & 0.38 \si{\second} & 15.3\% & 1.21\% & 29.2\% & 1.24\% & 0.51 \si{\second} \\
		
		\bottomrule 
	\end{tabular}
}
\end{table}

Lastly, we investigate the effect of the training sample size $M$ and the extrapolative generalization capabilities of the two methods. The results of this analysis are summarized in Table~\ref{tab:bifurcation ST--RB snapshots}. For $M=5, 10, 20, 50$, we report the velocity and pressure reduced basis sizes, the average relative errors and simulation durations for the usual ten testing parameters $\bm{\mu}^*_i, \ i=1, \dots, 10$, sampled within the training domain $\mathcal{P}_f$, and the velocity relative test errors and average simulation wall times for four parameters $\bm{\bar{\mu}}_i, \ i=1,\dots,4$, sampled outside the training domain. Specifically, we consider $\bm{\bar{\mu}}_1 = [\underline{2.0}, 0.2, 0.6]$, $\bm{\bar{\mu}}_2 = [6.0, \underline{0.4}, 0.6]$, $\bm{\bar{\mu}}_3 = [\underline{2.0}, \underline{0.4}, 0.5]$, and $\bm{\bar{\mu}}_4 = [6.0, 0.2, \underline{0.9}]$, where underlined values indicate entries lying outside the training envelope.
Regarding interpolation, we observe that the canonical RB method scales more favourably with the number of training snapshots compared to ST--GRB. In particular, SRB--TFO reaches a plateau in accuracy for both velocity and pressure already at $M=10$, with its computational cost remaining largely unaffected by the training sample size. Conversely, the accuracy of ST--GRB keeps improving as additional snapshots are added to the training dataset, although the gain becomes marginal beyond $M=20$. Moreover, the average wall times of ST--GRB increase with $M$, mainly due to the computation of the initial guess for the Newton iterations. Overall, we can claim that, for moderate accuracy targets, ST--GRB trades offline data efficiency for enhanced online performance.
The traditional RB method also appears to provide better generalization capabilities compared to its space--time--reduced counterpart. To better interpret the extrapolation results, we note that $\bm{\bar{\mu}}_1$, $\bm{\bar{\mu}}_2$, and $\bm{\bar{\mu}}_3$ lie outside $\mathcal{P}_f$ due to their first two entries, which influence the inflow waveform and thus the temporal diversity of the solution manifold. By contrast, only the last entry of $\bm{\bar{\mu}}_4$ --- related to spatial variability --- lies outside the training envelope, inducing a large flow imbalance across the two outlet branches. As it might be expected, the absence of dimensionality reduction in time leads to a (minor) deterioration of SRB--TFO accuracy, compared to the interpolation results, only for $\bm{\bar{\mu}}_4$, where mild overfitting effects can also be observed. On the other hand, the performances of ST--GRB are significantly compromised when extrapolating with respect to temporal dynamics. In particular, severe losses in accuracy are observed when the first parameter entry, controlling the flow perturbation frequency, lies outside the training range. High relative errors in velocity are consistently obtained for both $\bm{\bar{\mu}}_1$ and $\bm{\bar{\mu}}_3$, regardless of the value of $M$. Possible strategies to mitigate this limitation are discussed in Section~\ref{sec:conclusions}. Nonetheless, it is worth noting that ST--GRB demonstrates robust generalization capabilities (and no overfitting) with respect to both the flow perturbation amplitude ($\bm{\bar{\mu}}_2$) and the outflow imbalance factor ($\bm{\bar{\mu}}_4$). Finally, we observe that parametric extrapolation generally results in slightly increased computational times for ST--GRB --- for instance, $+34\%$ for $M=50$, from $0.38 \ \si{\second}$ to $0.51 \ \si{\second}$ --- since more Newton iterations are needed to achieve convergence.

\subsection{Test case 2: Coupled Momentum model in a bent tube} 
\label{subsec: membrane test 1}
In this test case, we numerically investigate the application of ST--GRB to the reduced FSI problem which stems from the use of the coupled momentum model to account for the vessel wall compliance. We consider a bent tube geometry, featuring a diameter--length aspect ratio of 1:4 (see Figure~\ref{fig: geometries} -- middle). We avoid the flow rate variability (hence $\mathcal{P}_f = \emptyset$), so that the problem parametrization entirely depends on the vessel wall properties and it is affine. The space of membrane parameters, whose generic element is $\bm{\mu_m} = [h_s, \rho_s, E, \nu]$, is defined as $\mathcal{P}_m := [0.05, 0.15] \times [1.08, 1.80] \times [2\cdot10^6, 6\cdot10^6] \times [0.35, 0.5] \subset \R^4$. Also, we do not take into account the elastic response of the surrounding tissues, thus setting $c_s = 0$ (see Eq.\eqref{eq: space-time FOM membrane matrices}). At the inlet $\Gamma_{\operatorname{IN}}^1$, we (weakly) prescribe a parabolic velocity profile that matches the flow rate $g^t = g_1^t(\cdot; \bm{0})$, where $g_1^t$ is defined as in Eq.\eqref{eq: flow rate 1}. At the outlet $\Gamma_{\operatorname{OUT}}^1$ we impose absorbing boundary conditions, with a resistance $R = 100 \ \si{\gram \ \per{(\centi\metre^4 \cdot \second)}}$. Furthermore, we enforce the following boundary conditions at the inlet/outlet rings:
\begin{equation}
\label{eq: BC rings}
\begin{alignedat}{3}
&\vec{u} \cdot \bm{n}_{\operatorname{IN}} = 0~; \quad &&\bm{n}_{\operatorname{IN}}^T\sigma(\vec{u}, p)(I - \bm{n}_{\operatorname{IN}} \otimes \bm{n}_{\operatorname{IN}}) = \vec{0}~; \qquad &&\text{on} \ \ \Gamma^{1,r}_{\operatorname{IN}} := \Gamma \cap \Gamma_{\operatorname{IN}}^1~; \\
&\vec{u} \cdot \bm{n}_{\operatorname{OUT}} = 0~; \quad &&\bm{n}_{\operatorname{OUT}}^T\sigma(\vec{u}, p)(I - \bm{n}_{\operatorname{OUT}} \otimes \bm{n}_{\operatorname{OUT}}) = \vec{0}~; \qquad &&\text{on} \ \ \Gamma^{1,r}_{\operatorname{OUT}} := \Gamma \cap \Gamma_{\operatorname{OUT}}^1~;
\end{alignedat}
\end{equation}
where $\bm{n}_{\operatorname{IN}}$, $\bm{n}_{\operatorname{OUT}} \in \R^3$ are the outward unit normal vectors to $\Gamma_{\operatorname{IN}}^1$, $\Gamma_{\operatorname{OUT}}^1$, respectively. Eq.\eqref{eq: BC rings} enables tangential deformations of the inlet/outlet faces, while it prevents any displacement in the normal direction. Unless otherwise specified, in the following tests we consider $M=50$ snapshots, $\varepsilon_u=\varepsilon_p=\varepsilon_\lambda^t=10^{-3}$, $\varepsilon_\lambda^s = 10^{-5}$ as POD tolerances, $n_{c,J}=0$ affine components for the Jacobian of the convective term, and an initial guess $\widehat{\bm{w}}^{(0)}$ that is computed with the PODI approach. Furthermore, we select $n_c$ as the number of velocity modes obtained after the truncated POD in space, thus neglecting the supremizers when computing the affine approximation of the convective term. \\

\begin{figure}[t!]
\centering
\includegraphics[width=0.495\textwidth]{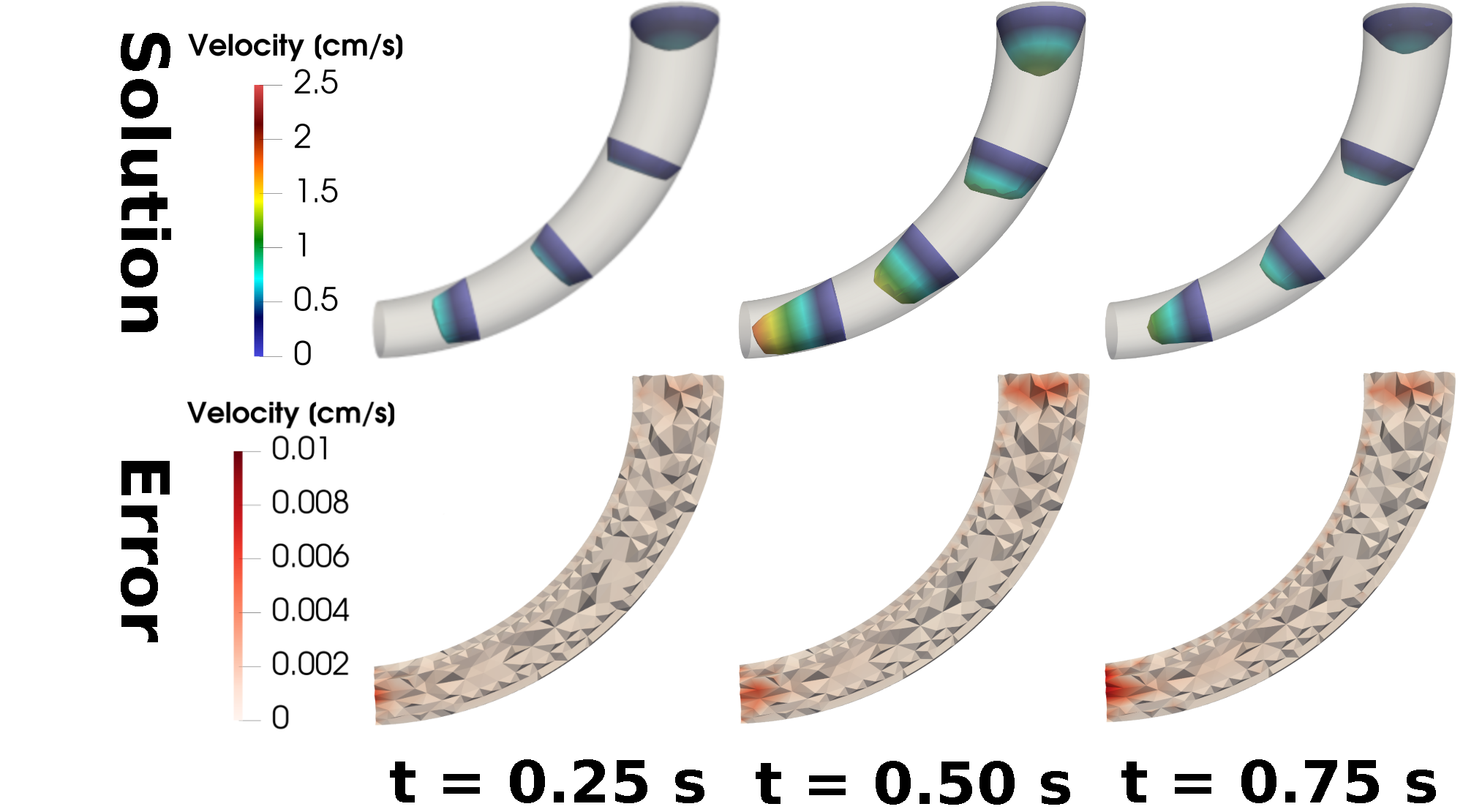}
\includegraphics[width=0.495\textwidth]{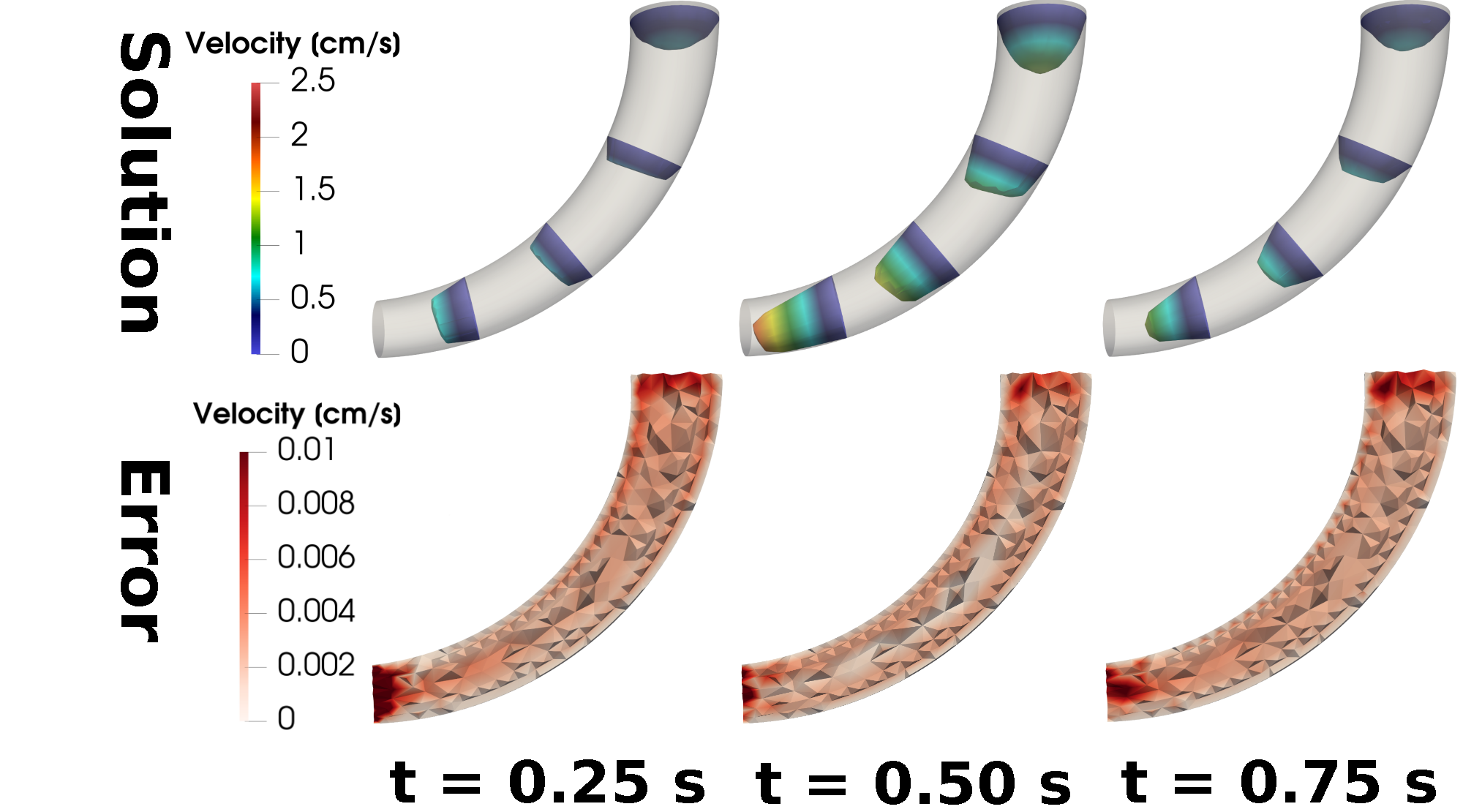}
\caption{Velocity field at $4$ different slices along the vessel centerline (top) and pointwise absolute error on the median slice (bottom) at three time instants, obtained with SRB--TFO (left) and ST--GRB (right) for $\bm{\mu}^* = [0.12, 1.44, 5.20 \cdot 10^6, 0.36]$.}
\label{fig: coupled momentum solution velocity}

\medskip

\includegraphics[width=0.495\textwidth]{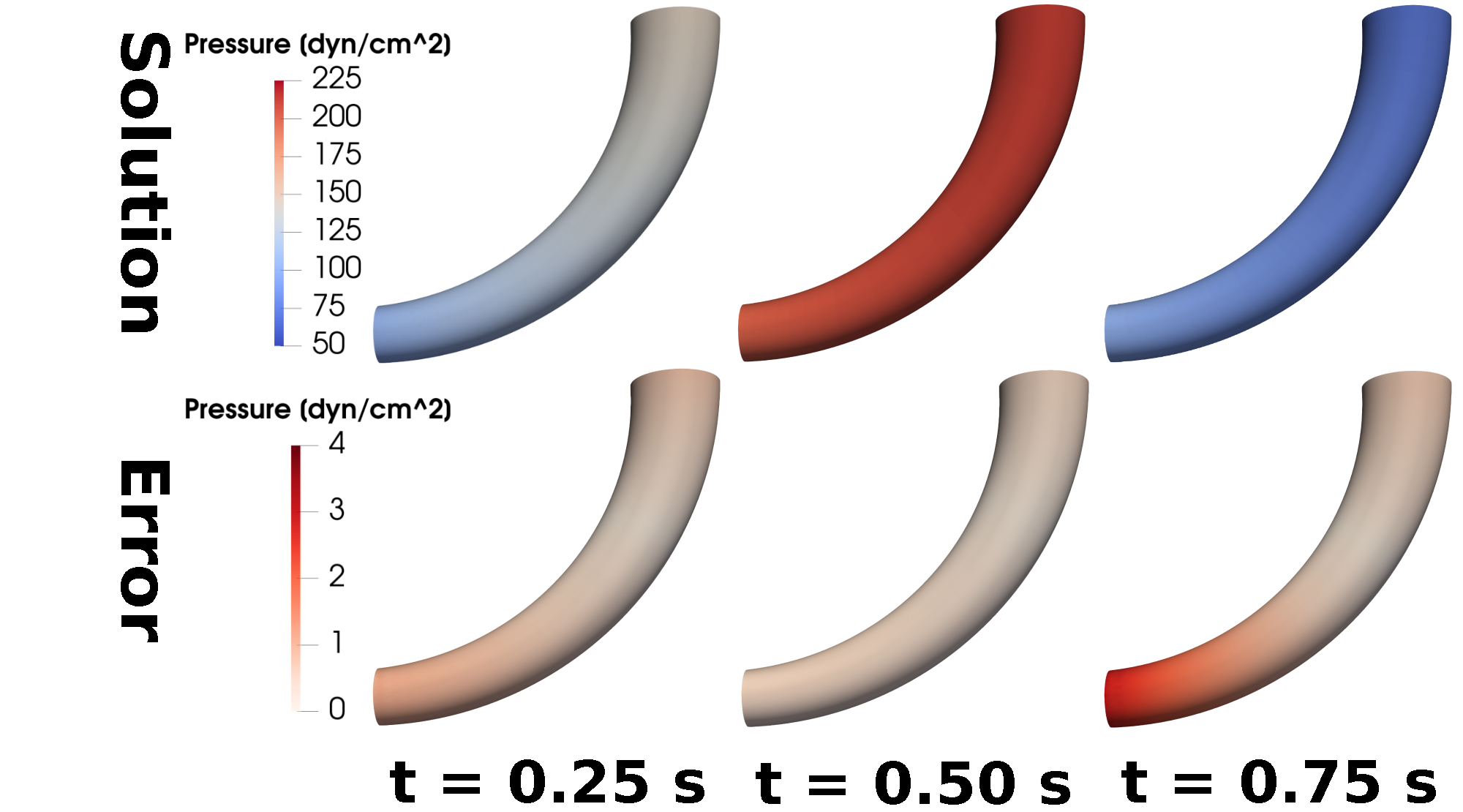}
\includegraphics[width=0.495\textwidth]{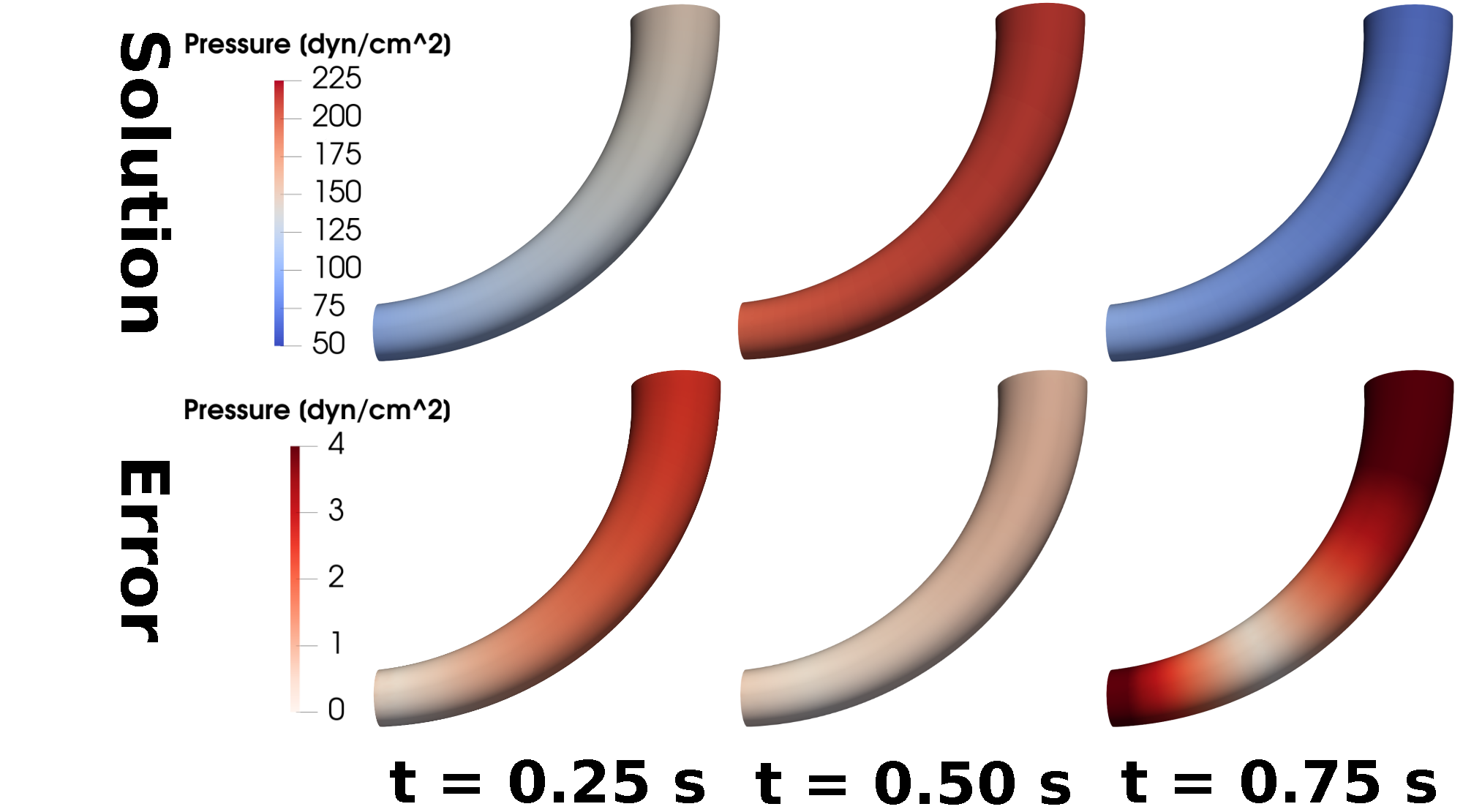}
\caption{Pressure field (top) and corresponding pointwise absolute error (bottom) at three time instants, obtained with SRB--TFO (left) and ST--GRB (right) for $\bm{\mu}^* = [0.12, 1.44, 5.20 \cdot 10^6, 0.36]$.}
\label{fig: coupled momentum solution pressure}
\end{figure}

Figures~\ref{fig: coupled momentum solution velocity},~\ref{fig: coupled momentum solution pressure} show, respectively, the velocity and the pressure field at three equispaced time instants obtained with SRB--TFO and ST--GRB for $\bm{\mu}^* = [0.12, 1.44, 5.20 \cdot 10^6, 0.36]$, and the corresponding pointwise absolute errors with respect to the FOM solution. This parameter was not considered for the computation of the reduced bases during the offline phase. The relative errors, measured in the natural spatio--temporal norms, for this particular test case are: $0.50 \%$ for the velocity, $0.58 \%$ for the pressure, $1.33 \%$ for the displacement with SRB--TFO; $1.26 \%$ for the velocity, $2.51 \%$ for the pressure, $1.73 \%$ for the displacement with ST--GRB. From a qualitative standpoint, we remark that the solution approximations yielded by the two methods are substantially indistinguishable by the eye. Nonetheless, as pointed out in Section~\ref{subsec: navier-stokes test}, we recognize that SRB--TFO attains better accuracy levels, since the temporal dynamics are processed in a high--fidelity fashion. For what concerns the velocity field, in particular, we notice larger errors close to the extremal boundaries. In the outlet region, this behaviour could be expected, since the flow is accelerated by the vessel lumen narrowing. At the inlet, instead, large errors can be linked to the weak imposition of the Dirichlet boundary conditions; in fact, the normalized relative errors on the Lagrange multipliers are relatively high both with SRB--TFO ($4.04 \%$) and with ST--GRB ($4.33 \%$). \\

Table~\ref{tab:tube ST--RB pod tolerances} reports the reduced bases sizes, the number of convective affine components and the performance metrics of SRB--TFO and ST--GRB for five different POD tolerances. 
From a general perspective, we observe that ST--GRB yields slightly less accurate results than SRB--TFO, but at a significantly lower computational cost. This suggests that the additional reduction layer introduced by temporal dynamics compression does not compromise approximation quality, while substantially enhancing efficiency. In fact, the advantages of ST--GRB over SRB--TFO are more pronounced compared to the first test case, because the temporal reduced basis sizes are smaller. This is primarily due to the imposition of a non--perturbed inflow rate and to the simpler domain topology. Furthermore, absorbing boundary conditions at the outlet prevent the onset of non--physiological reflecting waves.
From a quantitative standpoint, we note that both velocity and pressure errors are roughly one order of magnitude larger than the prescribed POD tolerance. Remarkably, the same consideration holds also for the displacement field, even if its values are roughly three orders of magnitude smaller than the velocity ones. Incidentally, we underline that the larger errors made by ST--GRB can be justified taking into account both the temporal dynamics compression and the different treatment of the structural displacement. Indeed, high--fidelity--in--time approaches, such as the FOM and SRB--TFO, feature an explicit embedding of the displacement dynamics within the Navier--Stokes equations, based on the kinematic coupling condition. Conversely, our implementation of ST--GRB implicitly encapsulates the wall mechanics into the fluid's equations (see Eq.\eqref{eq: space-time displacement}).
For what concerns computational efficiency, we highlight that ST--GRB is much faster than SRB--TFO for all the considered POD tolerances, even though its gains shrink as higher precision results are sought. For instance, ST--GRB is $600$ times faster than SRB--TFO for $\varepsilon=10^{-2}$ ($2.64 \ \si{\milli\second}$ vs. $1.57 \ \si{\second}$), $150$ times faster for $\varepsilon=10^{-3}$ ($16.9 \ \si{\milli\second}$ vs. $2.50 \ \si{\second}$), and ``only'' $60$ times faster for $\varepsilon=10^{-4}$ ($58.5 \ \si{\milli\second}$ vs. $3.68 \ \si{\second}$). 
Lastly, we highlight that considering the (approximate) convective term Jacobian in the Newton iterations does not substantially improve accuracy, while it moderately dilutes the computational gains. For example, if we employ ST--GRB and we set $\varepsilon = 10^{-3}$ and $n_{c,J} = 10$, the normalized average relative test errors on velocity and pressure decrease, respectively, of $2.40 \%$ (from $13.3$ to $13.0$) and $0.08 \%$ (from $26.4$ to $26.3$), while the average online time--to--solution increases of $32.2 \%$ (from $1.14 \ \si{\milli\second}$ to $1.69 \ \si{\milli\second}$). \\

\begin{table}[t!]
\centering
\caption{Results obtained with SRB--TFO and ST--GRB on the coupled momentum bent tube test case, for different POD tolerances $\varepsilon$. In particular: (left) size of the spatial and temporal reduced bases for velocity, pressure, and Lagrange multipliers, and number of affine components for the reduced convective term; (center) reduction factor and average speedup; (right) normalized average relative test errors on velocity, pressure, and displacement.} 
\label{tab:tube ST--RB pod tolerances}
\resizebox{.95\textwidth}{!}{
\begin{tabular}{ccccccccccc} 
	\toprule
	& &
	\multicolumn{4}{c}{\large\textbf{ROM size}}&
	\multicolumn{2}{c}{\large\textbf{Efficiency}}&
	\multicolumn{3}{c}{\large\textbf{Accuracy}}\\
	
	\cmidrule(lr){3-6} \cmidrule(lr){7-8} \cmidrule(lr){9-11}
	
	&
	\normalsize$\bm{\varepsilon}$ & 
	\normalsize$\bm{(n_u^s,n_u^t)}$ & 
	\normalsize$\bm{(n_p^s,n_p^t)}$ &
	\normalsize$\bm{(n_\lambda^s,n_{\lambda}^t)}$ &
	\normalsize$\bm{n_c}$ &
	\normalsize\textbf{RF} & 
	\normalsize\textbf{SU} & 
	\normalsize$\bm{E_u / \varepsilon_u}$ & 
	\normalsize$\bm{E_p / \varepsilon_p}$ &
	\normalsize$\bm{E_d / \varepsilon_u}$ \\
	
	\midrule
	
	\multirow{5}{*}{\rotatebox[origin=c]{90}{\centering \textbf{\normalsize SRB--TFO}}}
	& \normalsize\textbf{1e-2} & (10,$\ \cdot \ $) & (2,$\ \cdot \ $) & (3,$\ \cdot \ $) & 5 & 6.08e3  & 4.06e3 & 4.97 & 26.6 & 6.02 \\ 
	& \normalsize\textbf{5e-3} & (12,$\ \cdot \ $) & (2,$\ \cdot \ $) & (4,$\ \cdot \ $) & 6 & 5.07e3 & 3.86e3 & 9.72 & 39.5 & 9.04 \\
	& \normalsize\textbf{1e-3} & (17,$\ \cdot \ $) & (2,$\ \cdot \ $) & (5,$\ \cdot \ $) & 10 & 3.80e3 & 2.58e3 & 4.65 & 5.80 & 9.96 \\
	& \normalsize\textbf{5e-4} & (20,$\ \cdot \ $) & (2,$\ \cdot \ $) & (6,$\ \cdot \ $) & 12 & 3.26e3 & 2.92e3 & 10.9 & 11.6 & 17.3 \\
	& \normalsize\textbf{1e-4} & (32,$\ \cdot \ $) & (4,$\ \cdot \ $) & (8,$\ \cdot \ $) & 20 & 2.07e3 & 1.75e3 & 21.6 & 24.0 & 9.65  \\
	
	\midrule
	
	\multirow{5}{*}{\rotatebox[origin=c]{90}{\centering \textbf{\normalsize ST--GRB}}} 
	& \normalsize\textbf{1e-2} & (10,3)  & (2,2)  & (3,2)   & 5  & 2.28e6 & 2.44e6 & 6.00 & 12.9 & 6.23 \\
	& \normalsize\textbf{5e-3} & (12,3)  & (2,2)  & (4,2)   & 6  & 1.90e6 & 1.96e6 & 9.63 & 24.1 & 9.94 \\
	& \normalsize\textbf{1e-3} & (17,7)  & (2,3)  & (5,4)   & 10 & 6.29e5 & 3.82e5 & 13.3 & 26.4 & 14.7 \\
	& \normalsize\textbf{5e-4} & (20,9)  & (2,5) &  (6,5)   & 12 & 4.15e5 & 5.41e5 & 20.4 & 42.8 & 21.4 \\
	& \normalsize\textbf{1e-4} & (32,22) & (4,9)  & (8,9)   & 20 & 1.12e5 & 1.10e5 & 23.8 & 31.2 & 9.88 \\
	
	\bottomrule 
\end{tabular}
}
\end{table}

Since many h{\ae}modynamic problems are characterized by quasi--periodic behaviours due to the cardiac cycle, we further investigate the application of ST--GRB under conditions of repeated inflow over multiple cycles. Specifically, we adopt the same numerical setup described above, but replicate the inflow rate $g^t$ over $5$ periods. To reduce computational costs, in the FOM simulations we use $\mathbb{P}_1^b - \mathbb{P}_1$ elements, set $\Delta t = 4 \ \si{\milli\second}$, and consider $M = 25$ training snapshots and $M^* = 5$ testing snapshots.
Additionally, the first cycle is discarded during both the offline and online phases, as it represents a ramp--up transient whose associated solution fields differ significantly from those in the subsequent cycles, which exhibit a quasi--periodic behaviour. As a result, the problem features inhomogeneous and parameter--dependent initial conditions, whose efficient incorporation within the space--time--reduced setting is addressed via a temporal lifting procedure, described in detail in Appendix~\ref{app: incorporating inomogeneous initial conditions}.
Table~\ref{tab: tube ST-RB multibeat} summarizes the results obtained with the traditional RB method and the ST--GRB approach, reporting the temporal reduced basis sizes, average simulation times, and average relative test errors on velocity, pressure, and displacement. Figure~\ref{fig: tube ST-RB multibeat} displays the average test absolute errors in velocity and pressure over time, as obtained with the considered ROMs.
We investigate three variants of ST--GRB, distinguished by the number of periods $N_T$ used to construct the temporal reduced bases. In particular, for the cases $N_T = 1$ and $N_T = 2$, each ST--GRB cycle provides the initial condition for the subsequent one, thereby enforcing a global regularity of the solution across cycles. This leads us to propose a windowed ST--GRB method, similar to the approach introduced in~\cite{shimizu2021windowed} for a least--squares Petrov--Galerkin formulation.
We observe that, for all considered values of $N_T$, ST--GRB outperforms SRB--TFO, yielding comparably accurate approximations at reduced computational costs. As expected, using shorter time windows in the construction of the temporal reduced bases generally leads to a lower number of temporal modes, at least for the velocity field. However, this does not result in additional computational savings with respect to the baseline approach, due to the overheads introduced by the extraction and management of the initial conditions.
Incidentally, one could envision the use of periodic temporal modes to further shrink the reduced dimensions and, in turn, the computational effort. However, tentative approaches based on either orthonormal sinusoidal modes or data--adaptive periodic modes yielded unsatisfactory results in our setting, significantly worse than those reported in Table~\ref{tab: tube ST-RB multibeat}. The main reason is that, albeit the enforced inflow rate being periodic, the numerical solution is only quasi--periodic, gradually converging to a periodic regime over the cycles. Therefore, on the one hand, the adoption of recurring modes in time may provide computational gains if several initialization cycles are discarded. On the other hand, the use of non--periodic modes offers better generalization capabilities, since it seamlessly allows, for instance, to perform accurate simulations under variations of the physical parameters across cycles.

\begin{table}[t!]
\centering
\caption{Results obtained with SRB--TFO and ST--GRB on the coupled momentum bent tube test case, considering a periodic inflow rate over $5$ cycles; different variants of ST--GRB differentiate depending on the number of periods $N_T$ taken into account to compute the temporal reduced bases. In particular: (left) size of the temporal reduced bases for velocity, pressure, and Lagrange multipliers; (center) average wall time (WT) in \si{\milli\second}; (right) average test relative errors on velocity, pressure, and displacement.}
\label{tab: tube ST-RB multibeat}
\resizebox{.8\textwidth}{!}{
\begin{tabular}{ccccccccc} 
	\toprule
	\multirow{2}{*}{\large\textbf{Method}} & \multirow{2}{*}{$\bm{N_T}$} & \multicolumn{3}{c}{\large\textbf{ROM size}} & \large\textbf{Efficiency} & \multicolumn{3}{c}{\large\textbf{Accuracy}} \\
	
	\cmidrule(lr){3-5} \cmidrule(lr){6-6} \cmidrule(lr){7-9}
	
	& & $\bm{n_u^t}$ & $\bm{n_p^t}$ & $\bm{n_\lambda^t}$ & \textbf{WT} & $\bm{E_u}$ & $\bm{E_p}$ & $\bm{E_d}$ \\
	
	\midrule
	
	\textbf{SRB--TFO} & // & // & // & // & 366 \si{\milli\second}  & 0.72\% & 1.13\% & 2.62\% \\
	
	\midrule
	
	\multirow{3}{*}{\textbf{ST--GRB}} & 1 & 7 & 5 & 4 & 44 \si{\milli\second} & 0.59\% & 1.18\% & 1.66\% \\
	
	& 2 & 11 & 5 & 4 & 43 \si{\milli\second}  & 0.71\% & 1.04\% & 1.66\% \\
	
	& 4 & 14 & 5 & 4 & 16 \si{\milli\second}  & 0.66\% & 2.10\% & 1.77\% \\
	
	\bottomrule
	
\end{tabular}
}

\bigskip

\centering
\includegraphics[width=.95\textwidth]{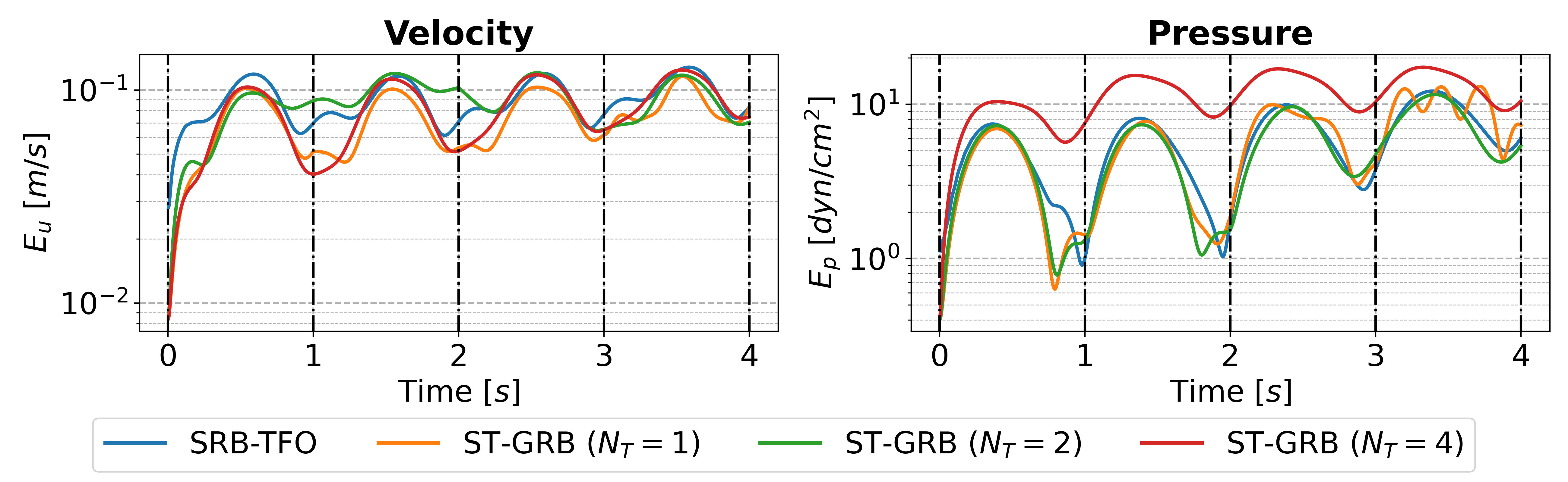}
\captionof{figure}{Velocity and pressure absolute errors trend over time, obtained with SRB--TFO and different variants of ST--GRB on the coupled momentum bent tube test case.}
\label{fig: tube ST-RB multibeat}

\end{table}

\subsection{Test case 3: Coupled Momentum model in aorta and iliac arteries} \label{subsec: membrane test 2}
This last numerical experiment analyses the application of ST--GRB in a realistic scenario. We consider the physiological geometry of an aorta with left and right external iliac arteries~\cite{wilson2013vascular} (see Figure~\ref{fig: geometries} -- right). At the inlet, we weakly enforce a parabolic velocity profile that matches a patient--specific supraceliac blood flow waveform, measured from PC--MRI and spanning over two heartbeats (see Figure~\ref{fig: flow rates} -- right). In fact, we enhance the solution manifold variability by parametrizing the available measurements. Specifically, we introduce three parameters that characterize the systolic flow peak, the time--to--peak--systole and the diastolic flow peak, and the space of flow parameters is defined as $\mathcal{P}_f = [12.0, 16.3] \times [0.11, 0.16] \times [1.57, 2.62] \subset \R^3$. The space of membrane parameters, whose generic element is $\bm{\mu_m} = [h_s, \rho_s, E, \nu]$, is defined as $\mathcal{P}_m := [0.075, 0.125] \times [1.08, 1.50] \times [3\cdot10^6, 5\cdot10^6] \times [0.4, 0.5] \subset \R^4$. Following \cite{malossi2013numerical, colciago2014comparisons}, we account for the constraint of the surrounding tissues and set $c_s = 6 \cdot 10^4 \ \si{\dyne~\per\centi\metre\cubed}$. At the two outlets, we impose equivalent resistance boundary conditions, setting $R = 1.6 \cdot 10^4 \  \si{\gram~\per(\centi\metre^4\cdot\second)}$. This value was taken from a \emph{SimVascular}~\cite{updegrove2017simvascular} tutorial~\footnote{\url{https://simvascular.github.io/clinical/aortofemoral1.html}} conducted on the same geometry (albeit featuring more downstream branches), and it was computed from subject--specific flow rate measurements. Additionally, we simply enforce zero velocity at the inlet/outlet rings. For what concerns the ROM simulations, we always consider $M=25$ snapshots and a space POD tolerance on Lagrange multipliers that is two orders of magnitude lower than the other tolerances. Preliminary tests --- analogous to the ones performed in Section~\ref{subsec: navier-stokes test} --- suggested that considering $n_c=15$ convective affine components provides the best trade--off between accuracy and efficiency; also, we neglect the convective Jacobian affine components in this test case. Finally, with ST--GRB, we compute the initial guess using the PODI approach.\\

Table~\ref{tab:aorta ST-RB results} reports the reduced basis sizes and the performance metrics of SRB--TFO and ST--GRB, considering three different POD tolerances. The outcomes are similar to the ones discussed in Section~\ref{subsec: membrane test 1} for the tube test case. In particular, for any fixed POD tolerance $\varepsilon$, ST--GRB is less accurate than SRB--TFO; for instance, the velocity error is roughly doubled in all cases. However, ST--GRB is always much faster than SRB--TFO and its computational gains get amplified at large values of $\varepsilon$. We highlight that the average online time--to--solution of ST--GRB simulations features a much stronger dependency on the POD tolerance than the one of SRB--TFO. Ultimately, it is possible to claim that ST--GRB outperforms SRB--TFO on the test case at hand. Indeed, by comparing the results of SRB--TFO for $\varepsilon=5 \cdot 10^{-3}$ and of ST--GRB for $\varepsilon=10^{-3}$, we observe similar relative errors, but ST--GRB is almost two times faster (6.76 \si{\second} vs. 12.2 \si{\second}). \\

\begin{table}[t!]
\centering
\caption{Results obtained with SRB--TFO and ST--GRB on the coupled momentum aorta test case. In particular: (left) size of the spatial and temporal reduced bases for velocity, pressure, and Lagrange multipliers; (center) average wall time (WT) in \si{\second}; (right) average relative test errors on velocity, pressure, and displacement.}
\label{tab:aorta ST-RB results}
\resizebox{.95\textwidth}{!}{
\begin{tabular}{ccccccccc} 
	\toprule
	& &
	\multicolumn{3}{c}{\large\textbf{ROM size}}&
	\multicolumn{1}{c}{\large\textbf{Efficiency}}&
	\multicolumn{3}{c}{\large\textbf{Accuracy}}\\
	
	\cmidrule(lr){3-5} \cmidrule(lr){6-6} \cmidrule(lr){7-9}
	
	& 
	\normalsize{$\bm{\varepsilon}$} &
	\normalsize$\bm{(n_u^s,n_u^t)}$ & 
	\normalsize$\bm{(n_p^s,n_p^t)}$ &
	\normalsize$\bm{(n_\lambda^s,n_{\lambda}^t)}$ &
	\normalsize\textbf{WT} &  
	\normalsize$\bm{E_u}$ & 
	\normalsize$\bm{E_p}$ &
	\normalsize$\bm{E_d}$ \\
	
	\midrule
	
	\multirow{3}{*}{\normalsize\textbf{SRB--TFO}} & \textbf{1e-2} & (39, $\ \cdot \ $) & (2, $\ \cdot \ $) & (2, $\ \cdot \ $) & 13.1 \si{\second} & 6.13\% & 1.98\% & 1.41\% \\
	
	& \textbf{5e-3} & (61, $\ \cdot \ $) & (3, $\ \cdot \ $) & (3, $\ \cdot \ $) & 13.2 \si{\second} & 3.25\% & 0.47\% & 0.91\% \\
	
	& \textbf{1e-3} & (108, $\ \cdot \ $) & (6, $\ \cdot \ $) & (9, $\ \cdot \ $) & 14.1 \si{\second}  & 2.40\% & 0.10\% & 0.29\% \\
	
	\midrule
	
	\multirow{3}{*}{\normalsize\textbf{ST--GRB}} & \textbf{1e-2} & (39, 32) & (2, 6) & (2, 3) & 0.27 \si{\second} & 15.5\% & 3.14\% & 4.91\%  \\
	
	& \textbf{5e-3} & (61, 44) & (3, 8) & (3, 5) & 0.50 \si{\second} & 7.96\% & 1.64\% & 6.66\%  \\
	
	& \textbf{1e-3} & (108, 79) & (6, 19) & (9, 9) & 6.76 \si{\second} & 3.84\% & 0.48\% & 0.98\% \\
	
	\bottomrule 
\end{tabular}
}
\end{table}

\begin{figure}[t!]
\centering
\includegraphics[width=0.495\textwidth]{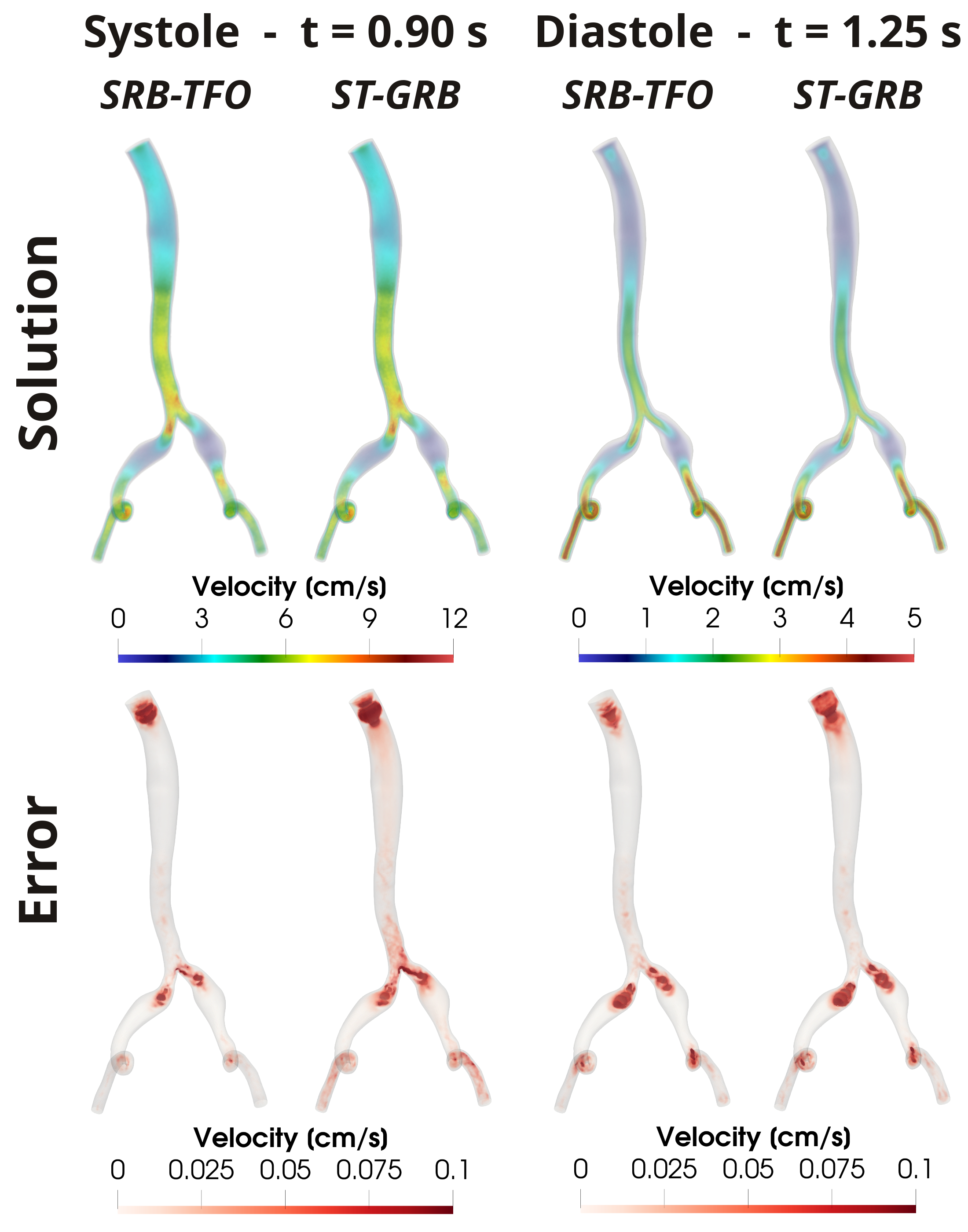}
\hfill
\includegraphics[width=0.495\textwidth]{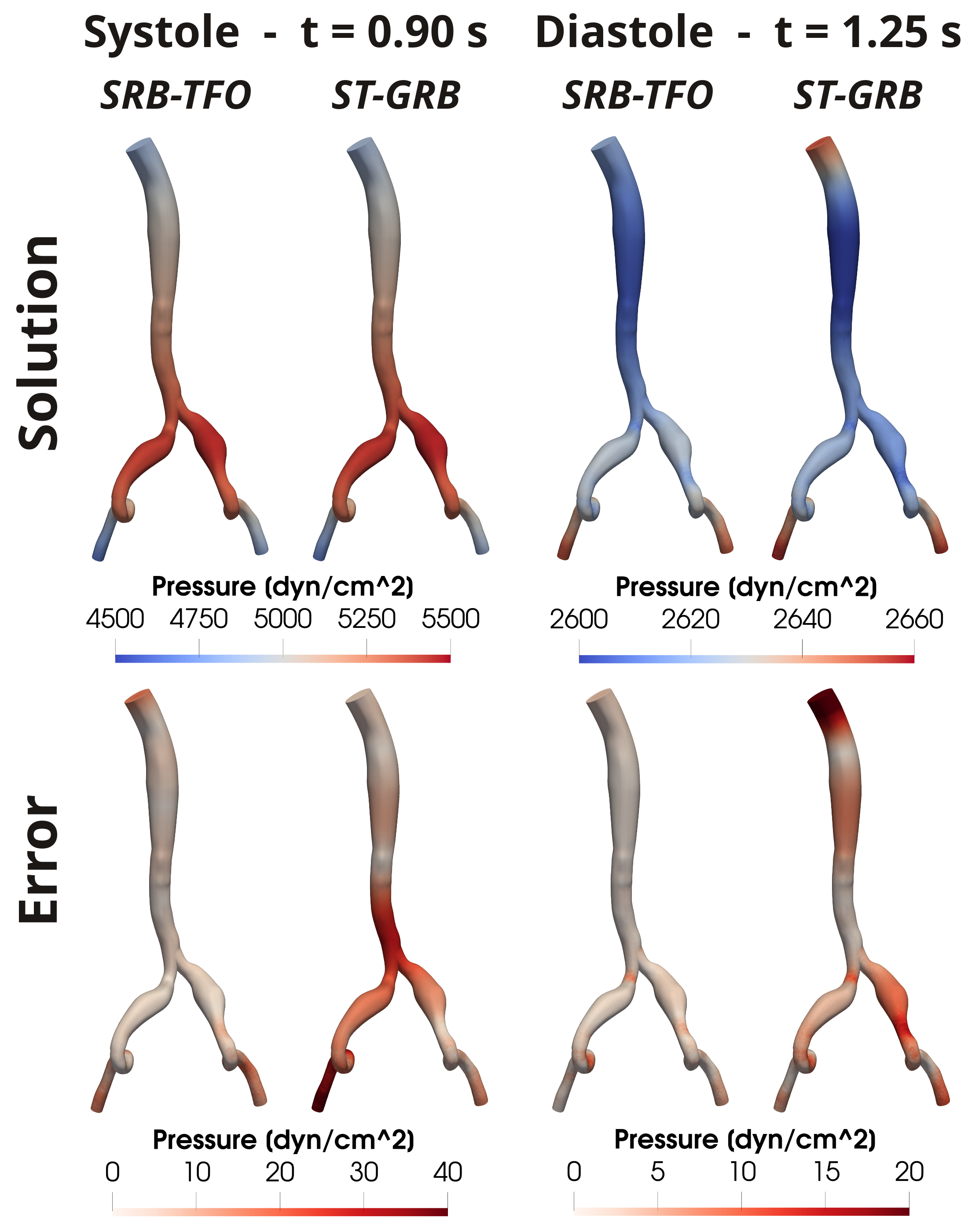}
\caption{Systolic and diastolic velocity (left) and pressure (right) solutions (top) and errors (bottom), obtained with SRB--TFO and ST--GRB for $\varepsilon= 10^{-3}$, $\bm{\mu_f^*} = [0.11, 13.2,  2.40]$, $\bm{\mu_m^*} = [0.097, 1.50, 3.07 \cdot 10^6, 0.43]$.}
\label{fig: results aorta}
\end{figure}

Figure~\ref{fig: results aorta} shows the velocity and pressure fields at systole and diastole of the second heartbeat, obtained with SRB--TFO and ST--GRB for $\varepsilon = 10^{-3}$, $\bm{\mu_f^*} = [0.11, 13.2,  2.40]$, $\bm{\mu_m^*} = [0.097, 1.50, 3.07 \cdot 10^6, 0.43]$; the corresponding pointwise errors with respect to the FOM solution are also represented. The relative errors in the natural spatio--temporal norms for this test case are: $2.53 \%$ for the velocity, $0.11 \%$ for the pressure and $0.28 \%$ for the displacement with SRB--TFO; $4.24 \%$ for the velocity, $0.43 \%$ for the pressure and $1.00 \%$ for the displacement with ST--GRB. For what concerns velocity, the largest errors are observed, at both timesteps, at the inlet and in the regions where complex flow patterns arise, namely the bifurcation and the final portion of the common iliac trunk. From a qualitative standpoint, the enhanced accuracy carried by SRB--TFO is nearly invisible by the eye, at the cost of more important computational costs. Regarding pressure, instead, both methods are extremely accurate, attaining accuracies aligned with the prescribed POD tolerance. Nonetheless, the higher precision level of SRB--TFO is more evident, as ST--RB produces relatively large errors in the bifurcation region during systole and at the inlet during diastole. \\

\begin{figure}[t!]
\begin{minipage}{.475\textwidth}
	\centering
	\includegraphics[width=0.99\textwidth]{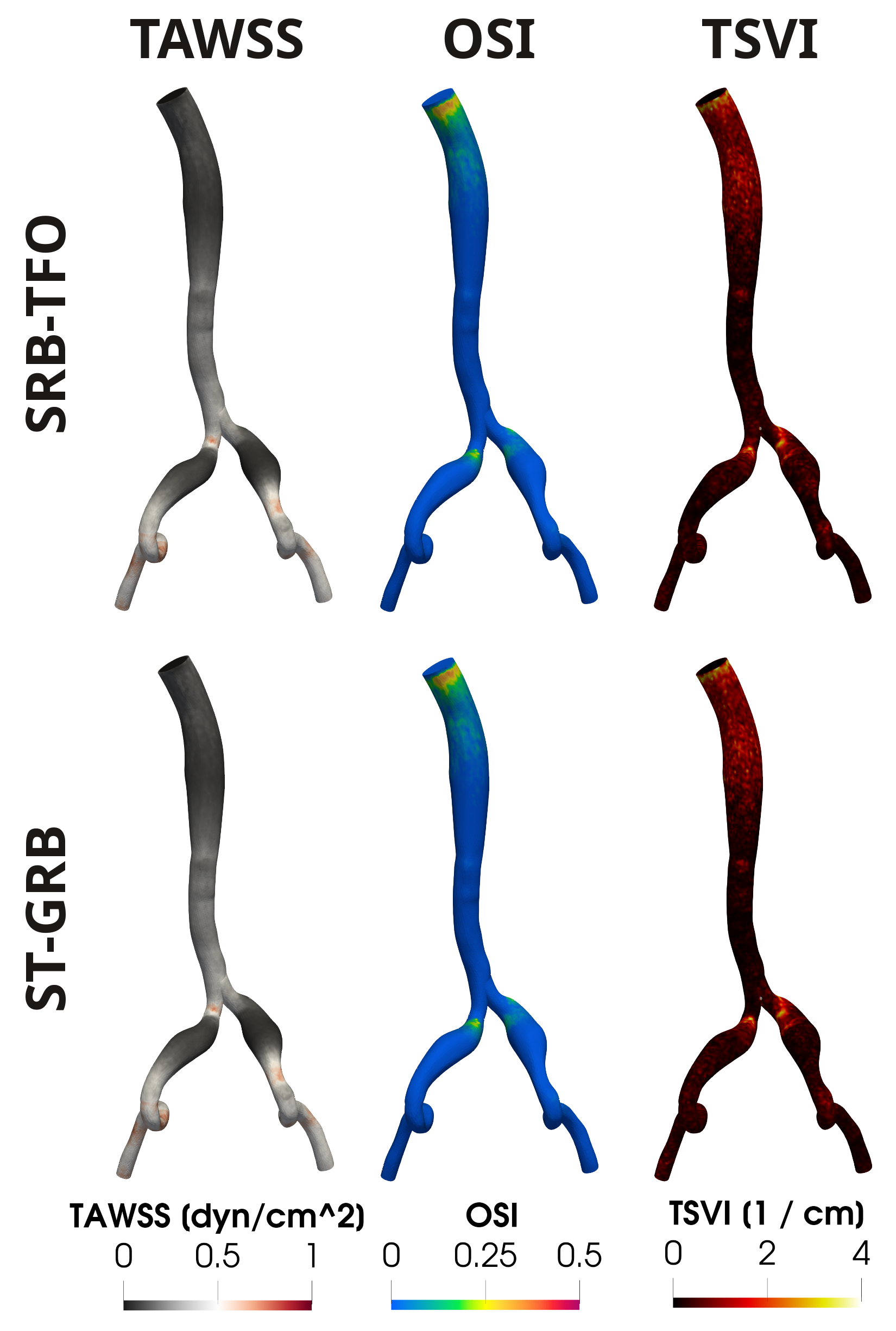}
\end{minipage}
\hfill
\begin{minipage}{.475\textwidth}
	\centering
	\includegraphics[width=0.99\textwidth]{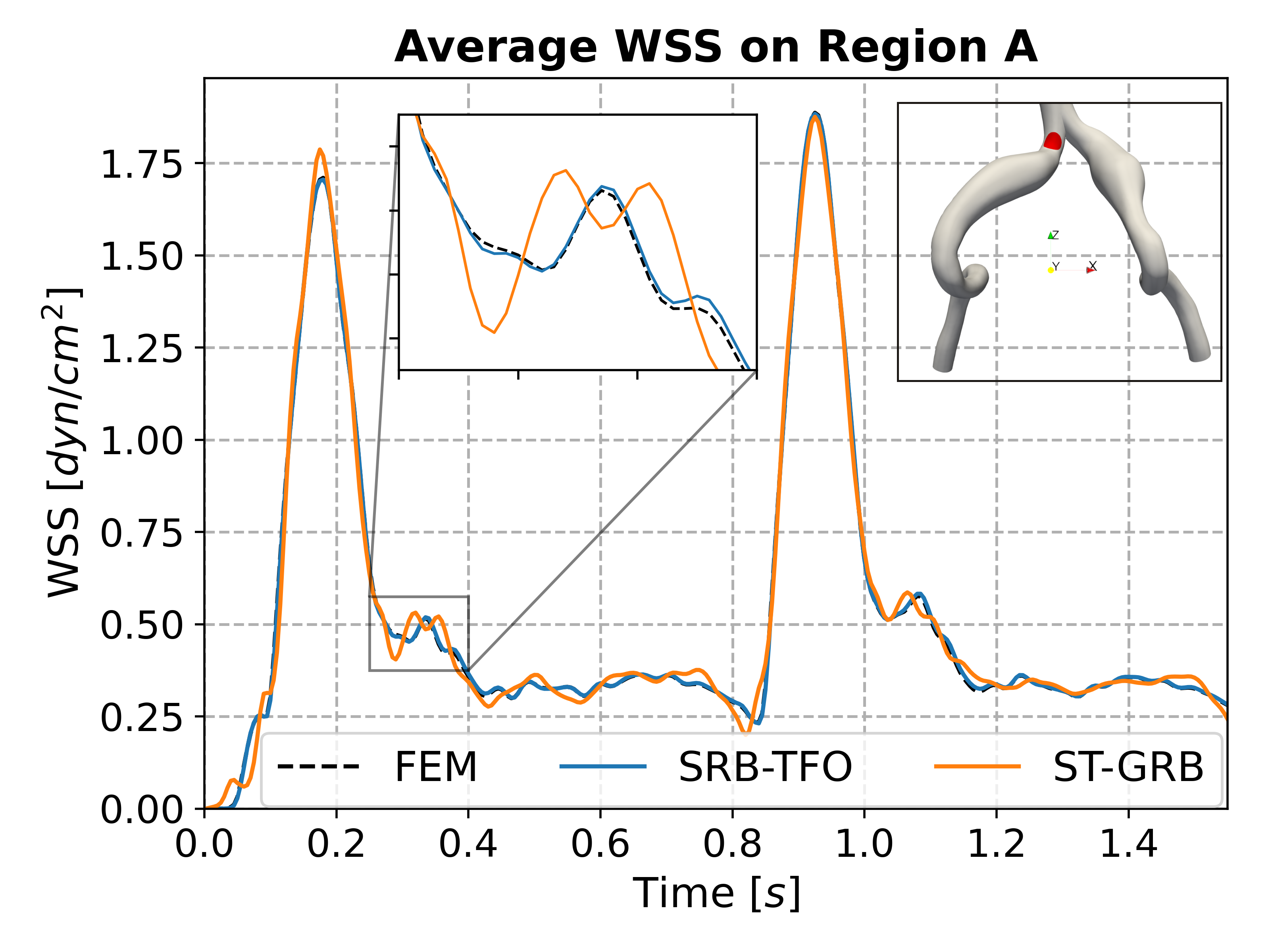}
	\includegraphics[width=0.99\textwidth]{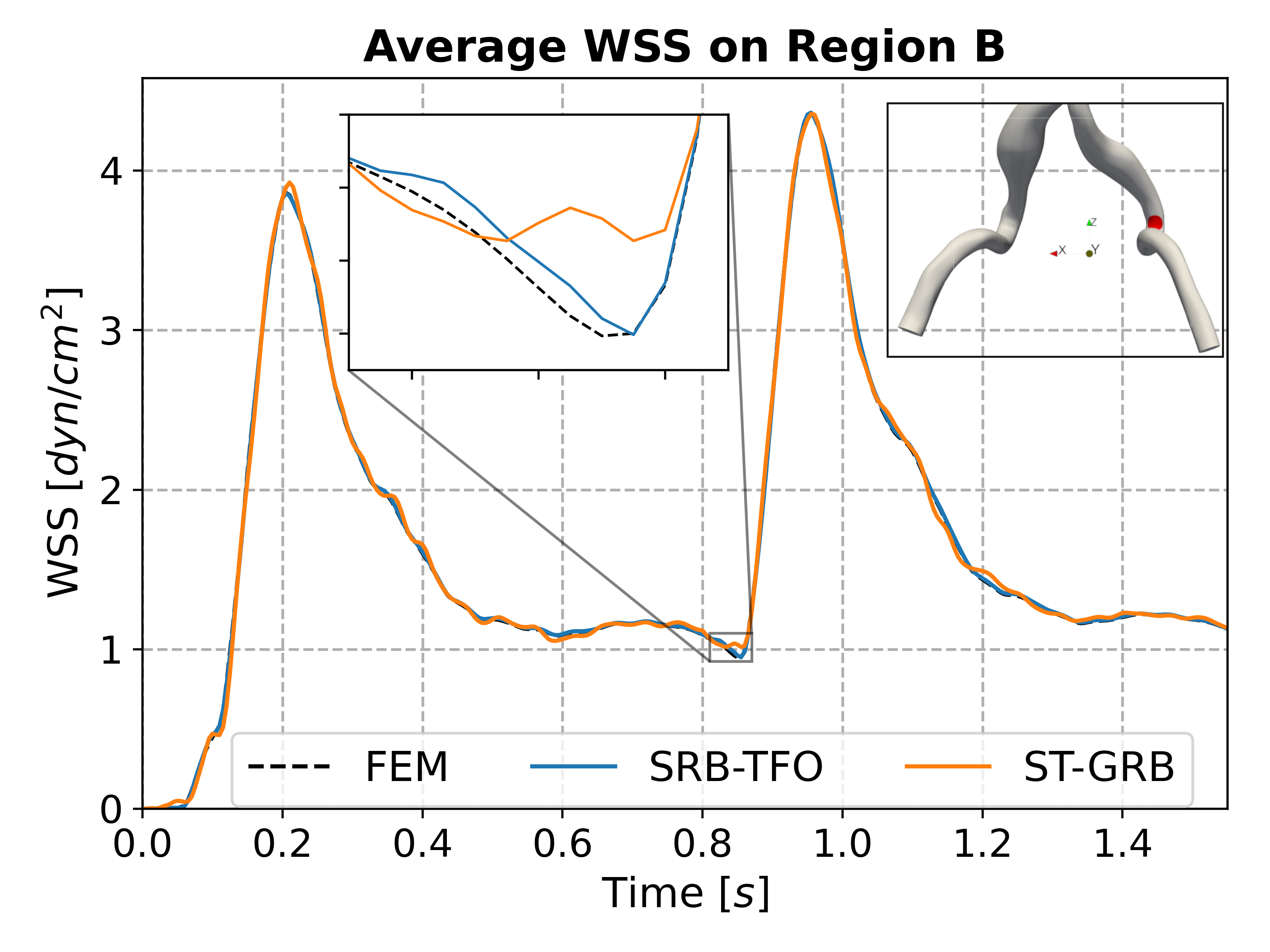}
\end{minipage}
\caption{Wall shear stress approximation results, obtained with SRB--TFO and ST--GRB for $\varepsilon= 10^{-3}$, $\bm{\mu_f^*} = [0.11, 13.2,  2.40]$, $\bm{\mu_m^*} = [0.097, 1.50, 3.07 \cdot 10^6, 0.43]$. In particular: (left) approximations of the time--averaged wall shear stress (TAWSS), oscillatory shear index (OSI) and topological shear variation index (TSVI); (right) trends over time of WSS, averaged in space over two different vessel wall portions, whose locations are depicted in the top--right corner.}
\label{fig: results aorta WSS}
\end{figure}

Finally, Figure~\ref{fig: results aorta WSS} reports the approximation results for the wall shear stress (WSS) --- the frictional force exerted by blood on the endothelium --- obtained with both SRB--TFO and ST--GRB. In particular, on the right we report the spatial distribution of three aggregated--in--time WSS--based quantities, namely the time averaged WSS (TAWSS), the oscillatory shear index (OSI)~\cite{ku1985pulsatile} and the topological shear variation index (TSVI)~\cite{denisco2020deciphering, mazzi2021early}. We remark that these quantities were found to be reliable predictors of atherosclerotic plaques growth and, in coronary arteries, of myocardial infarction~\cite{kumar2018high, candreva2022risk}. Notably, SRB--TFO and ST--GRB produce very similar results, which implies that the compression of the temporal dynamics does not affect accuracy on aggregated--in--time indicators. Incidentally, it is work underlining that taking into account vessel wall elasticity is crucial for accurately simulating the evolution of WSS over time; however, time--averaged indicators can be precisely identified also under the rigid wall assumption~\cite{eslami2020effect}.
On the left, instead, we report the space--averaged WSS over time, computed over two conveniently located portions of the vessel wall. While the results attained with SRB--TFO are closer to the high--fidelity target, particularly on \emph{Region A}, ST--GRB is nonetheless able to nicely capture the macroscopic trends.   

\section{Conclusions} \label{sec:conclusions}
In this work, we discussed the application of the ST--GRB method for the efficient numerical simulation of h{\ae}modynamics in compliant vessels. We considered three--dimensional vascular anatomies (either idealized or realistic) and exploited linear projections onto \emph{ad hoc} low--dimensional spatio--temporal subspaces to effectively solve the unsteady incompressible Navier--Stokes equations. Relying on the coupled momentum model, the latter have been endowed with non--standard Robin boundary conditions, that seamlessly allow to take into account the effect of wall compliance on the flow. Therefore, compared to classical CFD solvers, the proposed model is characterized by reduction processes that happen both at the physical and at the numerical level. To the best of our knowledge, this work represents the first application of ST--RB methods to reduced FSI models, particularly in the h{\ae}modynamics context. The numerical results empirically demonstrated the ability of our model to quickly produce rather precise solution approximations. Nonetheless, we acknowledge that the classical RB method should be preferred over ST--GRB in high--precision applications or if the temporal dynamics arising within the high--fidelity solutions manifold are too complex. Indeed, the computational gains of ST--GRB over SRB--TFO vanish if the number of temporal modes is too large, compared to the number of discrete time instants. Furthermore, SRB--TFO appears to offer better data efficiency and scalability property compared to ST--GRB, and  improved robustness in extrapolation regimes strongly linked to temporal dynamics. \\

Several further developments could be pursued to address the current limitations of the ST--GRB method. Computational gains at improved precision levels could be attained by associating each spatial mode to a tailored reduced basis in time, as proposed and demonstrated in~\cite{choi2019space}. 
To enhance data efficiency and scalability, other than generalization capabilities, one promising avenue is to relax the reliance on a fixed global reduced basis to capture the entire solution manifold. To this aim, we foresee the extension to space--time--reduced scenarios of localized variants of the RB method, either based on clustering with respect to the parametric domain~\cite{eftang2010hp, haasdonk2011training, maday2013locally} or on adaptive enrichments through \emph{a posteriori} error indicators~\cite{drohmann2011adaptive, ohlberger2015error}. In cases where online intrusivity does not pose major implementation challenges, improved computational performances could be achieved by adopting alternative hyper--reduction techniques for the nonlinear convective term, such as DEIM~\cite{chaturantabut2010nonlinear} or S--OPT~\cite{lauzon2024sopt}. Lastly, efficiently integrating more physiological zero--dimensional models of the downstream vasculature~\cite{vignon2006outflow, kim2010patient}, as well as accounting for inter--patient and intra--patient geometrical variability~\cite{guibert2014group, ballarin2016fast, pegolotti2021model}, are key steps toward bridging the gap with clinical applications.

\section*{Acknowledgements}
This work was supported by the Swiss National Science Foundation (grant agreement No 200021\_197021, ``Data--driven approximation
of hemodynamics by combined reduced order modeling and deep neural networks'').

\section*{Data availability}
Our implementation of ST--RB models for h{\ae}modynamics is publicly available at \url{https://github.com/riccardotenderini/STRB4hemo}, under \emph{BSD 3-Clause license}. High--fidelity simulation data can be made available upon request.

\bibliographystyle{unsrtnat}
\bibliography{references}

\newpage

\appendix

\section{ST--GRB problem with inhomogeneous initial conditions} \label{app: incorporating inomogeneous initial conditions}

In this appendix section, we extend the ST--RB problem formulation presented in Section~\ref{subsec: ST-RB problem definition} by considering non--zero initial conditions. Although none of the test cases presented in the manuscript falls within this scenario, a suitable model derivation and a convenient implementation pipeline become crucial when performing sequential ST--RB simulations over multiple cycles, as discussed in Section~\ref{subsec: membrane test 1}.\\

So, we introduce the vectors $\bm{u_0} \in \R^{N_u^s}$, $\bm{d_0} \in \R^{N_d^s}$, $\bm{p_0} \in \R^{N_p^s}$, and $\bm{\lambda_0} \in \R^{N_\lambda}$, that contain the initial conditions for the velocity, displacement, pressure and Lagrange multipliers fields. Furthermore, since we employ a second--order multistep time marching scheme, we suppose to know the velocity fields $\{\bm{u_{-s}}\}_{s=1}^{S-1} \in \R^{N_u^s}$ at time $t = -\Delta t$. Incidentally, we acknowledge that the initial data may be parameter dependent, but we intentionally omit this aspect in the following formulation for simplicity. 
The initial conditions are then used as temporal lifting terms, so to approximate the spatio--temporal solution fields as
\begin{equation*}
\label{eq: space-time reduced solutions with IC}
\begin{alignedat}{2}
&\bm{u}^{st} \approx \ \bm{u_0} \otimes \bm{1_t} \ + \ \sum_{\ell=1}^{n_u^{st}} \widehat{\bm{u}}_\ell \bm{\pi}_\ell^u \ = \ \bm{u_0}^{st} + \bm{\Pi}^u \widehat{\bm{u}}; \qquad 
&&\bm{d}^{st} \approx \ \bm{d_0} \otimes \bm{1_t} \ + \ \bm{u_0} \otimes \bm{t_t} \ + \ \sum_{\ell=1}^{n_u^{st}} \widehat{\bm{u}}_\ell \bm{\pi}_\ell^d \ = \ \bm{d_0}^{st} + \tilde{\bm{d_0}}^{st} + \bm{\Pi}^d \widehat{\bm{u}}~;
\\
&\bm{p}^{st} \approx \ \bm{p_0} \otimes \bm{1_t} \ + \ \sum_{k=1}^{n_p^{st}} \widehat{\bm{p}}_k \bm{\pi}_k^p \ = \ \bm{p_0}^{st} + \bm{\Pi}^p \widehat{\bm{p}}~; 
\qquad 
&&\bm{\lambda}^{st} \approx \ \bm{\lambda_0} \otimes \bm{1_t} \ + \ \sum_{i=1}^{n_\lambda^{st}} \widehat{\bm{\lambda}}_i \bm{\pi}_i^{\lambda} \ = \ \bm{\lambda_0}^{st} + \bm{\Pi}^\lambda \widehat{\bm{\lambda}}~;
\end{alignedat}
\end{equation*}
where $\bm{1_t} \in \R^{N^t}$ is a vector of ones of size $N^t$ and $\bm{t_t} := \mathcal{P}_0 \left(\mathcal{E}_S^t\left(\bm{1_t}\right)\right) \in \R^{N^t}$ is its discrete primitive, rescaled so that it nullifies at $t=0$. We remark that the displacement field features an ``extra'' terms that depends on $\bm{u_0}$, which allows to exactly guarantee the kinematic coupling condition with the velocity field. It is fundamental to note that this approach allows to strongly enforce the initial conditions if all temporal modes are equal to zero at $t=0$. This condition can be fulfilled by subtracting the initial data from the snapshots prior to the bases construction process. \\

Upon the lifting procedure, the application of the ST--GRB method requires to find $\widehat{\bm{w}} = [\widehat{\bm{u}}, \widehat{\bm{p}}, \widehat{\bm{\lambda}}]\in \R^{n^{st}}$ such that
\begin{equation}
\label{eq: ST-RB definition with IC}
\begin{alignedat}{2}
\widehat{\bm{r}}(\widehat{\bm{w}}) := \ \bm{\Pi}^T \Bigg(
&\bm{F}^{st}(\bm{\mu_f}) + \bm{F}^{st}_0(\bm{\tilde{\mu}_m}, \bm{u_0}, \bm{u_{-1}}, \ldots, \bm{u_{-S+1}}) + \bm{F}^{st}_{0,L}(\bm{\tilde{\mu}_m}, \bm{u_0}, \bm{d_0}, \bm{p_0}, \bm{\lambda_0}) 
\\
&- \left( \bm{A}^{st} + \bm{A}^{\Gamma, st}(\bm{\tilde{\mu}_m})\right) \bm{\Pi}\widehat{\bm{w}} - \mathcal{E}_u \Big(\bm{c}^{st}\left( \bm{u_0}^{st} + \bm{\Pi}^u\widehat{\bm{u}}\right) + \left(\bm{A_s}^{st} (\bm{\tilde{\mu}_m}) + \bm{E_s}^{st}\right) \bm{\Pi}^d \widehat{\bm{u}} \Big) \Bigg) = \bm{0}~,
\end{alignedat}
\end{equation}
where $\bm{F}^{st}_0$ takes into account the contributions of the velocity initial conditions --- conveniently denoted as $\bm{U}_0$ in the following --- stemming from the time marching scheme, while $\bm{F}^{st}_{0,L}$ models the presence of the temporal lifting terms. \\

The vector $\bm{F}^{st}_0 \in \R^{N_{st}}$ only features a non--zero velocity block $\bm{F}^{st}_{0,u} \in \R^{N_u^s N^t}$, which takes the following form: 
\begin{equation*}
\bm{F}^{st}_{0,u}(\bm{\tilde{\mu}_m}, \bm{U}_0) \ := \ \bm{F}_{0,u}^{M,st}(\bm{U}_0) \ + \ \bm{F}_{0,u}^{\Gamma,st}(\bm{\tilde{\mu}_m}, \bm{U}_0)~,
\end{equation*}
where
\begin{align*}
\bm{F}_{0,u}^{M,st}(\bm{U}_0) &= \Bigg[ 
\left(\sum_{s=1}^S \alpha_s \bm{M} \bm{u}_{1-s}\right)^T, \ \cdots, \left(\sum_{s=S}^S \alpha_s \bm{M} \bm{u}_{S-s}\right)^T, \ \underbrace{\left(\bm{0}_{N_u^s}\right)^T, \ \cdots, \ \left(\bm{0}_{N_u^s}\right)^T}_{N^t - \ S} \Bigg]^T~; \\ 
\bm{F}_{0,u}^{\Gamma,st}(\bm{\tilde{\mu}_m}, \bm{U}_0) &= \Bigg[ 
\left(\sum_{s=1}^S \alpha_s (\bm{\tilde{\mu}_m})_1 \ \bm{M_s}\right) \bm{u}_{1-s}^T, \ \cdots, \ \left(\sum_{s=S}^S \alpha_s (\bm{\tilde{\mu}_m})_1 \ \bm{M_s} \bm{u}_{S-s}\right)^T, \ \underbrace{\left(\bm{0}_{N_u^s}\right)^T, \ \cdots, \left(\bm{0}_{N_u^s}\right)^T}_{N^t - \ S} \Bigg]^T~.
\end{align*}
Within the ST--RB formulation, these terms translate into two non--zero velocity right--hand side blocks of the form
\begin{align*}
\widehat{\bm{F}}_{0,u}^M\bigl(\bm{U_0}\bigr) &= \left(\bm{\Pi}^u\right)^T \bm{F}_1^{st}\bigl(\bm{U_0}\bigr) \in \R^{n_u^{st}} \ &&\text{where} \ \left(\widehat{\bm{F}}_{0,u}^M\right)_\ell = \sum_{s=1}^S \sum_{s'=s}^S \alpha_{s'} \left(\bm{\widebar{M}}  \bm{\widebar{u}}_{s-s'}\right)_{\ell_s} \left(\bm{\psi}^u_{\ell_t}\right)_s~; \\
\widehat{\bm{F}}_{0,u}^{\Gamma}\bigl(\bm{\tilde{\mu}_m}, \bm{U_0}\bigr) &= \left(\bm{\Pi}^u\right)^T \bm{F}_1^{\Gamma, st}\bigl(\bm{\tilde{\mu}_m}, \bm{U_0}\bigr) \in \R^{n_u^{st}} \ &&\text{where} \ \left(\widehat{\bm{F}}_{0,u}^{\Gamma}\big(\bm{\tilde{\mu}_m}\big)\right)_\ell = (\bm{\tilde{\mu}_m})_1 \ \sum_{s=1}^S \sum_{s'=s}^S \alpha_{s'} \left( \bm{\widebar{M}_s}  \bm{\widebar{u}}_{s-s'}\right)_{\ell_s} \left(\bm{\psi}^u_{\ell_t}\right)_s~;
\end{align*}
where $\widebar{\bm{u}}_s := (\bm{\Phi}^u)^T \bm{u_s} \in \R^{n_u^s}$, for $s = -S+1, \ldots, 0$. Analogously, we also define $\widebar{\bm{p}}_0 := (\bm{\Phi}^p)^T \bm{p_0} \in \R^{n_p^s}$ and $\widebar{\bm{\lambda}}_0 := (\bm{\Phi}^\lambda)^T \bm{\lambda_0} \in \R^{n_\lambda^s}$. Hence, the space--time--reduced counterpart of $\bm{F}^{st}_0$ is $\widehat{\bm{F}}^{st}_0 := \bm{\Pi}^T \bm{F}^{st}_0 = \widehat{\mathcal{E}}_u(\widehat{\bm{F}}_{0,u}^M + \widehat{\bm{F}}_{0,u}^\Gamma)$, where $\mathcal{\widehat{E}}_u: \R^{n_u^{st}} \to \R^{n^{st}}$ is space--time--reduced velocity zero--padding operator, introduced in Eq.\eqref{eq: space-time-reduced residual}.\\

The vector $\bm{F}^{st}_{0,L} \in \R^{N_{st}}$, instead, writes as follows:
\begin{equation*}
\bm{F}^{st}_{0,L}(\bm{\tilde{\mu}_m}, \bm{u_0}, \bm{d_0}, \bm{p_0}, \bm{\lambda_0}) = - \left( \bm{A}^{st} + \bm{A}^{\Gamma, st}(\bm{\tilde{\mu}_m})\right) \bm{w_0}^{st} - \mathcal{E}_u \Big(\left(\bm{A_s}^{st} (\bm{\tilde{\mu}_m}) + \bm{E_s}^{st}\right) \bm{d_0}^{st} \Big)~,
\end{equation*}
where $\bm{w_0}^{st} := [\bm{u_0}^{st}, \bm{p_0}^{st}, \bm{\lambda_0}^{st}] \in \R^{N_{st}}$ is the global spatio--temporal vector of initial conditions. The assembly of its space--time--reduced counterpart $\widehat{\bm{F}}^{st}_{0,L} := \bm{\Pi}^T \bm{F}^{st}_{0,L}$, that appears within the ST--GRB formulation, can be efficiently performed exploiting the space--time factorization of the reduced bases. For instance, the $\ell$--th component of the velocity block stemming from the term $\bm{\Pi}^T \bm{A}^{st} \bm{w_0}^{st}$ is such that:
\begin{equation*}
\left(\bm{\Pi}^T \bm{A}^{st} \bm{w_0}^{st}\right)_{u, \ell} = \left(\widebar{\bm{M}} \widebar{\bm{u}}_0\right)_{\ell_s} \left(\bm{1_t}^T \bm{\psi}_{\ell_t}^u + \sum_{s=1}^S \alpha_s \left(\bm{1_t}^T\right)_{:-s} \left(\bm{\psi}_{\ell_t}^u\right)_{s:}\right) 
+ \beta \Delta t 
\left(\left(\bigl(\widebar{\bm{A}} + \widebar{\bm{R}}\bigr) \widebar{\bm{u}}_0 + \widebar{\bm{B}} \widebar{\bm{p}}_0 + \widebar{\bm{C}} \widebar{\bm{\lambda}}_0\right)_{\ell_s} \left(\bm{1_t}^T \bm{\psi}_{\ell_t}^u\right) \right)~,
\end{equation*}
where we employ the usual index notation $\ell = \mathcal{F}_u(\ell_s, \ell_t)$. Incidentally, we note that all terms deriving from the coupled momentum model are affinely parametrized.\\

Lastly, we handle the presence of the velocity initial condition for the evaluation of convective term. Exploiting the affine components introduced in Section~\ref{subsubsec: ST-reduced convective term}, the space--time--reduced convective term entries become
\begin{equation*}
\label{eq: convective term affine decomposition ST-ROM with IC}
\begin{alignedat}{2}
\left(\widehat{\bm{c}}(\widehat{\bm{u}}; \widebar{\bm{u}}_0)\right)_m \approx
\beta \Delta t \sum_{\ell_s'=1}^{n_c} \sum_{\ell_s''=1}^{n_c} \Biggl( 
&\left(\bm{\widebar{k}}_{\ell_s'\ell_s''}\right)_{m_s} \sum_{\ell_t'=1}^{n_u^t} \sum_{\ell_t''=1}^{n_u^t} \widehat{\bm{u}}_{\ell'} \widehat{\bm{u}}_{\ell''} \ \psi^{u,3}_{\ell_t',\ell_t'',m_t} \\
& + \left(\widebar{\bm{u}}_0\right)_{\ell_s'} \left(\bm{\widebar{k}}_{\ell_s'\ell_s''} + \bm{\widebar{k}}_{\ell_s''\ell_s'}\right)_{m_s} \sum_{\ell_t''=1}^{n_u^t} \widehat{\bm{u}}_{\ell''} \left({\bm{\psi}^u_{\ell_t''}}^T \bm{\psi}^u_{m_t}\right)  \\
&  + \left(\widebar{\bm{u}}_0\right)_{\ell_s'} \left(\widebar{\bm{u}}_0\right)_{\ell_s''} \left(\bm{\widebar{k}}_{\ell_s'\ell_s''}\right)_{m_s} \left(\bm{1_t}^T \bm{\psi}^u_{m_t}\right)~,~
\Biggr)
\end{alignedat}
\end{equation*}
where the quantity $\psi^{u,3}_{\ell_t',\ell_t'',m_t}$ is defined as in Eq.\eqref{eq: convective term affine decomposition ST-ROM}. We note that the last two contributions are, respectively, linear and constant with respect to the unknown $\widehat{\bm{u}}$. Similarly, the space--time--reduced convective Jacobian entries can be conveniently approximated as follows:
\begin{equation*}
\label{eq: convective Jacobian affine decomposition ST-ROM with IC}
\left(\widehat{\bm{J}}_{\bm{c}}(\widehat{\bm{u}}; \widebar{\bm{u}}_0)\right)_{m'm''} \approx
\beta \Delta t \ \sum_{\ell_s=1}^{n_{c,J}} \left( \left(\bm{\widebar{K}}_{\ell_s}\right)_{m_s' m_s''} \sum_{\ell_t=1}^{n_u^t} \widehat{\bm{u}}_\ell \ \psi^{u,3}_{\ell_t, m_t', m_t''} 
+ \left(\widebar{\bm{u}}_0\right)_{\ell_s} \left(\bm{\widebar{K}}_{\ell_s}\right)_{m_s' m_s''} \left({\bm{\psi}^u_{m_t'}}^T \bm{\psi}^u_{m_t''}\right)
\right)~.
\end{equation*}

Ultimately, the incorporation of inhomogeneous initial conditions modifies the expression of the space--time--reduced residual and of its Jacobian as follows:
\begin{align*}
\widehat{\bm{r}}(\widehat{\bm{w}}) &= 
\widehat{\bm{F}}(\bm{\mu_f}) + \widehat{\bm{F}}_0^{st}(\bm{\tilde{\mu}_m}) + \widehat{\bm{F}}_{0,L}^{st}(\bm{\tilde{\mu}_m})
- \left(\widehat{\bm{A}} + \widehat{\bm{A}}^{\Gamma}(\bm{\tilde{\mu}_m})\right) \widehat{\bm{w}} 
- \widehat{\mathcal{E}}_u \big(\widehat{\bm{c}}(\widehat{\bm{u}}; \widebar{\bm{u}}_0)\big) - \widehat{\mathcal{E}}_u \big( (\widehat{\bm{A}}_{\bm{s}} (\bm{\tilde{\mu}_m}) + \widehat{\bm{E}}_{\bm{s}}) \ \widehat{\bm{u}} \big)~;\\
\widehat{\bm{J}}_{\widehat{\bm{r}}}\big(\widehat{\bm{w}}\big) &=
- \left(\widehat{\bm{A}} + \widehat{\bm{A}}^{\Gamma}(\bm{\tilde{\mu}_m})\right) - \widehat{\mathcal{E}}_{u,u}\big(\widehat{\bm{J}}_{\bm{c}}(\widehat{\bm{u}}; \widebar{\bm{u}}_0)\big) - \widehat{\mathcal{E}}_{u,u}\big(\widehat{\bm{A}}_{\bm{s}} (\bm{\tilde{\mu}_m}) + \widehat{\bm{E}}_{\bm{s}}\big)~.
\end{align*}

\end{document}